\documentclass[a4paper,11pt]{article}
\pdfoutput=1
\usepackage{amsmath}
\usepackage{amssymb}
\usepackage{amsfonts}
\usepackage[utf8]{inputenc}
\usepackage{url}
\usepackage[final]{graphicx, epsfig}
\usepackage{amsthm}
\usepackage{pstricks, tikz}
\usepackage{textgreek}
\usepackage{hyperref}
\usepackage{authblk}

\theoremstyle{definition}
\newtheorem{Algorithm}{Algorithm}
\numberwithin{Algorithm}{section}

\newtheorem{Example}[Algorithm]{Example}
\newtheorem{Test problem}[Algorithm]{Test problem}
\newtheorem{Test Configurations}[Algorithm]{Test configurations}

\newtheorem{remark}[Algorithm]{Remark}

\newtheorem{Building block}[Algorithm]{Building block}

\title{Detecting and approximating decision boundaries in low dimensional spaces}
\author[1,3]{M. Grajewski\thanks{corresponding author. Email: grajewski@fh-aachen.de}}
\author[2,1]{A. Kleefeld}

\affil[1]{FH Aachen University of Applied Sciences, Faculty of Medical Engineering and Technomathematics, Heinrich-Mußmann-Str. 1, 52428 Jülich, Germany}
\affil[2]{Forschungszentrum Jülich GmbH, Jülich Supercomputing Centre, Wilhelm-Johnen-Str., 52425 Jülich, Germany}
\affil[3]{Institute for Data-Driven Technologies, FH Aachen University of Applied Sciences, Heinrich-Mußmann-Str. 1, 52428 Jülich, Germany}

\begin{document}
    \maketitle
    \begin{abstract}
        A method for detecting and approximating fault lines or surfaces, respectively, or decision curves in two and three
        dimensions with guaranteed accuracy is presented.
        Reformulated as a classification problem, our method starts from a set of scattered points along with the
        corresponding classification algorithm to construct a representation of a decision curve by points with
        prescribed maximal distance to the true decision curve.
        Hereby, our algorithm ensures that the representing point set covers the decision curve in its entire extent
        and features local refinement based on the geometric properties of the decision curve.
        We demonstrate applications of our method to problems related to the detection of faults, to Multi-Criteria
        Decision Aid and, in combination with Kirsch's factorization method, to solving an inverse acoustic
        scattering problem.
        In all applications we considered in this work, our method requires significantly less pointwise classifications
        than previously employed algorithms.
    \end{abstract}

    {\bf keywords}: fault detection, fault approximation, inverse scattering problem, MCDA

    \section{Introduction}\label{sec:introduction}

    Let us consider a piecewise constant function $f: \mathbb{R}^m \supset \Omega \to \{1,2, \hdots, n\}$ with $\Omega$ being compact,
    simply connected and equipped with a piecewise smooth boundary, $m \in \{2, 3\}$.
    Such $f$ subdivides $\Omega$ into mutually disjoint subsets $\Omega_i := f^{-1}(i)$ with $\Omega = \cup_{i=1}^n \overline{\Omega_i}$ and
    $\Omega_i \cap \Omega_j = \emptyset$, if $i \neq j$.
    We assume that each $\Omega_i$ features a piecewise smooth boundary.
    We are interested in approximating $\Gamma_{i,j} := \overline{\Omega_i} \cap \overline{\Omega_j}$ relying on as few evaluations
    of $f$ as possible and present in this work an algorithm for this task.
    We choose this quantity as a measure of efficiency because evaluating $f$ can be arbitrarily costly in
    applications and dominates the runtime in such a case.
    Our problem can be immediately understood as a classification problem, such that we identify
    $\Omega_i$ with a class $i$ and interpret the curves or surfaces of discontinuity $\Gamma_{i,j}$ as decision curves or surfaces.

    One field of application is economics and operations research.
    Multicriteria Decision Aid (MCDA) methods can help a decision maker to choose the best one from a finite number of
    alternatives based on different, even conflicting criteria.
    MCDA methods assume that a decision depends on quantifiable parameters $(x_1, \hdots, x_m)$ (``input factors'') and
    is drawn deterministically.
    For an overview over various MCDA approaches, applications and case studies, we refer
    to~\cite{Figueira.2005,PAPATHANASIOU.2019} among many others and the references cited therein.
    In this context, $\Omega_i$ is the set of all input factors that lead to the decision for the $i$-th alternative in the
    MCDA method.
    Analysing the decision process with respect to the input factors means consistently describing
    all $\Omega_i$ based upon evaluating the MCDA method for arbitrary combinations of input factors.
    This can be achieved by approximating all $\Gamma_{i,j}$.\\
    Computing a reconstruction of an obstacle in three dimensions is a field of application in acoustic
    scattering theory.
    More precisely, one wants to determine the support of an inhomogeneous object (its boundary to be exact) from
    measured far-field data which typically is a desired task in non-destructive testing.
    The far-field data are obtained for different incident waves and measured points on the unit sphere.
    Several reconstruction algorithms to find the boundary of the unknown inhomogeneity are available such as
    iterative methods~\cite{kleefeldlin3D}, decomposition methods~\cite{jiguangsun} and sampling/probe methods
    (see~\cite{Potthast2006} for a detailed overview).
    The latter ones can be further categorised into the linear sampling method, the generalized linear sampling method,
    the factorization method, the probe method and variants of it (refer
    to~\cite{Colto1996,audibert,Kirsch2008,Ikeha2000a,Potthast2006}, respectively).
    However, here we will focus on the classical factorization method for the acoustic transmission problem, refer also
    to~\cite{anachakle}, with which one can decide if a given point is located inside or outside the obstacle.
    Therefore, the factorization method transfers the reconstruction of an obstacle to a classification problem and thus
    into a field of application of our method.
    Note that the far-field data within~\cite{anachakle} has also been used in~\cite{bazan}
    and~\cite{Harris2019}.

    Finding and approximating the sets $\Gamma_{i,j}$, sometimes called fault lines, is important in the exploration of
    natural resources.
    The presence and the location of $\Gamma_{i,j}$ can provide useful insights for exploring and later on exploiting
    of oil reservoirs~\cite{Gutzmer.1997} and play a significant role in some geophysical applications~\cite{Gout2008}.
    The underlying mathematical problem is closely related to ours, albeit not the same, as usually, the function $f$
    considered does not provide integer values as in our case.
    Therefore, an additional algorithm for detecting fault lines is required then, and the classification of a single
    point may be not trivial anymore as it is in our case.
    Moreover, algorithms for fault detection may need to deal with noisy data (e.g.~\cite{Bozzini.2013}),
    whereas we consider certain data only.

    Many algorithms for detecting and approximating the sets $\Gamma_{i,j}$ have been proposed,
    like~\cite{Arge.1994,Gutzmer.1997, Gout2008, Allasia2010EfficientAA, Bozzini.2013} among many others, which all feature
    strengths and weaknesses.
    However, the vast majority of these approaches restrict to the 2D case, whereas we present a method for 2D and 3D.
    The algorithm proposed in~\cite{Allasia2010EfficientAA} and the work cited therein was the starting point for our
    research.
    Classification is one of the standard problems in Machine Learning.
    There are a lot of powerful and versatile algorithms available which could readily applied to our problem;
    we refer to~\cite{Bishop2006} for an overview.
    These methods are however designed for uncertain data in high dimensional spaces, whereas we consider secure data
    in low dimensional spaces, a completely different use case.

    Our method approximately describes the $\Gamma_{i,j}$ by providing a set of points with a guaranteed maximal normal
    distance to $\Gamma_{i,j}$.
    These points are intended for constructing a polygon (2D) or a surface triangulation in 3D.
    While there are more sophisticated and elegant ways of describing these sets, it allows us to (approximately)
    replace an actual classification by a simple and fast point-in-polygon test.

    This article is organised as follows: We describe our algorithm for 2D and 3D in Section~\ref{sec:algorithm} and
    elaborate on the direct and inverse acoustic scattering problem, one of our applications, in
    Section~\ref{sec:direct_acoustic_scattering}.
    We present results and applications of our algorithm to the detection of faults, to decision modeling and to inverse
    scattering in Section~\ref{sec:applications} and finally conclude in Section~\ref{sec:conclusions}.

    \section{The algorithm}\label{sec:algorithm}
    \subsection{Detection and Approximation of $\Gamma_{i,j}$ in 2D}\label{subsec:2d_description}
\begin{figure}
        \begin{center}
            \includegraphics[height=0.7cm]{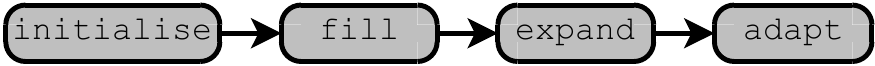}\hspace{1cm}
            \caption{General flow chart.}\label{fig_flow_chart}
        \end{center}
    \end{figure}

\begin{figure}
        \begin{center}
            \includegraphics[height=4.5cm]{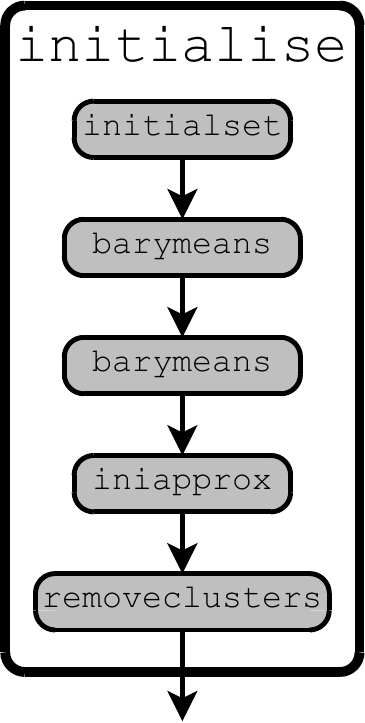}\hspace{0.7cm}
            \includegraphics[height=7.5cm]{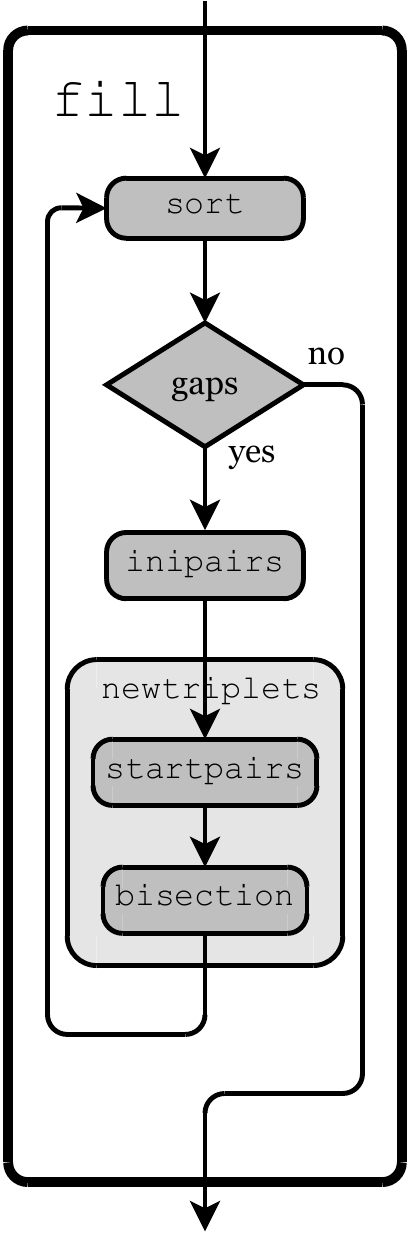}\hspace{0.7cm}
            \includegraphics[height=7.5cm]{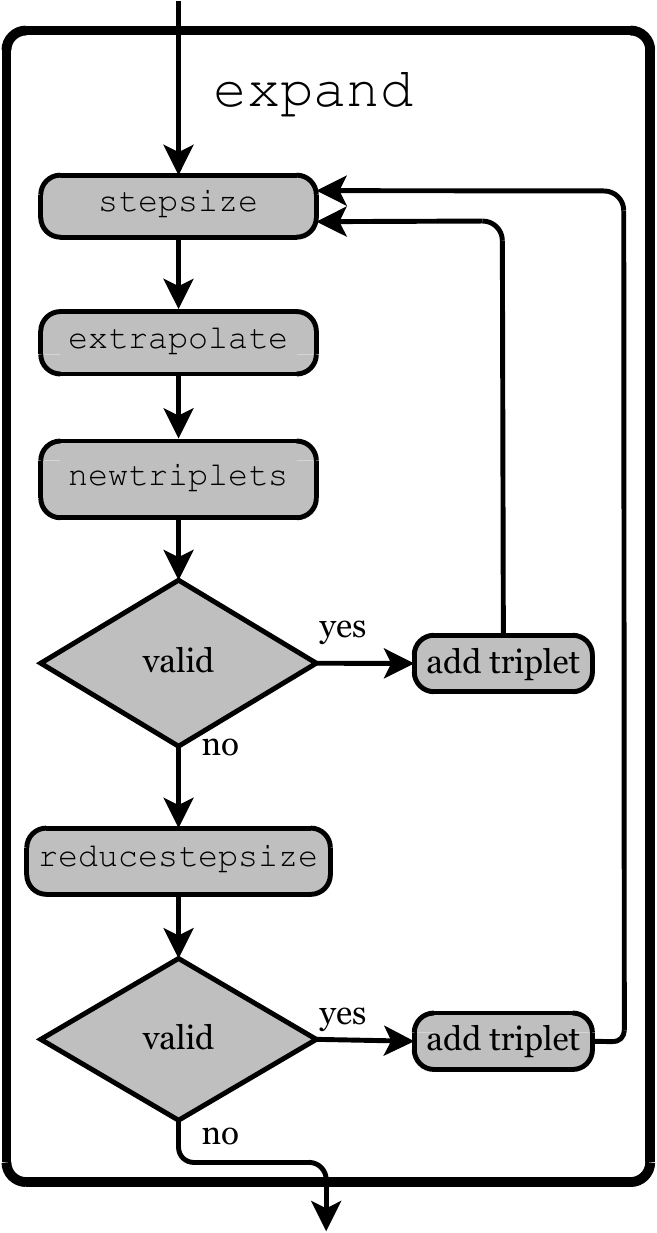}\hspace{0.7cm}
            \includegraphics[height=7.5cm]{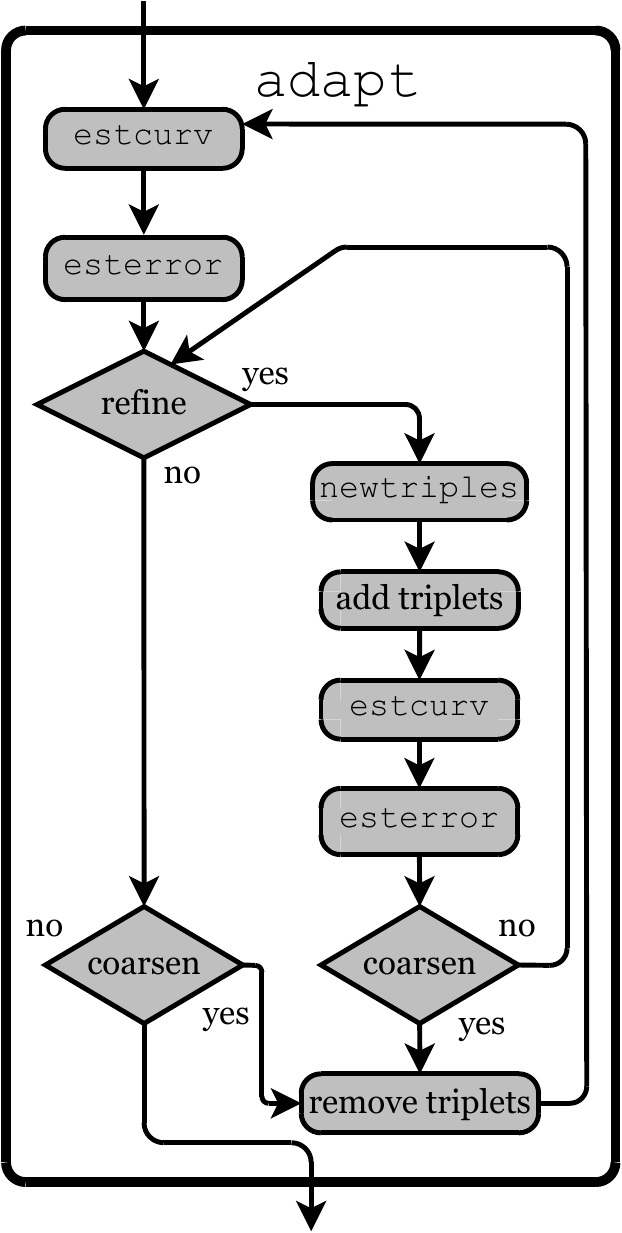}
            \caption{Flow charts of the Algorithms \texttt{initialise}, \texttt{fill}, \texttt{expand} and
            \texttt{adapt}.
            Note that Building block \texttt{barymeans} inside \texttt{initialise} is employed twice, albeit to
            different sets.}\label{fig_bb1_2D}
        \end{center}
    \end{figure}

    In this section, we present our algorithm for approximating $\Gamma_{i,j} \neq \emptyset$ for $m=2$ by
    sufficiently many well distributed and ordered points sufficiently close $\Gamma_{i,j}$.
    This implicitly constitutes a polygonal description of $\Gamma_{i,j}$.
    For convenience, we provide flow charts (Figs.\ \ref{fig_flow_chart} and~\ref{fig_bb1_2D}) and describe selected
    building blocks in detail.
    In the flow charts, actual building blocks are typed in monospace lettering.

    \begin{Algorithm}[\texttt{initialise}]
        \label{alg:initialise}
        This algorithm aims at providing initial approximations to any $\Gamma_{i,j}$; we refer to
        Fig.\ \ref{fig_flow_chart} (right) for an overview.
        From now on, we assume $\Gamma_{i,j} \neq \emptyset$.
        In Building block \texttt{initialset}, we sample $f$ on $\Omega$ rather coarsely and obtain an
        initial point set $X$ along with the corresponding classification information.
        Building block \texttt{barymeans} creates additional sampling points in the vicinity of any $\Gamma_{i,j}$.
        Following Allasia et al.\ \cite{Allasia2010EfficientAA}, we employ a $k_{\mathrm{near}}$-nearest neighbour
        approach: For any $x \in X$, let $N(x)$ be the set of the $k_{\mathrm{near}}$-
        nearest neighbours of $x$ in $X$.
        If $N(x) \cap \Omega_i \neq \emptyset$ and $N(x) \cap \Omega_j \neq \emptyset$, we consider $x$ close to $\Gamma_{i,j}$.
        Let $N_{\ell}(x) = N(x) \cap \Omega_{\ell}$ for some $\ell$.
        Hence, $N(x) = \bigcup_{i=1}^r N_{c_i}(x)$ for certain indices $c_1, \hdots, c_r$, $N_{c_i}(x) \neq \emptyset$.
        We compute the barycentres $b_{c_i}$ of all $N_{c_i}(x)$, $1 \leq i \leq r$, and then their arithmetic means
        $y_{c_i,c_j} = 0.5(b_{c_i} + b_{c_j}), 1 \leq i < j \leq r$.
        Let $M(x)$ be the set of all $y_{c_i,c_j}$ generated from $N(x)$.
        If $N(x) \subset \Omega_{\ell}$ for some $\ell$, we set $M(x) = \emptyset$.
        We end up with $\mathcal{M} = \bigcup_{x \in X} M(x)$.
        By definition of $\mathcal{M}$, there are no duplicate points; however, a practical implementation requires removing duplicates.
        After classifying the points in $\mathcal{M}$, we repeat
        \texttt{barymeans} on $\mathcal{M}$ obtaining sets $M^2(x)$ for any $x \in \mathcal{M}$, then
        $\mathcal{M}^2 = \bigcup_{x \in \mathcal{M}} M^2(x)$ and ultimately
        a further enriched set of sampling points $\overline{X} = X \cup \mathcal{M} \cup \mathcal{M}^2$.\\
        Building block \texttt{iniapprox} computes initial approximations for all $\Gamma_{i,j}$.
        For any $x \in (\mathcal{M} \cup \mathcal{M}^2) \cap \Omega_i$, we search the nearest point $x' \in \overline{X} \cap \Omega_j, j > i$.
        We use $x$ and $x'$ as starting points for a bisection algorithm on the line $\overline{x'x}$.
        If the bisection algorithm is successful, we end up with a point pair $x^{(i)} \in\Omega_i$ and $x^{(j)} \in \Omega_j$ with
        $\| x^{(i)} - x^{(j)}\|\leq 2\varepsilon_b$, where $\varepsilon_b$ is a user-prescribed threshold.
        Then, the distance of $x^{(i,j)} = 0.5(x^{(i)} + x^{(j)})$ to $\Gamma_{i,j}$ is at most $\varepsilon_b$.
        We subsume the bisection process up to $\varepsilon_b$ and computing $x^{(i,j)}$ from its results in Building block
        \texttt{bisection}.
        From now on, we consider point triplets for approximating
        $\Gamma_{i,j}$ only; for any such triplet $x$, the superscript $(i)$ denotes the point in $\Omega_i$, the superscript
        $(j)$ its counterpart in $\Omega_j$ and the superscript $(i,j)$ the arithmetic mean of the two points.
        We end up with a set of triplets $\widetilde{S}_{i,j}$.
        We moreover set $\widetilde{S}_{i,j}^{(i)} = \{ x^{(i)} \:|\: x \in \widetilde{S}_{i,j} \}$ and
        $\widetilde{S}_{i,j}^{(j)} = \{ x^{(j)} \:|\: x \in \widetilde{S}_{i,j} \}$
        It may occur that some triplets in $\widetilde{S}_{i,j}$ are tightly clustered.
        We thin such clusters as they add to complexity but not to accuracy by removing appropriate triplets
        (Fig.~\ref{Fig_ClusterRemoved}).
        After cluster removal, \texttt{initialise} provides sets of triplets $S_{i,j}$.
    \end{Algorithm}
    \begin{remark}
        Building block \texttt{initialise} detects $\Gamma_{i,j} \neq \emptyset$ by $\widetilde{S}_{i,j} \neq \emptyset$.
        However, depending on $X$, it may happen that $\widetilde{S}_{i,j} = \emptyset$ even if $\Gamma_{i,j} \neq \emptyset$.
        Reliably detecting all non-empty $\Gamma_{i,j}$ depends on a sufficiently large $X$ and thus ultimately on the user.
    \end{remark}

\begin{figure}
        \begin{center}
            \includegraphics[height=5cm]{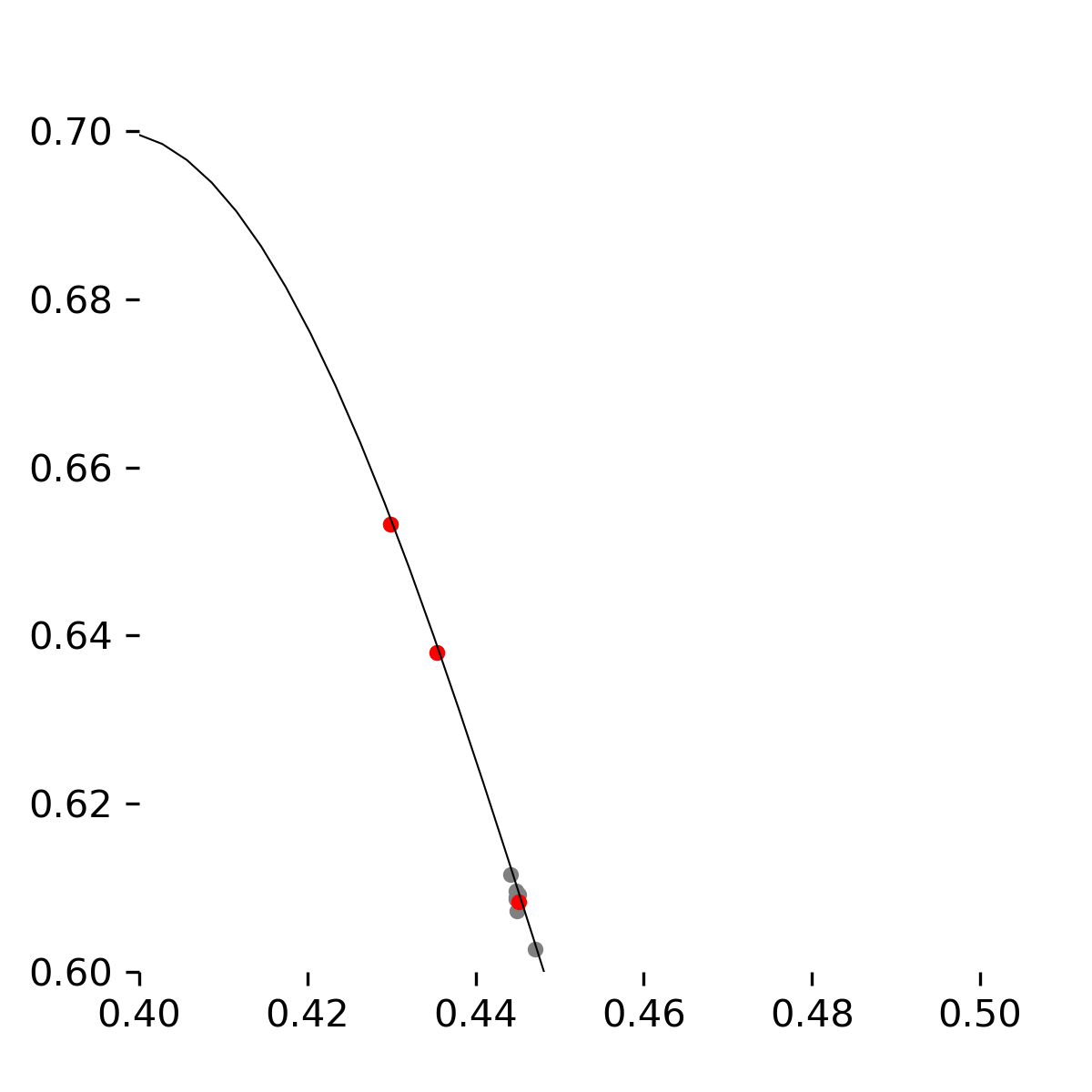}
            \caption{Cluster removal by \texttt{removeclusters} in Example~\ref{example1} for $\widetilde{S}_{1,2}$ (greyed out); points in $S_{1,2}^{(1)}$ are displayed
            in red. This figure is an excerpt of Fig.\ \ref{Fig_Subdivision1}, right.} \label{Fig_ClusterRemoved}
        \end{center}
    \end{figure}

\begin{figure}
        \begin{center}
            \includegraphics[height=6cm]{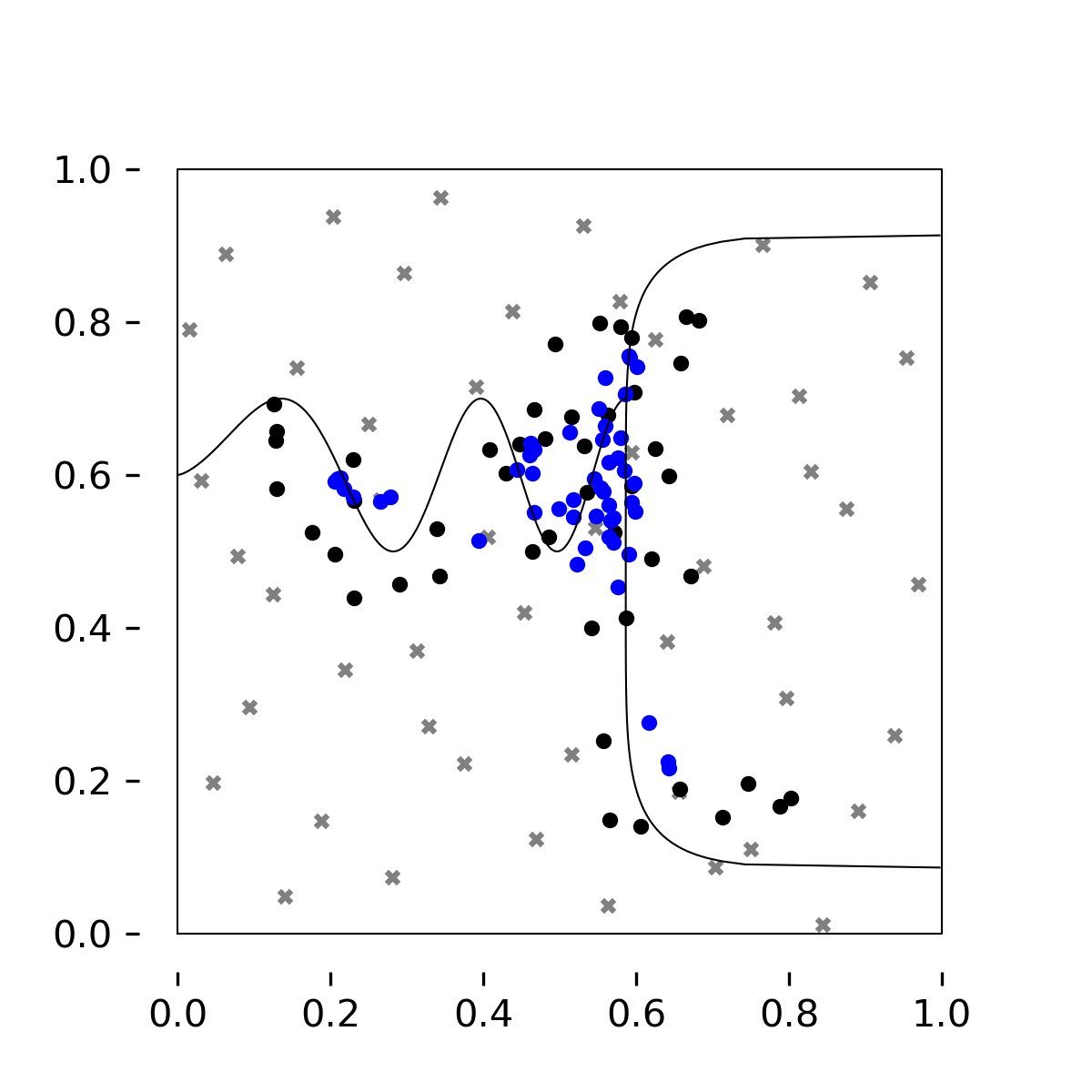}
            \includegraphics[height=6cm]{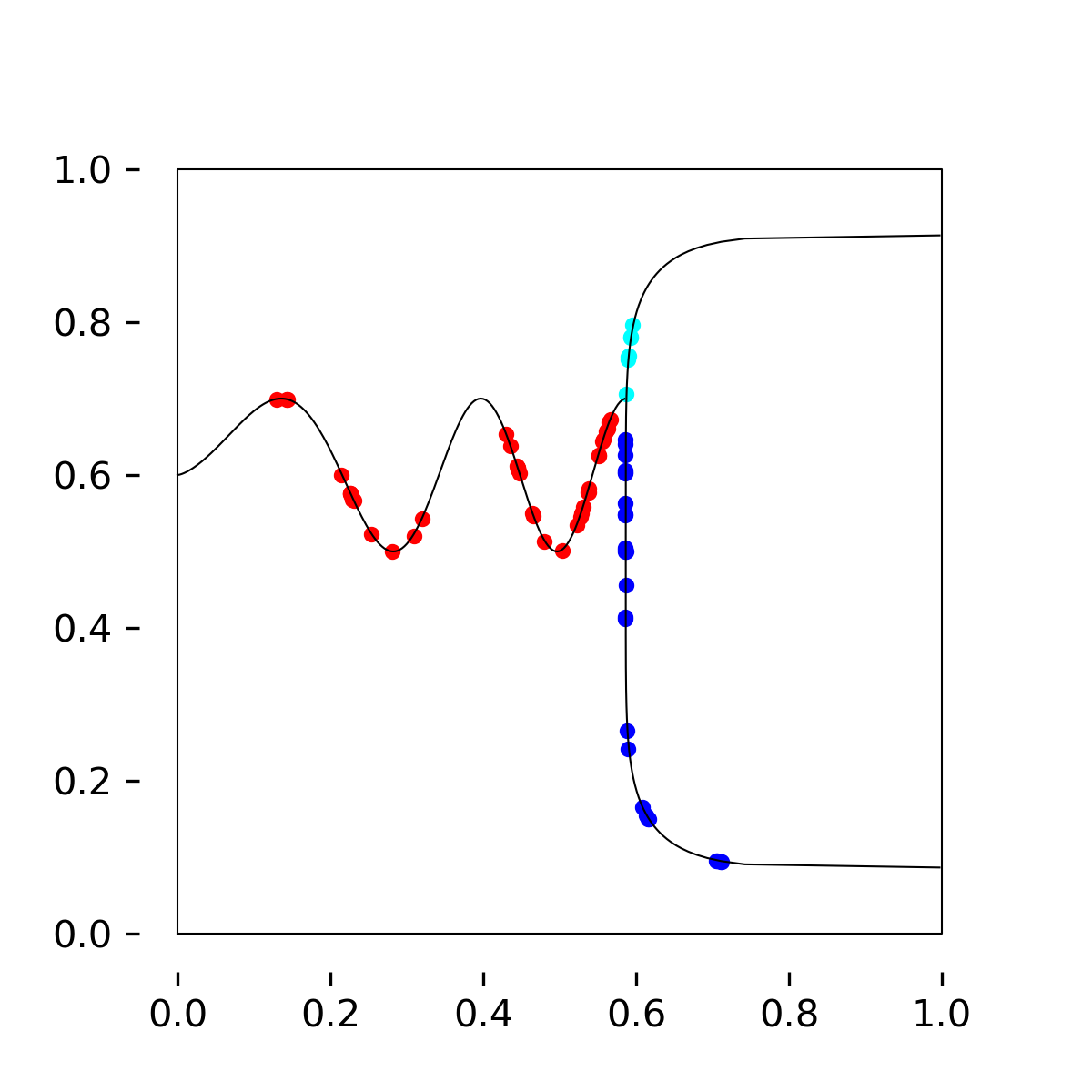}
            \caption{Left: Partition of $\Omega = [0,1]^2$ according to Example~\ref{example1} indicated by grey
            solid lines with the initial point set $X$ displayed as grey crosses.
            We moreover show $\mathcal{M}$ (black points) and $\mathcal{M}^2$ (blue points);
            Right: $S_{1,2}^{(1)}$ (red points), $S_{1,3}^{(1)}$ (blue points), and $S_{2,3}^{(2)}$ (cyan-coloured points).
            }\label{Fig_Subdivision1}
        \end{center}
    \end{figure}

    \begin{Test problem}
        \label{TestProb1}
        For $\Omega = [0,1]^2$, we consider the following partition: Let $\Omega_3 := \{ (x-1)^6 + (y-0.5)^6 < 0.005 \}
        \cap \Omega$, $\Omega_1' := \{ y \leq 0.7 + 0.1\sin(10\pi x^{1.5})\} \cap \Omega$ and
        $\Omega_2' := \Omega \setminus \Omega_1'$.
        Then, we set $\Omega_1 := \Omega_1' \setminus \Omega_3$ and $\Omega_2 := \Omega_2'\setminus \Omega_3$ and study
        the partition $\Omega = \Omega_1 \cup \Omega_2 \cup \Omega_3$ (Fig.~\ref{Fig_Subdivision1}, left).
    \end{Test problem}
    \begin{Example}
        \label{example1}
        For Test problem~\ref{TestProb1}, we choose an initial sampling set $X$ consisting of 50 Halton-distributed points
        (Fig.~\ref{Fig_Subdivision1}, left) and set $k_{\mathrm{near}} = 10$ and $\varepsilon_b = 0.001$.
        The set $\mathcal{M}$ consists of 47 additional points
        represented as black dots; $\mathcal{M}^2$ contains 49 points displayed in blue.
        $\widetilde{S}_{1,2}$ contains 43 triplets, $\widetilde{S}_{1,3}$ contains
        26, and $\widetilde{S}_{2,3}$ contains 7 (Fig.\ \ref{Fig_Subdivision1}, right).
        After cluster removal, we obtain $|S_{1,2}| = 23$, $|S_{1,3}| =
        13$ and $|S_{2,3}| = 4$
        (Fig.~\ref{Fig_ClusterRemoved} and Fig.~\ref{Fig_Subdivision1}, right).
    \end{Example}

    The triplets in $S_{i,j}$ usually provide an incomplete approximation of $\Gamma_{i,j}$ only, feature gaps and lack
    ordering (e.g.\ Fig.~\ref{Fig_Subdivision1}, right).

\begin{figure}
        \begin{center}
            \includegraphics[height=2.5cm]{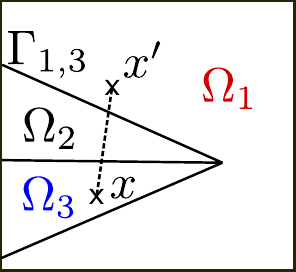}
            \caption{Building block~\texttt{bisection} can fail: Finding points near $\Gamma_{1,3}$ fails for the
            starting points $x \in \Omega_3$ and $x' \in \Omega_1$ shown, as
                $\overline{x'x}$ does not intersect $\Gamma_{1,3}$.}\label{Fig_BisectionFails}
        \end{center}
    \end{figure}

    \begin{remark}
        The bisection-based approach in Building block \texttt{iniapprox} implicitly assumes that $\overline{x'x}$ intersects $\Gamma_{i,j}$ and therefore may fail
        if this does not hold (Fig.\ \ref{Fig_BisectionFails}).
        However, such failure modes occur only rarely in practical computations, and we implemented fallbacks in that case.
    \end{remark}
    Algorithm \texttt{fill} (for an overview, we refer to the flow chart in Fig.\ \ref{fig_bb1_2D}) provides triplets with a maximal
    user-prescribed distance $\varepsilon_b$ to $\Gamma_{i,j}$ and maximal mutual distance $\varepsilon_{\mathrm{gap}}$
    based upon $S_{i,j}$, i.e.\ the result of Algorithm~\ref{alg:initialise} (\texttt{initialise}).
    As we assume that $\varepsilon_b \ll \varepsilon_{\mathrm{gap}}$, we define the distance of two triplets $x,y \in S_{i,j}$
    as $\| x^{(i)} - y^{(i)}\|$ and the distance of $x$ to $\Gamma_{i,j}$ as the distance of $x^{(i,j)}$ to
    $\Gamma_{i,j}$.
    Following a bottom-up approach, we start with discussing selected Building blocks employed in
    \texttt{fill}.

    \begin{Building block}[\texttt{sort}]
        \label{BB_OrderPoints}
        This Building block is to sort a set of triplets $S$ according to their position along $\Gamma$.
        We omit the indices $i$ and $j$ for clarity.
        Let us assume that $\Gamma$ is piecewise smooth and fulfills an inner cone condition with angle $\beta_{\mathrm{angle}}$.
        We first search a triplet $x_{\mathrm{start}} \in S$ closest to the boundary of $\Omega$
        and assume $x_{\mathrm{start}}$ to be the first triplet in the sorted set.
        We initially set $\check{S} = \{x_{\mathrm{start}} \}$.
        Let now the triplets in $\check{S} \subset S$ be already sorted with $x_r$ being the last of those, $r >1$.
        We consider the $k_{\mathrm{sort}}$ nearest triplets $y_1, \ldots, y_{k_{\mathrm{sort}}}$ to $x_r$ in
        $S \setminus \check{S}$, sorted by increasing distance to $x_r$.
        If for $s=1$, $\angle (x_r - x_{r-1}, x_i - y_s) <
        \beta_{\mathrm{angle}}$, we set $x_{r+1} = y_s$ and add $x_{r+1}$ to $\check{S}$; otherwise, we repeat
        with $s+1$.
        If $s > k_{\mathrm{sort}}$, we store $\check{S}$, set $S = S\setminus \check{S}$, and repeat
        the sorting procedure until $S = \emptyset$ ending up with a finite number of disjoint ordered subsets.
        At the end, we combine all sorted subsets based upon the Euclidean distance between first and last points of
        the subsets reversing the order of a subset if necessary.
    \end{Building block}
    \begin{figure}
        \begin{center}
            \includegraphics{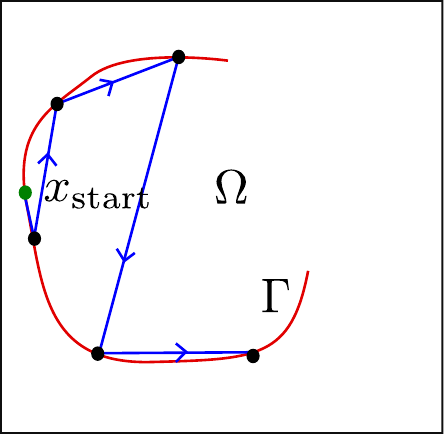}\hspace{1cm}\includegraphics{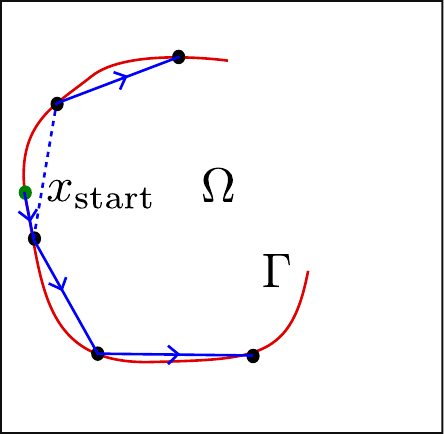}
            \caption{Sorting due to Allasia may fail (left; sorting indicated in blue), whereas our method succeeds in
            the present situation (right). We obtain two ordered subsets (direction of sorting is indicated by arrows
            which after combination yields the correct ordering. The dashed line represents a connection to the
            nearest neighbour rejected due to angle in our approach.}\label{fig_compareSorting}
        \end{center}
    \end{figure}
    \begin{remark}
        Allasia et al.~\cite{Allasia2010EfficientAA} present a simpler sorting method than ours as they do not
        enforce the condition $\angle (x_r - x_{r-1}, x_i - y_s) <
        \beta_{\mathrm{angle}}$.
        However, it may fail if $x_{\mathrm{start}}$ is not the true starting point and if additionally the points are
        unevenly distributed along $\Gamma$ (Fig.\ \ref{fig_compareSorting}) in contrast to ours.
        Of course, our sorting method can fail as well, but due to our
        experience, it is more reliable than Allasia's method and works sufficiently well.
        There are many more sophisticated sorting methods based e.g.\ upon graph
        theory~\cite{Althaus.2001,Amenta.1998,Dey.2000,Ohrhallinger.2013}, which are more reliable than our
        approach, but more time-consuming and much harder to implement.
    \end{remark}
    For describing Building block \texttt{inipairs}, which is part of \texttt{fill}, we
    introduce another two Building blocks, which will be used in \texttt{adapt} as well.

    \begin{Building block}[\texttt{estcurv}]
        \label{bb:EstimateCurv}
        Let be $S_{\mathrm{loc}} = \{x_1, \hdots, x_r\}, r > 2$, a set of ordered points near $\Gamma$ up to $\varepsilon_b$.
        This Building block estimates the curvature $c_{\ell}$ of $\Gamma$ in $x_{\ell}, 1 < \ell < r$ by
        least-squares fitting an approximation using Gaussian radial basis functions (RBFs) and then $c_{\ell}$ by
        the curvature of that approximation in $x_{\ell}$.
        We hereby assume that after shifting and suitable rotation, $\Gamma$ can be locally represented as a graph of an
        unknown function.
        As the points in $X_{\mathrm{loc}}$ are located on $\Gamma$ only up to
        $\varepsilon_b$, we penalise the second derivative of the RBF approximation subject to a
        maximal residual of $\varepsilon_b$.
        This coincides with the maximal deviation in the value at an approximation point.
        We employ Tikhonov regularization with parameter estimation using Morozov's discrepancy principle.

        If $\Gamma$ cannot be considered a graph of a function even after rotation, we draw as a fallback a circle through
        the points with indices $\ell-1$, $\ell$ and $\ell+1$ and use the inverse of its radius for estimating $c_{\ell}$.
        We estimate $c_{\ell}$ in the first or last point of $S_{\mathrm{loc}}$ by drawing a circle through the first or last three
        points in $S_{\mathrm{loc}}$.
    \end{Building block}
    \begin{Building block}[\texttt{esterror}]
        \label{bb:estimate_error}
        This Building block estimates the maximal deviation $\delta$ of a smooth curve from a straight line between
        two points on the curve with distance $d$.
        A straightforward calculation reveals that
        \begin{equation}
            \delta = 0.25cd^2 + 1/16c^3 d^4 + {\cal{O}}(d^6), \label{eq:startpairs}
        \end{equation}
        where $c$ denotes the maximal curvature of the curve between the two points.
        For $S_{\mathrm{loc}}$ as in Building block~\ref{bb:EstimateCurv}, we estimate the maximal deviation $\delta$
        of $\Gamma$ from the straight line between consecutive points $x_{\ell}$ and $x_{\ell +1}$ by replacing $c$ with
        $\max\{c_{\ell}, c_{\ell+1}\}$ from Building block~\ref{bb:EstimateCurv} (\texttt{estcurv}).
        Hereby, we rely on~\eqref{eq:startpairs} and neglect higher order terms.
    \end{Building block}

\begin{figure}
        \begin{center}
            \includegraphics[scale=0.85]{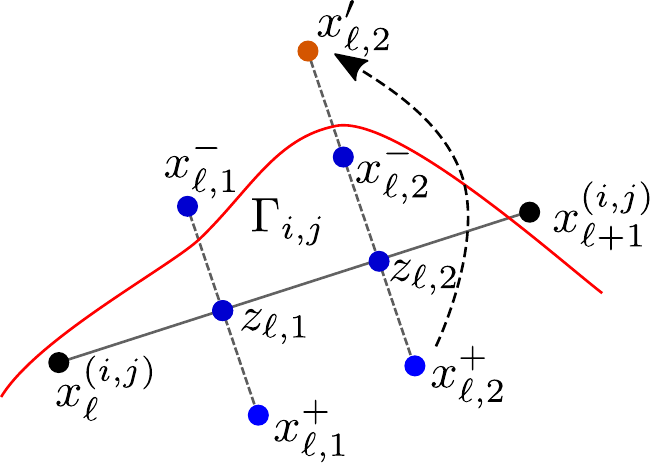}
        \end{center}
        \caption{Creating valid starting pairs in \texttt{fill} with \texttt{startpairs}: While
            $(x_{\ell, 1}^+, x_{\ell, 1}^-)$ is valid, $(x_{\ell, 2}^+, x_{\ell, 2}^-)$ is not, as both points
            belong to the same class. However, $(x_{\ell, 2}^-, x'_{\ell, 2})$ is valid.}\label{fig:AddPointsOnLine}
    \end{figure}

    Now we are prepared to discuss \texttt{fill}.
    \begin{Algorithm}[\texttt{fill}]\label{alg:fill2D}
        Building block~\ref{BB_OrderPoints} (\texttt{sort}) sorts all triplets in $S_{i,j}$ according to their
        position near $\Gamma_{i,j}$.
        We detect gaps in the representation of $\Gamma_{i,j}$ by $S_{i,j}$ by considering subsequent triplets $x_{\ell}$ and
        $x_{\ell+1}$.
        If the distance $d_\ell$ of $x_{\ell}$ to $x_{\ell+1}$ is larger than a user-prescribed threshold
        $\varepsilon_{\mathrm{gap}}$, we consider this a gap and aim to
        equidistantly add $R = \lceil d_\ell/\gamma \rceil$ triplets near $\Gamma_{i,j}$ between $x_{\ell}$ and $x_{\ell+1}$.
        To do so, Building block \texttt{inipairs} places
        new points $z_{\ell, r}, 1\leq r \leq R$, equidistantly on $\overline{x_\ell^{(i,j)} x_{\ell+1}^{(i,j)}}$
        and computes from these initial point pairs
        \begin{equation*}
            x^+_{\ell, r} = z_{\ell, r} + \alpha n, \quad x_{\ell, r}^- = z_{\ell, r} - \alpha n, \quad 1 \leq r \leq R.
        \end{equation*}
        Here, $n$ denotes the (estimated) outer normal unit vector of $\Omega_i$ near $x_{\ell,r}$.
        Applying Building block~\ref{bb:estimate_error} (\texttt{esterror}) to
        $S_{\mathrm{loc}} = \{x_{\ell -2}^{(i,j)}, \hdots x_{\ell +2}^{(i,j)}\}$ and some safeguarding leads to
        \begin{equation}
            \label{eq:comp_alpha}
            \alpha = \min\{\varepsilon_{\mathrm{safemax}} d_{\ell}, \max\{\delta, \varepsilon_{\mathrm{safemin}} \varepsilon_b\}\}
        \end{equation}
        with user-prescribed safety factors $\varepsilon_{\mathrm{safemax}}$ and $\varepsilon_{\mathrm{safemin}}$.
    \end{Algorithm}
    \begin{remark}
        Having a local RBF approximation of $\Gamma$ at hand when computing $c$ in \texttt{estcurv}, it
        seems to be straightforward for efficiency reasons to choose points
        $z_{\ell,r}$ on that
        RBF approximation instead of just subdividing a straight line.
        Numerical experiments did not show any significant advantage of that approach compared to ours for
        \texttt{fill}.
        This could be related to the uneven distribution of points on $\Gamma$ near gaps to fill, which may decrease the quality of approximation.
        Therefore, we stick to the easier approach presented here.
        However, estimating the curvature using Building block~\ref{bb:EstimateCurv} (\texttt{estcurv}) is sufficiently
        reliable for efficiently computing starting pairs.
    \end{remark}
    However, the pairs of starting points (aka starting pairs) obtained from \texttt{inipairs} are not necessarily valid.
    We call a starting pair for approximating $\Gamma_{i,j}$ valid, if one of its points belongs to $\Omega_i$ and the other
    one to $\Omega_j$.
    Therefore, we introduce Building block~\ref{alg_StartValues} (\texttt{startpairs}).
    \begin{Building block}[\texttt{startpairs}]
        \label{alg_StartValues}
        This Building block obtains a valid starting pair from a pair of points $(x_{\ell, r}^+, x_{\ell,r}^-)$.
        If the starting pair is already valid, we return it as the result.
        If $x_{\ell, r}^+$ or $x_{\ell,r}^-$ belongs to a third class, we stop without result.
        If both points belong to the same class, we reflect $x_{\ell, r}^+$ on $z_{\ell, r}$
        obtaining $x_{\ell,r}'$ (Fig.\ \ref{fig:AddPointsOnLine}).
        If both $(x_{\ell, r}^+, x'_{\ell,r})$ and $(x_{\ell, r}^-, x'_{\ell,r})$ still belong to the same class, we repeat
        this process with changing roles and escalating distances at most $k_{\mathrm{rep}}$ times (typically $k_{\mathrm{rep}} = 3$) and
        stop, if any of the resulting point pairs is valid or one of the points belongs to a third class.
    \end{Building block}
    \begin{remark}
        \noindent
        \begin{enumerate}
            \item
            We iterate the process of filling gaps in \texttt{fill}, as the arclength of $\Gamma$
            between two subsequent triplets may considerably exceed the length of the straight line between
            them, such that after a first pass, the mutual distance of subsequent triplets may still exceed
            $\varepsilon_{\mathrm{gap}}$ in some cases.
            \item If only two or even less triplets are known on $\Gamma_{i,j}$, estimating
            curvature is impossible with \texttt{estcurv}, and computing $\alpha$ by~\eqref{eq:comp_alpha} fails.
            For this case, we implemented fallbacks.
        \end{enumerate}
    \end{remark}
    Using \texttt{bisection}, \texttt{fill} creates new triplets, yielding new sets
    $\overline{S}_{i,j}$.

    \begin{figure}
        \begin{center}
            \includegraphics[height=6cm]{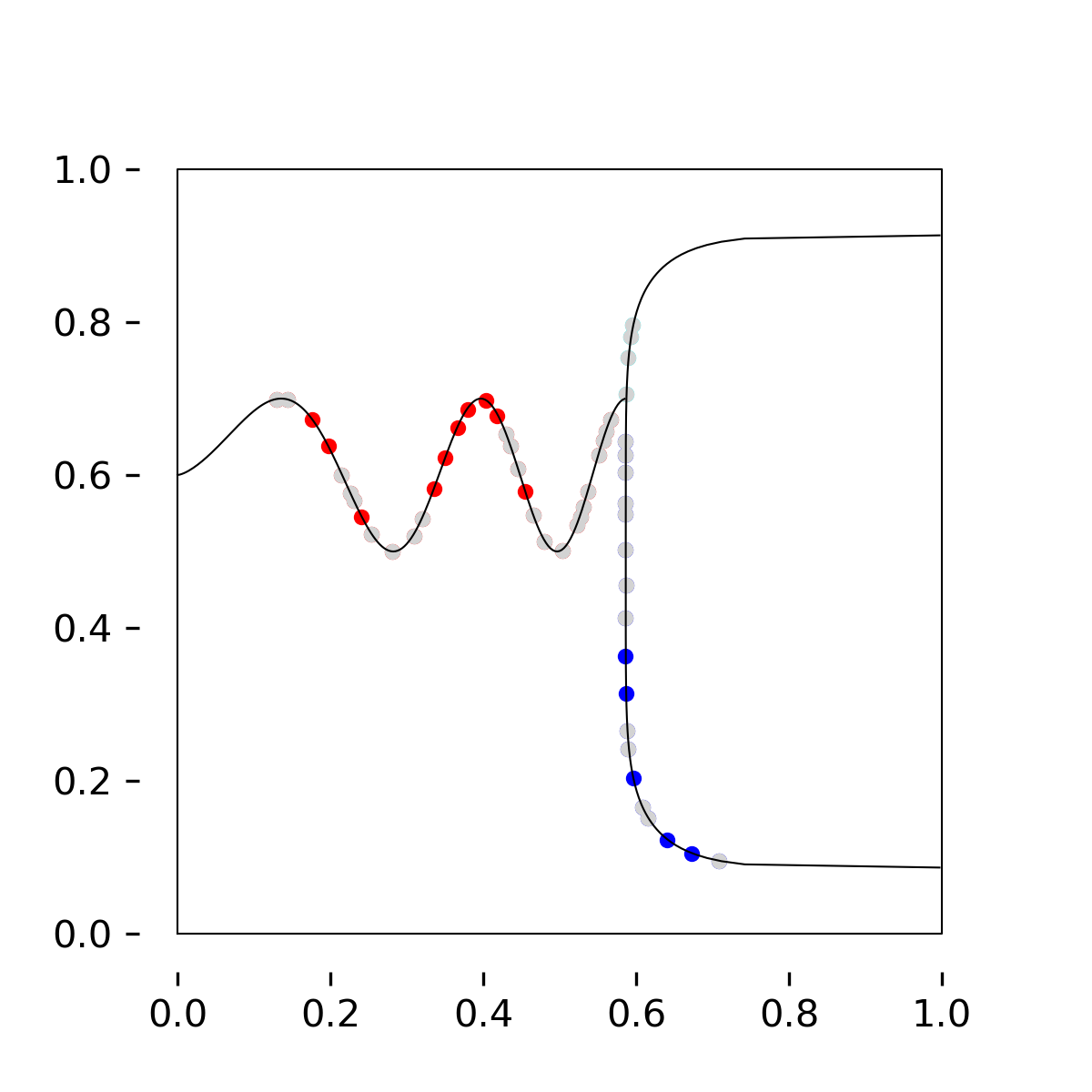}
        \end{center}
        \caption{Algorithm~\ref{alg:fill2D} (\texttt{fill}) for Test problem~\ref{example1}. Previously existing points
        are greyed out in $\overline{S}_{i,j}^{(i)}$. Otherwise, we stick to the
        colouring scheme of Fig.\ \ref{Fig_Subdivision1}.}\label{fig:AddPointsOnLinesExample}
    \end{figure}
    \begin{Example}
        For Test problem~\ref{TestProb1}, we ordered the sets $S_{i,j}$ according to Building block~\ref{BB_OrderPoints}
        with $\beta_{\mathrm{angle}} = \arccos(-0.9) \approx 154^{\circ}$ and $k_{\mathrm{sort}} = 5$.
        By filling gaps with $\varepsilon_{\mathrm{gap}} = 0.05$, $\varepsilon_{\mathrm{safemin}} = 0.95$ and
        $\varepsilon_{\mathrm{safemax}} = 0.25$, we obtain $|\overline{S}_{1,2}| = 33$,
        $|\overline{S}_{1,3}| = 18$, $|\overline{S}_{2,3}| = 4$ (Fig.\ \ref{fig:AddPointsOnLinesExample}).
        We will stick to the values of the parameters given here for all subsequent numerical examples.
        \label{Ex_TestCase01_prolo}
    \end{Example}
    \begin{remark}
        \label{rem_several_components}
        If $\Omega_i$ is not simply connected, some $\Gamma_{i,j}$ may consist of several components.
        We detect this using \texttt{fill}.
        If there are some significant gaps which can not be filled, it indicates
        the presence of several components.
        We then subdivide $\overline{S}_{i,j}$ correspondingly and proceed with every
        subset separately.
    \end{remark}

\begin{figure}
        \begin{center}
            \includegraphics[scale=0.85]{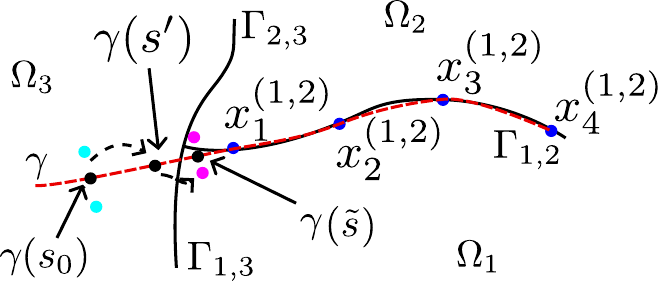}
        \end{center}
        \caption{Scheme of expanding $S_{1,2}$ until its end. Starting from
            $x_1^{(1,2)}, \hdots, x_4^{(1,2)}$ (displayed as blue dots), we construct $\gamma$ by
            \texttt{extrapolate}. As $(x_{s_0}^+, x_{s_0}^-)$, displayed as cyan-coloured dots, fulfills
            \eqref{eq:GammaExceeded}, we apply bisection with respect to the parameter $s$ until a valid starting pair based upon
            $\gamma(\tilde{s})$ can be constructed
            (displayed in magenta), from which we compute the final triplet in $S_{1,2}$ by \texttt{bisection}.}\label{fig:end_of_curve}
    \end{figure}

    Fig.\ \ref{fig:AddPointsOnLinesExample} indicates that even with filling gaps, $\overline{S}_{i,j}$ may not appropriately represent
    $\Gamma_{i,j}$, as parts of $\Gamma_{i,j}$ before the first known triplet $x_1$ and after the last known
    may be neglected.
    Algorithm \texttt{expand} is used to expand $\overline{S}_{i,j}$ to a representation of the complete curve $\Gamma_{i,j}$.
    It relies on several building blocks, which we discuss first.
    \begin{Building block}[\texttt{extrapolate}]\label{bb:extrapolate2D}
        For a given ordered set $S$ with $n \geq 2$ triplets close up to $\varepsilon_b$ to a curve $\Gamma$ and with
        average distance $d_{\mathrm{avg}}$, we fit a polynomial with degree $n-1$ in local coordinates.
        We compute these coordinates by least-squares-fitting a line to $S$.
        As points in $S$ are located on $\Gamma$ up to $\varepsilon_b$ only, we do not interpolate, but penalise the
        second derivative in a least-squares approximation.
        Following Morozov's discrepancy principle, we regularise such that the maximal residual is
        approximately\ $\varepsilon_b$.
    \end{Building block}
    \begin{Building block}[\texttt{stepsize}]\label{bb:stepsize}
        Provided that $\Gamma$ extends sufficiently far before $x_1 \in S$, it seems to be
        straightforward to seek for a new triplet with distance $\varepsilon_{\mathrm{gap}}$
        to $x_1$.
        However, we limit the step size for extrapolation based upon the curvature $c$ of $\Gamma$ in the vicinity of
        $x_1$.
        Extrapolating far is unreliable in case of large curvature $c$, and \texttt{adapt} will insert additional
        points afterwards in that region anyway for accuracy reasons, such that it is much more efficient to adjust the
        step length to the local properties of $\Gamma$ beforehand.
        Let us assume that a polygonal final approximation of $\Gamma$ may deviate at most by
        $\varepsilon_{\mathrm{err}}$ from $\Gamma$.
        We compute the step size which would lead to a deviation of $\varepsilon_{\mathrm{err}}$ from a straight line
        segment between $x_1$ and the new triplet yet to compute.
        This is a natural upper bound for the step size $l_{\mathrm{extra}}$ in extrapolation.

        Rearranging~\eqref{eq:startpairs} and neglecting higher order terms leads to
        \begin{equation}
            l_{\max} = \frac{2}{c^2} \left(\sqrt{1+4c\varepsilon_{\mathrm{err}}} - 1 \right)\label{eq:maxsteplength}
        \end{equation}
        for the maximal admissible step length $l_{\max}$.
        However, evaluating~\eqref{eq:maxsteplength} is numerically unstable if $c \varepsilon_{\mathrm{err}}$ is small.
        We set $(c\varepsilon_{\mathrm{err}})^2 = v$
        and search for the roots $v_{\min}$ and $v_{\max}$ of $v^2 + 4 v - 16cd$.
        According to Vieta, $v_{\min} = - \left(2 + \sqrt{4+ 16cd}\right)$ and $v_{\max} = -16cd/v_{\min}$.
        Resubstituting $v$ yields $l_{\max} = 4 \sqrt{d/(-c v_{\min})}$.
        Some safeguarding of this result leads to a step length of
        \begin{equation*}
            l_{\mathrm{extra}} = \min\{\varepsilon_{\mathrm{gap}},  \beta_{\mathrm{growth}} d_{\mathrm{avg}}, l_{\max} \} \: ,
        \end{equation*}
        where we estimate $c$ using Building block \ref{bb:EstimateCurv} (\texttt{estcurv}) applied to
        $S_{\mathrm{loc}} = \{x_1^{i,j)}, \hdots, x_{k_{\text{extra}}}^{(i,j)} \}$
        with a user-defined parameter $k_{\text{extra}}$.
        The term $\beta_{\mathrm{growth}} d_{\mathrm{avg}}$ increases robustness, as extrapolation is reliable only
        near the points to extrapolate, and it may happen that $d_{\mathrm{avg}} \ll \varepsilon_{\mathrm{gap}}$.
    \end{Building block}
    \begin{Algorithm}[\texttt{expand}]\label{alg:expand2D}
        This algorithm aims at finding triplets near $\Gamma_{i,j}$ beyond the first or last known
        in $\overline{S}_{i,j}$ until the start or end of $\Gamma_{i,j}$ is reached or $\Gamma_{i,j}$ turns out to be a
        closed curve.
        For the sake of simplicity, we refer in what follows to finding triplets before the first one in
        $\overline{S}_{i,j}$.
        Finding triplets near $\Gamma_{i,j}$ beyond the last one in $\overline{S}_{i,j}$ works analogously.
        We add a new triplet before $x_1$ by extrapolating an approximation $\gamma$ of $\Gamma_{i,j}$ with
        Building block~\ref{bb:extrapolate2D} (\texttt{extrapolate}) setting
        $S = \left\{x_1^{(i,j)}, \hdots, x_{k_{\mathrm{extra}}}^{(i,j)}\right\}$ with $x_{\ell} \in \overline{S}_{i,j}$
        and choose the step size according to Building block~\ref{bb:stepsize} (\texttt{stepsize}).
        This way, we obtain an extrapolating curve $\gamma$
        and some $s_0$ such that $\|x_1 - \gamma(s_0)\| \approx l_{\mathrm{extra}}$.
        We then create a point pair $(x_{s_0}^+, x_{s_0}^-)$ based on $\gamma(s_0)$ in a similar way as in \texttt{fill}.
        Computing a valid starting pair with Building block~\ref{alg_StartValues} (\texttt{startpairs}) and subsequent
        bisection yields a new triplet in $\overline{S}_{i,j}$.
        We repeat this process until we reach the true starting point of $\Gamma_{i,j}$.
        As heuristic criterion for $\gamma(s_0)$ exceeding this starting point, we consider
        \begin{equation}
            \label{eq:GammaExceeded}
            x_{s_0}^+ \not\in (\Omega_i \cup \Omega_j) \vee x_{s_0}^- \not\in (\Omega_i \cup \Omega_j)
        \end{equation}
        (Fig.~\ref{fig:end_of_curve}).
        In this case we employ Building block~\ref{bb:reducestepsize} (\texttt{reducestepsize}) for obtaining a
        valid pair of points $(x_{\tilde{s}}^+, x_{\tilde{s}}^-)$, which represents the starting point of
        $\Gamma_{i,j}$. After adding it to $\overline{S}_{i,j}$, \texttt{expand} terminates.
        If $\Gamma_{ij}$ is closed, expanding $\overline{S}_{i,j}$ as described above would lead to an endless loop.
        Therefore, we start expanding $\overline{S}_{i,j}$, but detect after every addition of a triplet, if
        $\Gamma_{i,j}$ is closed.
        If so, we stop expanding and resort.
        We skip the details to keep the presentation uncluttered.
        Algorithm \texttt{expand} yields approximating sets $\breve{S}_{i,j}$.
    \end{Algorithm}
    \begin{Building block}[\texttt{reducestepsize}]\label{bb:reducestepsize}
        Using the parameter value $s_1$ corresponding to $x_1^{(i,j)}$ as lower bound and $s_0$ as upper bound, we
        obtain $\tilde{s}$, which leads to a valid starting pair $(x_{\tilde{s}}^+, x_{\tilde{s}}^-)$ based upon
        $\gamma(\tilde{s})$ and fulfils
        $\|\gamma(\tilde{s}) - \gamma(s')\| < \varepsilon_b$ by bisection with respect to
        condition~\eqref{eq:GammaExceeded}.
        Here, $s'$ denotes the second to last parameter in the bisection process (compare Fig.~\ref{fig:end_of_curve}).
    \end{Building block}

    \begin{figure}
        \begin{center}
            \includegraphics[height=6cm]{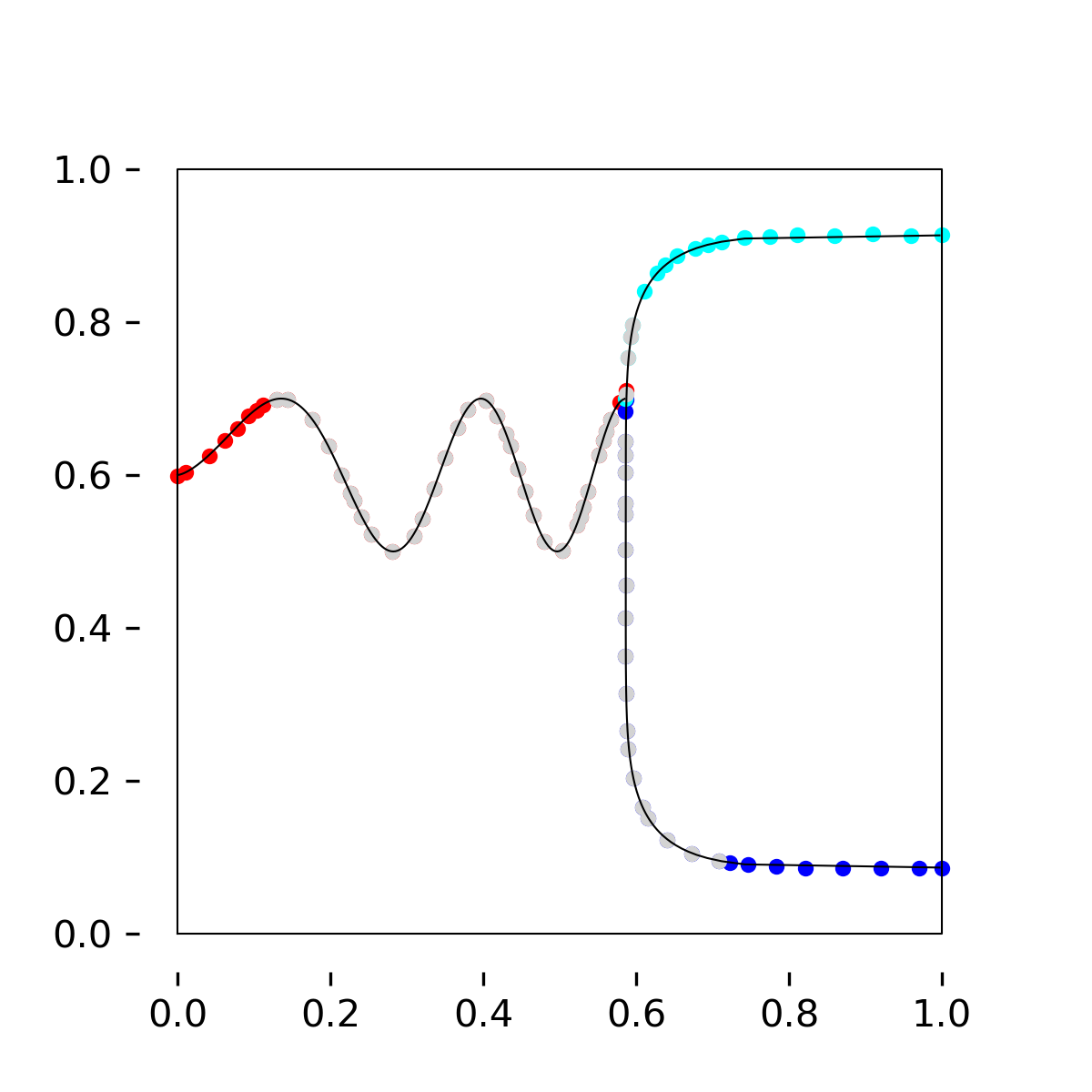}
            \includegraphics[height=6cm]{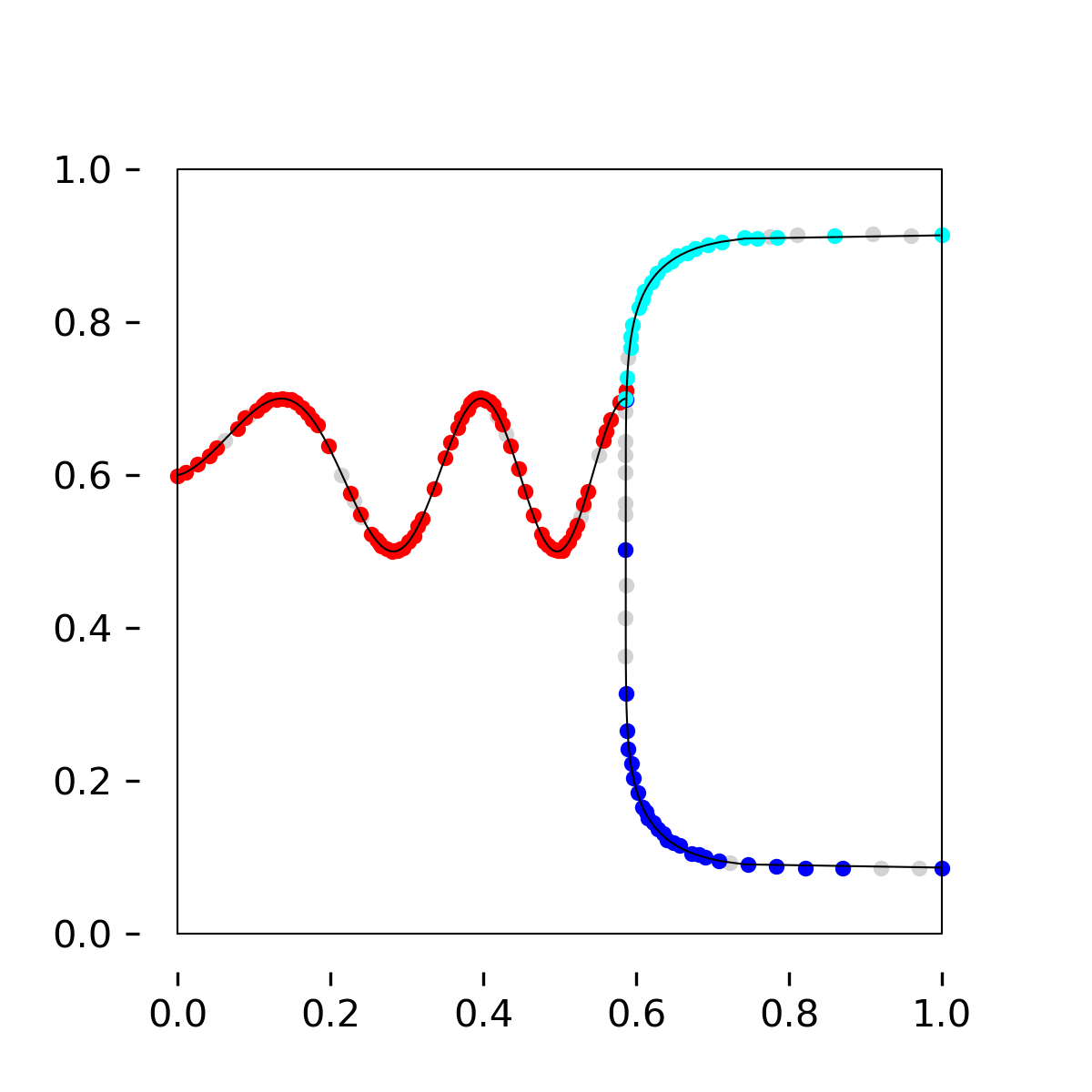}
            \caption{Sets $\breve{S}_{i,j}^{(i)}$ (left) and
                $\hat{S}_{i,j}^{(i)}$ (right) for Test problem~\ref{example1}. We stick to the
                colouring scheme of Fig.\ \ref{Fig_Subdivision1}.}\label{Fig_Prolongated_Lines}
        \end{center}
    \end{figure}
    \begin{Algorithm}[\texttt{adapt}]
        Based on \texttt{esterror}, this algorithm inserts a triplet (approximately) halfway between
        consecutive triplets $x_{\ell}$ and $x_{\ell+1}$, if \texttt{esterror} indicates an
        error larger than $\varepsilon_{\mathrm{err}}$ and removes a triplet, if the estimated error of both line segments the
        triplet belongs to is smaller than $\varepsilon_{\mathrm{coarse}}$.
        In contrast to \texttt{fill}, we employ the local RBF
        approximation from \texttt{estcurv} for computing an initial point $x_{\mathrm{new}}$ between
        $x_{\ell}$ and $x_{\ell+1}$ when refining.
        With $x_{\mathrm{new}}^{\pm} = x_{\mathrm{new}} \pm \alpha' n$, we then proceed as in \texttt{fill}.
        We compute $\alpha'$ according to~\eqref{eq:comp_alpha}, but replace $\delta$
        by $\delta' = 1/16c^3 d^4$, as the error due to~\eqref{eq:startpairs} refers to an approximation by line
        segments and is overly pessimistic for an RBF approximation.
        For robustness, we never delete consecutive triplets in one pass of the adaptive loop.
        After at most $k_{\mathrm{adap}}$ refinement and coarsening sweeps, we end up with final sets $\hat{S}_{i,j}$.
    \end{Algorithm}
    \begin{Example}
        For Test problem~\ref{TestProb1}, we continue our calculations from Example~\ref{Ex_TestCase01_prolo}.
        We set $k_{\mathrm{extra}} = 4, \varepsilon_{\mathrm{err}} = 0.001$, $\varepsilon_{\mathrm{coarse}} = 0.0001$
        and obtain $|\breve{S}_{1,2}| = 43$,
        $|\breve{S}_{1,3}| = 28$, and $|\breve{S}_{2,3}| = 18$
        (Fig.\ \ref{Fig_Prolongated_Lines}, left).
        Algorithm \texttt{adapt} with $k_{\text{adap}} = 4$ yields $|\hat{S}_{1,2}| = 77$,
        $|\hat{S}_{1,3}| = 26$, and $|\hat{S}_{2,3}| = 22$
        (Fig.\ \ref{Fig_Prolongated_Lines}, right).
        Computing these sets requires 1088 classifications.
    \end{Example}
    \begin{figure}
        \begin{center}
            \includegraphics[height=4.5cm]{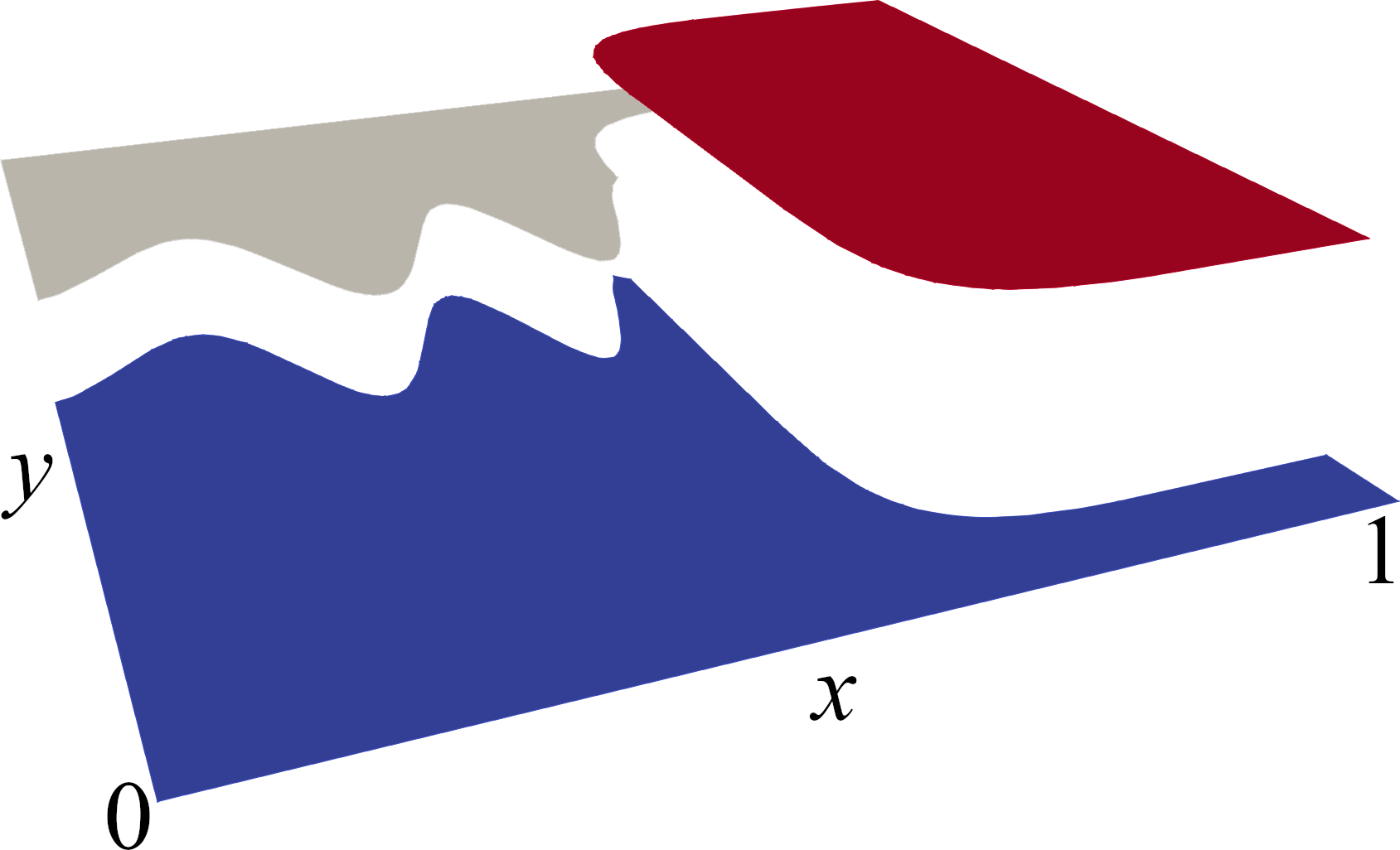}
            \caption{Reconstructed subdivision for Test problem~\ref{example1}.}
        \end{center}\label{fig:ExampleFinal}
    \end{figure}

    \subsection{Approximating $\Gamma_{i,j}$ in 3D}\label{subsec:approximation_in_3D}
\begin{figure}
        \begin{center}
            \includegraphics[height=6cm]{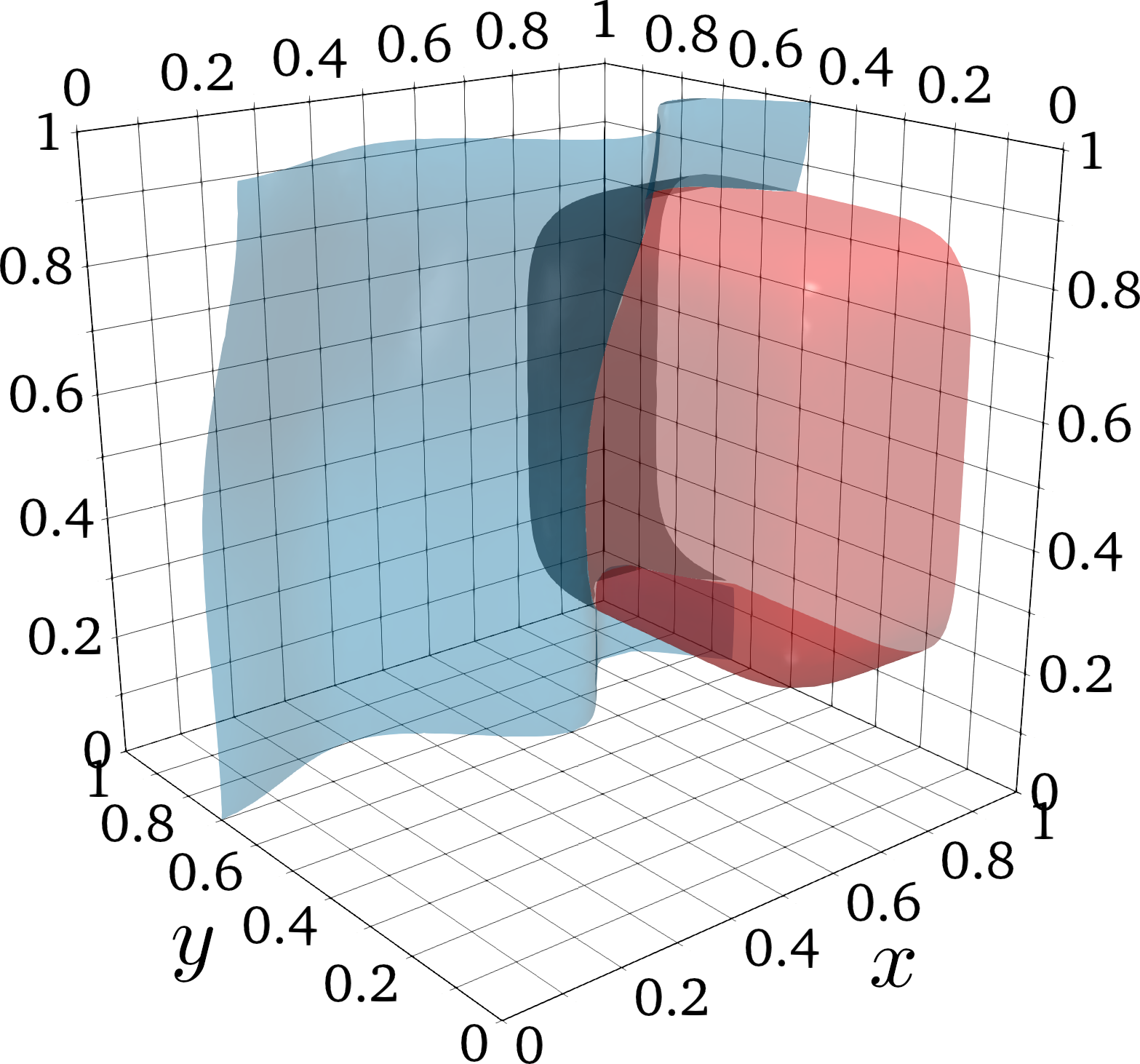}\hspace{1cm}
            \includegraphics[height=6cm]{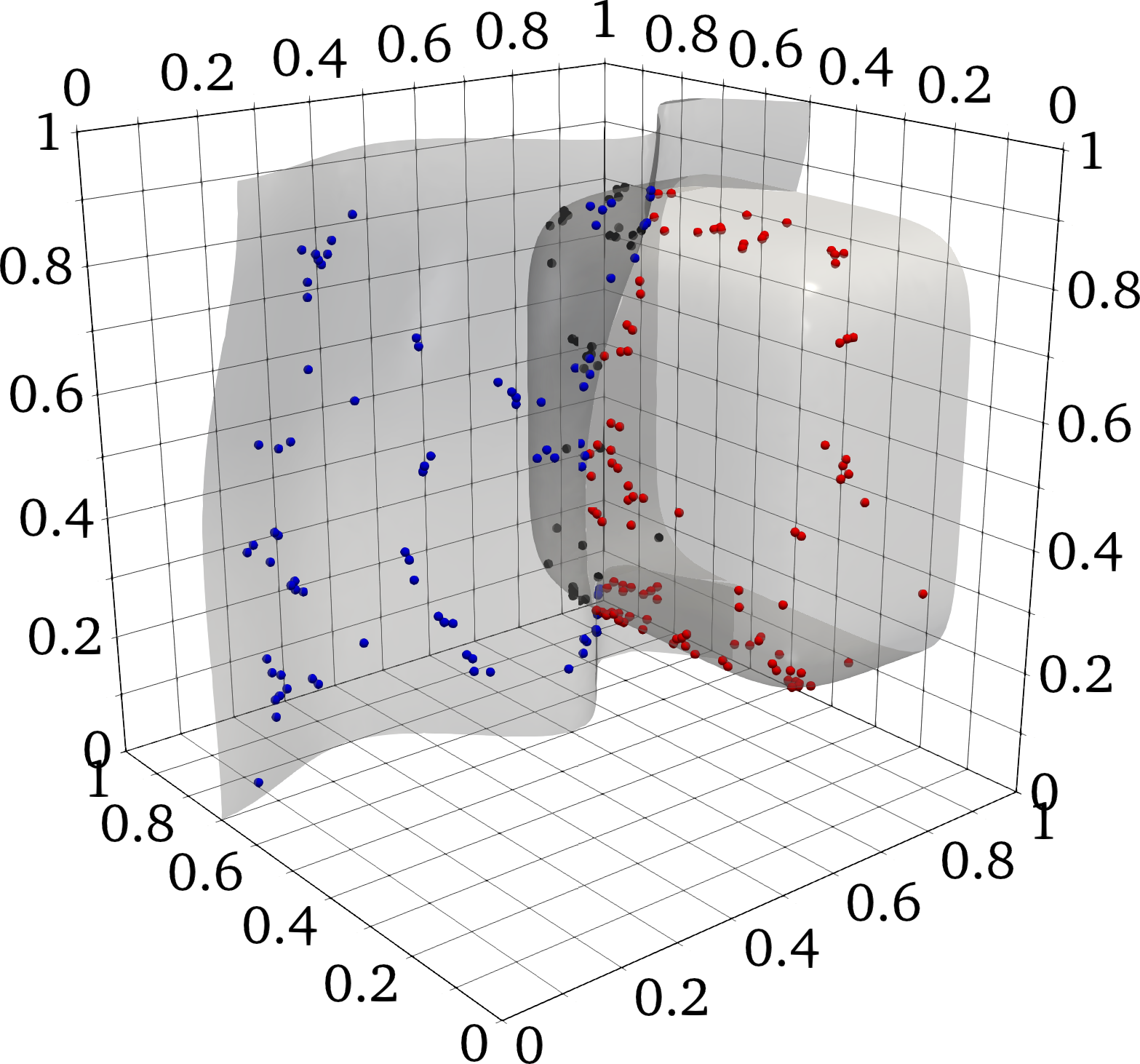}
            \caption{Left: Partition of $[0,1]^3$ for Test problem~\ref{prob:3D_1} displayed by $\Gamma_{1,2}$
                (red), $\Gamma_{1,3}$ (blue) and $\Gamma_{2,3}$ (grey; partially obscured by $\Gamma_{1,3}$); Right:
                approximating sets $S_{i,j}^{(i)}$ after removing clusters.}\label{Fig:TestCase3D}
        \end{center}
    \end{figure}

    For approximating $\Gamma_{i,j}$ in three dimensions, we stick to the general procedure for approximating these sets
    in two dimensions (Fig.\ \ref{fig_flow_chart}).
    While \texttt{initialise} remains unchanged, \texttt{fill}, \texttt{expand} and \texttt{adapt} differ from their 2D
    counterparts, as these exploit ordering of the points on a decision curve or fault line.
    There is however no straightforward ordering of points on a surface.
    Because \texttt{fill} and \texttt{adapt} rely on \texttt{esterror} as in two dimensions, we discuss error estimation
    first.
    \begin{Building block}[\texttt{esterror}]\label{bb:estError3D}
        This Building block aims at estimating the maximal error $e$ of a linear triangulation of $\Gamma$ based upon a
        finite set of points $S$ near $\Gamma$ up to $\varepsilon_b$.
        Let $\mathcal{T}$ be a Delaunay triangulation of $S_{\text{loc}} \subset S$ and let us assume that $\Gamma$
        can be represented on the support of $\mathcal{T}$ by an unknown function $g$ after appropriate change of
        coordinates.
        For some triangle $T \in \mathcal{T}$, it holds according to~\cite[Theorem 4.1]{Waldron2006}
        \begin{equation*}
            \|g - I_T g \|_{T, \infty} \leq \frac{1}{2}\left(R^2-d^2\right) |g|_{2,\infty,T}.
        \end{equation*}
        Here, $R$ describes the radius of the circumcircle of $T$, $d$ the distance from its center to $T$ and $I_T g$
        the Lagrange interpolant of $g$ on $T$; with $|g|_{2,\infty,T}$, we denote the $L_{\infty}$-norm of the second
        derivative of $g$ on $T$.
        It remains to estimate $|g|_{2,\infty,T}$.
        As $g$ is unknown, we employ an RBF approximation $\varphi$ as in Building block \ref{bb:EstimateCurv} (\texttt{estcurv})
        instead and approximate $|g|_{2,\infty,T}$ by evaluating its second derivative in the vertices of $T$ and its center.
        The maximum of these four values yields the desired approximation $\phi$ of $|g|_{2,\infty,T}$ and therefore
        \begin{equation}
            e \approx \frac{1}{2}\left(R^2-d^2\right) \phi
        \end{equation}
    \end{Building block}
    \begin{Algorithm}[\texttt{fill}]\label{alg:fill3D}
        For any triplet $x$ in $S_{i,j}$, we search the $k_{\mathrm{near}}$ nearest neighbours
        $x_1, \hdots, x_{k_{\text{near}}}$ and switch to a local 2D
        coordinate system by computing the optimal fitted plane in the sense that the sum of
        the squared distances between the points and the plane is minimal (see~\cite{Shakarji.1998}).
        We then compute in local coordinates a Delaunay triangulation of the patch
        $S_{\text{loc}} = \{x, x_1, \hdots, x_{k_{\text{near}}}\}$.
        If the maximal edge length of a triangle in this triangulation exceeds $\varepsilon_{\text{gap}}$, we mark
        the centre of that triangle as starting point for constructing a valid starting pair similar as in Building
        \texttt{inipairs}, as long as the triangle is not too anisotropic, i.e.\ does not contain angles approaching
        $0^{\circ}$ or $180^{\circ}$.
        However, large gaps in $S_{i,j}$ are not covered by this local approach (Fig.~\ref{Fig:FillGapsin3D}).
        If there is a gap in $S_{i,j}$ and if the current point is at the boundary of this gap, it will be at the
        boundary of the current patch.
        We detect this by exploiting the local coordinate system.
        We assume that a patch in the vicinity of a large gap is at least slightly elongated in tangential direction to
        the patch boundary.
        Therefore, we consider the first local coordinate.
        If the first local coordinate of the current point is almost the minimum or maximum of all respective
        first coordinates of the patch, then it is at its boundary, and we assume a large gap to fill.
        In this case, we construct $x_{\text{new}}$ on the elongated line between the centre of gravity of the patch
        and $x$, enrich $S_{\text{local}}$ by $x$ and continue as described above.
        Looping over all triplets in $S_{i,j}$, removing duplicate starting points and starting points which are
        extremely close enables us to compute triplets using Building block \texttt{bisection} similar as in
        Algorithm~\ref{alg:fill2D} (\texttt{fill}).
        This yields an enriched representation of $\Gamma_{i,j}$.
        Repeating this procedure with the enriched set until no gaps are detected anymore leads to $\overline{S}_{i,j}$.
    \end{Algorithm}
    \begin{figure}
        \begin{center}
            \includegraphics{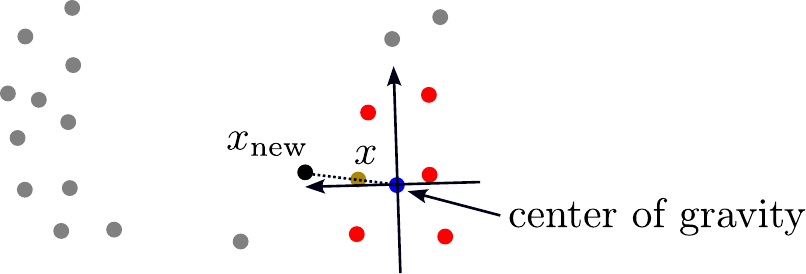}
            \includegraphics[height = 4.5cm]{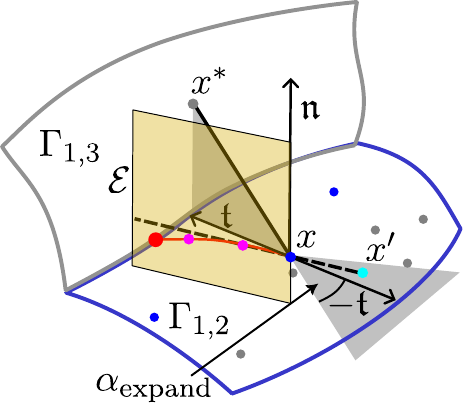}
            \caption{Left: Gradually filling large gaps in $S_{i,j}$ (displayed as grey dots).
            The triplet $x$ is represented in golden colour; the triplets in
                $S_{\text{local}}$, displayed in red, define a local coordinate system sketched with black arrows.
                The centre of gravity of $S_{\text{local}}$ is shown as a blue dot. Right:
                Expansion of $\overline{S}_{1,2}$ to the inner boundary $\Gamma_{1,2} \cap \Gamma_{1,3}$. Triplets
            in $\overline{S}_{1,2} \setminus \overline{S}_{1,2}^I$ are displayed as grey dots, triplets in
                $\overline{S}_{1,2}^I$ as blue dots. 2D-Expanding $\{x, x'\}$ on $\mathcal{E} \cap \Gamma_{1,2}$
            yields a new triplet near $\Gamma_{1,2} \cap \Gamma_{1,3}$, displayed as red dot.}\label{Fig:FillGapsin3D}
        \end{center}
    \end{figure}
    \begin{Test problem}\label{prob:3D_1}
        For $\Omega = [0,1]^3$, we consider the following partition: Let $x_M = (1, 0.5, 0.5)^{\top}$ and
        $\Omega_3 := \{ \|x - x_M\|_6 < 0.002 \} \cap \Omega$,
        $\Omega_1' := \{ y + 0.1z < 0.7 + 0.1\sin(10x^{1.5}) + 0.05\sin(5z^{1.5})\} \cap \Omega$ and
        $\Omega_2' := \Omega \setminus \Omega_1'$.
        Then, we set $\Omega_1 := \Omega_1' \setminus \Omega_3$ and $\Omega_2 := \Omega_2'\setminus \Omega_3$ and study
        the partition $\Omega = \Omega_1 \cup \Omega_2 \cup \Omega_3$ (Fig.~\ref{Fig:TestCase3D}).
    \end{Test problem}

    \begin{Algorithm}[\texttt{expand}]\label{alg:expand3D}
        When expanding a representation $\overline{S}_{i,j}$ of $\Gamma_{i,j}$, we distinguish between expansion towards
        inner boundaries $\Gamma_{i,j} \cap \Gamma_{k,\ell}$ and expansion towards outer boundaries
        $\Gamma_{i,j} \cap \partial \Omega$.
        We consider expanding to inner boundaries first.
        Let $\overline{S}_{i,j}^I \subseteq \overline{S}_{i,j}$ be all triplets
        closer than $0.75\varepsilon_{\text{gap}}$ to a triplet in another $S_{k,\ell}$.
        For some $x \in \overline{S}_{i,j}^I$, let $x^*$ be the triplet in $S_{k,\ell}$ closest to $x$.
        We estimate the normal vector $\mathfrak{n}$ of $\Gamma_{i,j}$ in $x$ and project $\overline{xx^*}$ on the
        corresponding tangential plane at $x$, yielding $\mathfrak{t}$.
        Let $x' \in \overline{S}_{i,j}$ be the nearest triplet to $x$ fulfilling
        $\angle (-\mathfrak{t}, \overline{xx'}) < \alpha_{\text{expand}}$
        and $\mathcal{E}$ the plane which contains $x$ and is spanned by $\overline{xx'}$ and $\mathfrak{n}$.
        We reduce the expansion to the two-dimensional case by expanding the curve $\Gamma_{i,j} \cap \mathcal{E}$
        with  Algorithm~\ref{alg:expand2D} (\texttt{expand}) from Section~\ref{subsec:2d_description}, using $x$ and $x'$ as initial set of
        triplets (compare Fig.~\ref{Fig:FillGapsin3D}, right).

        For expanding to outer boundaries, we construct $\overline{S}_{i,j}^B$ as follows: For each coordinate
        direction $i$, we select $n_{\text{expand}, i}$
        points in $S_{i,j}$ which are minimal or maximal with respect to this coordinate and set
        $\overline{S}_{i,j}^B = \overline{S}_{i,j}^B \setminus \overline{S}_{i,j}^I$ in order to avoid duplicate triplets.
        For determining $n_{\text{expand}, i}$, we rely on the axis-parallel bounding box of $\overline{S}_{i,j}$ with
        sizes $b_1, b_2, b_3$.
        Then, $n_{\text{expand}, i} = \lceil \max\{b_{i +1 \mod 3}, b_{i \mod 3 + 1}\}/ \varepsilon_{\text{gap}} \rceil$.
        From these triplets, we select the ones which are either closer than $\varepsilon_{\text{gap}}$ to one of the
        boundary facets or fulfil
        \begin{equation}
            \angle (\mathfrak{n}, \mathfrak{n}_{\text{outer}}) > \alpha_{\text{expbound}}\label{eq:angle_cond},
        \end{equation}
        where $\mathfrak{n}_{\text{outer}}$ stands for the outer normal vector of the assigned boundary facet.
        Figure~\ref{fig:expand_to_boundary} illustrates the purpose of condition \eqref{eq:angle_cond}.
        For any triplet in $\overline{S}_{i,j}^B$, we proceed as for expanding to inner boundaries, but with
        $\mathfrak{t} = \mathfrak{n}_{\text{outer}}$.
        All new triplets in the same facet $F$ of $\partial \Omega$ constitute the set $\overline{S}_{i,j}^{F}$.
        As this set represents the curve $\Gamma_{i,j} \cap F$ in the plane $F$, we are confronted with approximating a
        decision curve or fault line in two dimensions.
        Therefore, we apply Algorithms~\ref{alg:fill2D} (\texttt{fill}) and \ref{alg:expand2D} (\texttt{expand}) to
        $\overline{S}_{i,j}^{F}$.
        After adding the new triplets on the boundary of $\Gamma$, we apply Algorithm \ref{alg:fill3D} again and end up
        with an enlarged set $\breve{S}_{i,j}$.
    \end{Algorithm}

    \begin{figure}
        \begin{center}
            \includegraphics[height=6cm]{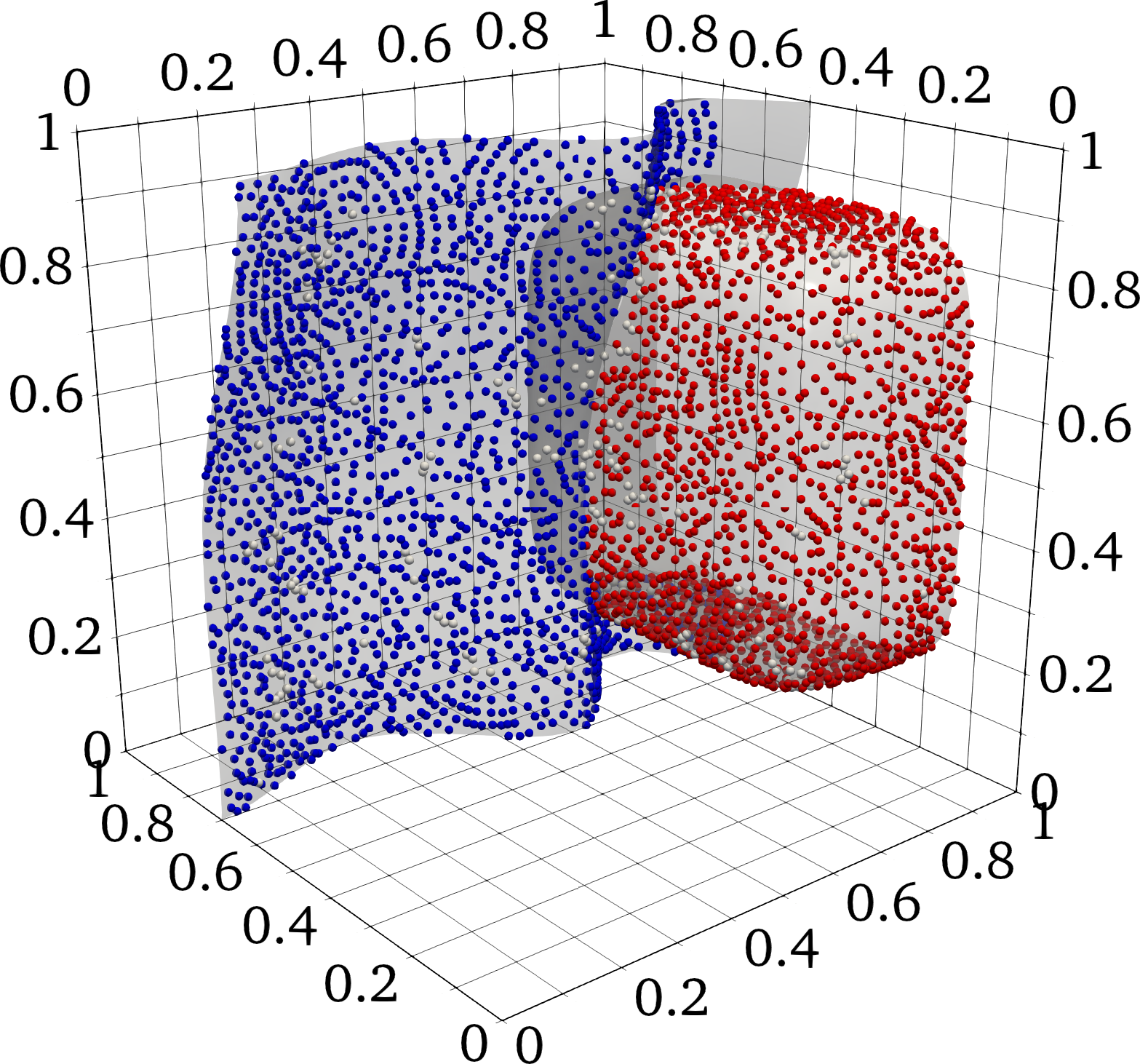}\hspace{1cm}
            \includegraphics[height=6cm]{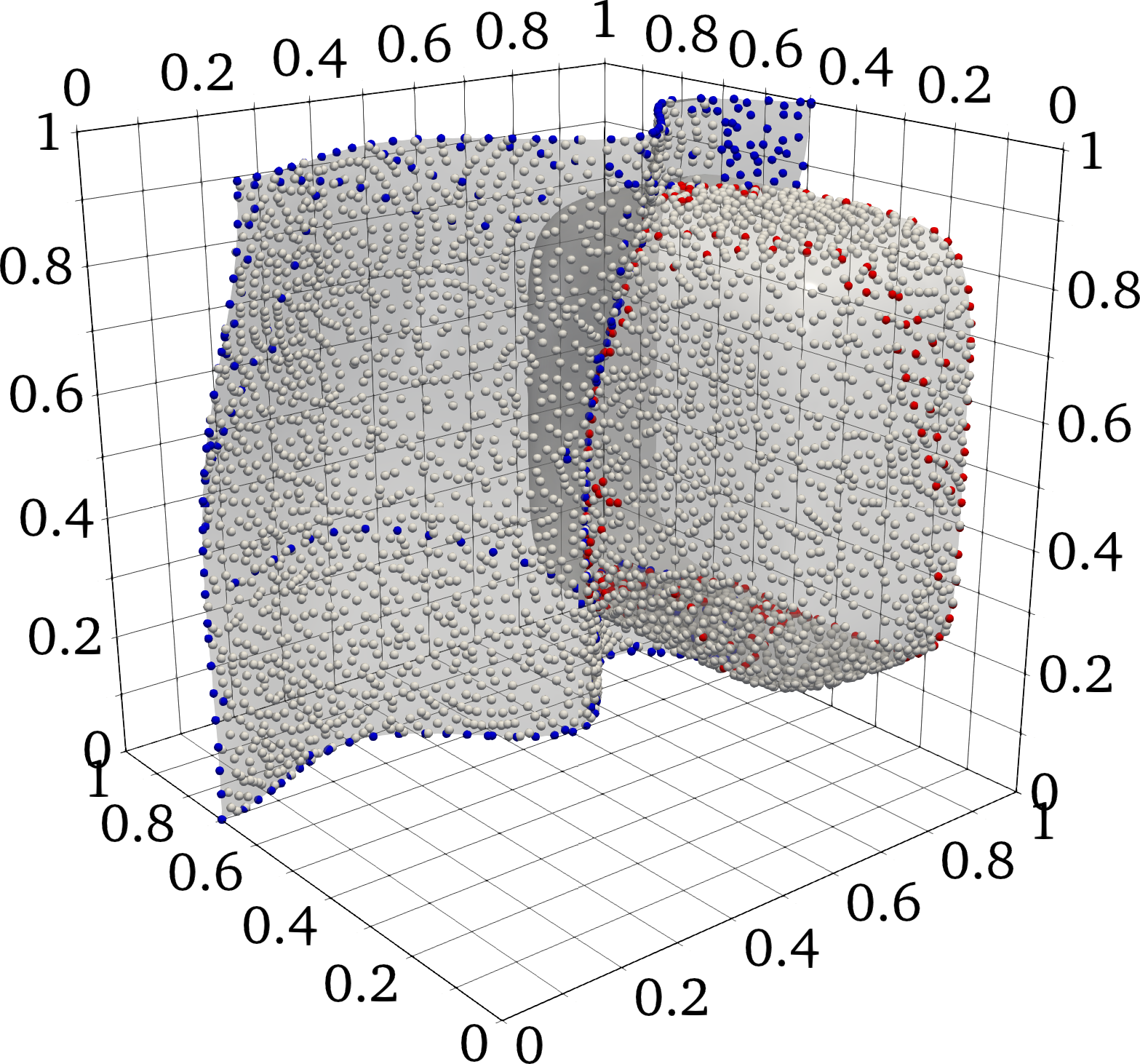}
            \caption{Left: Sets $\overline{S}_{1,2}^{(1)}$, $\overline{S}_{1,3}^{(1)}$ obtained applying
            Algorithm~\ref{alg:fill3D} (\texttt{fill}) for Test
                problem~\ref{prob:3D_1}. Right: $\breve{S}_{1,2}^{(1)}$, $\breve{S}_{1,3}^{(1)}$
                after Algorithm \ref{alg:expand3D} (\texttt{expand}) for the same Test problem.
            Previously existing points are displayed in light grey;
            we omit $\overline{S}_{2,3}^{(2)}$ and $\breve{S}_{2,3}^{(2)}$ in the visualisation for the sake of clarity.}\label{Fig:TestCase3D_fill}
        \end{center}
    \end{figure}
    \begin{figure}
        \begin{center}
            \includegraphics{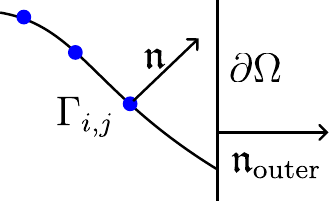}\hspace{1cm}
            \includegraphics{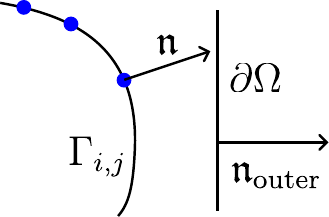}
            \caption{Suitable (left) and unsuitable (right) triplets for expanding to a domain boundary.}\label{fig:expand_to_boundary}
        \end{center}
    \end{figure}

    \begin{Algorithm}[\texttt{adapt}]
        In contrast to the two-dimensional case, we do not adaptively coarse the set of triplets.
        While this is possible and does not harm accuracy, it complicates computing a surface mesh from the set of
        triplets representing $\Gamma_{i,j}$.

        For adaptive refinement, we do not rely on a global triangulation of $\breve{S}_{i,j}$.
        Instead, for a given triplet $x \in \breve{S}_{i,j}$, let $\breve{S}_{\text{loc}}$ consist of the
        $k_{\text{near}}$ nearest triplets in $\breve{S}_{i,j}$ to $x$.
        We least-squares fit a plane $\mathcal{E}$ to $\breve{S}_{i,j}^{\text{loc}}$ as in Algorithm \ref{alg:fill3D}
        (\texttt{fill}) and create a Delaunay triangulation $\mathcal{T}$ of $\breve{S}_{i,j}^{\text{loc}}$ projected
        to $\mathcal{E}$.
        For each non-degenerated triangle in $\mathcal{T}$, we estimate the error applying \texttt{esterror} to
        $\breve{S}_{\text{loc}}^{(i,j)}$.
        If the estimated error exceeds $\varepsilon_{\text{err}}$, we employ the center of the triangle as a starting
        point for adding a new triplet.
        Looping over all $x \in \breve{S}_{i,j}$ naturally leads to many duplicate or very close starting points which
        need to eliminated before adding triples via bisection.
        However, these duplicates do not harm the efficiency of our method, as no function evaluations are required
        before computing triplets from starting points.
        We enrich $\breve{S}_{i,j}$ with the new triples created and repeat the refinement procedure at most
        $k_{\text{adap}}$ times or until no more starting points for computing triplets have been created.
    \end{Algorithm}

\begin{figure}
        \begin{center}
            \includegraphics[height=6cm]{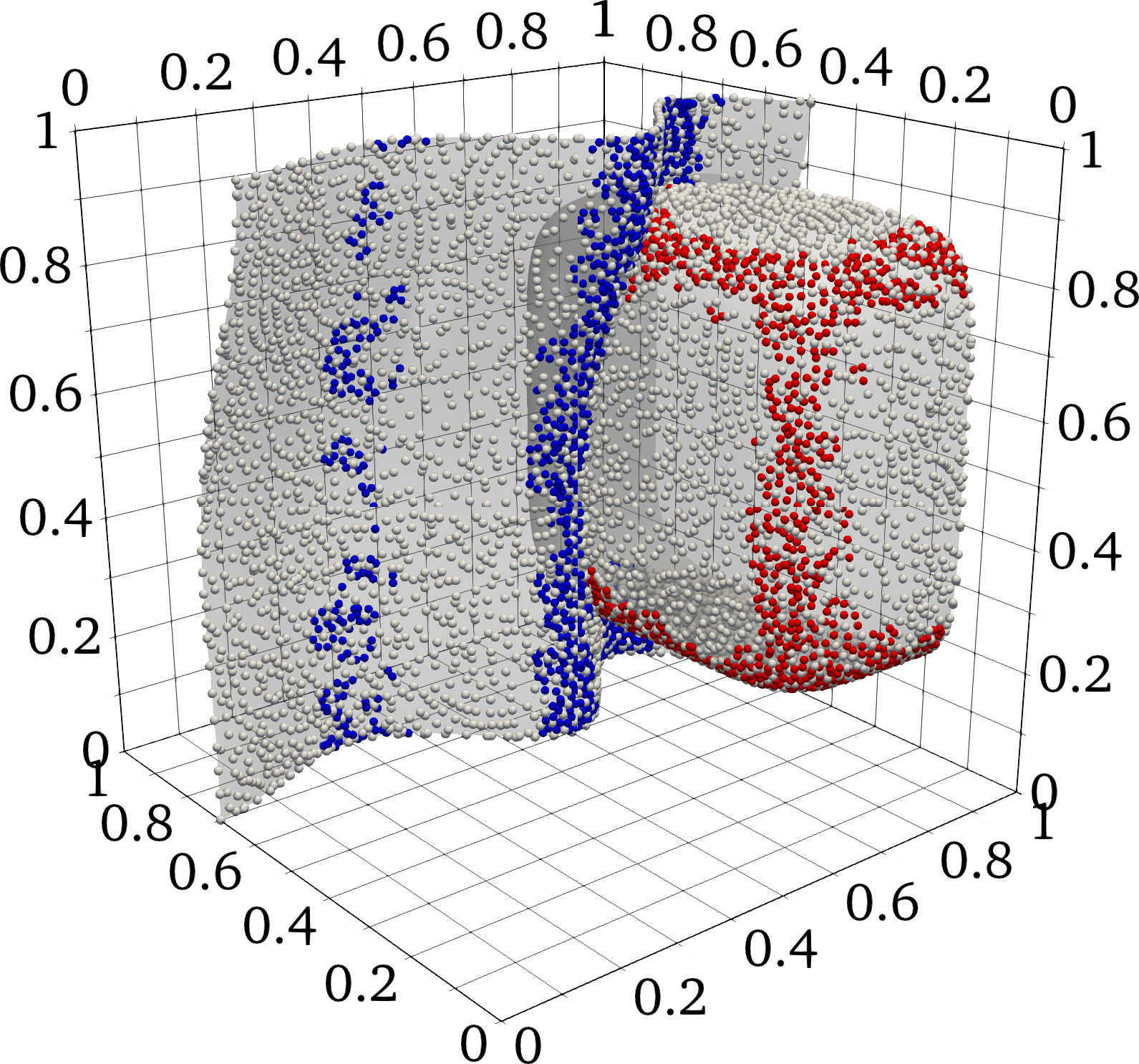}
            \caption{Final sets $\hat{S}_{1,2}^{(1)}$ and $\hat{S}_{1,3}^{(1)}$ for Test problem~\ref{prob:3D_1}.
            }\label{Fig:TestCase3D_adaptive}
        \end{center}
    \end{figure}

    \begin{Example}\label{ex:TestCase3D}
        We compute Test problem~\ref{prob:3D_1} starting with $200$ Halton-distributed points and $k_{\text{near}} = 10$.
        All other parameters are set as in the computation of Test problem~\ref{example1} in
        Section~\ref{subsec:2d_description}.
        For the corresponding numbers of triplets and the number of function evaluations, we refer to
        Table~\ref{tab:TestCase3D}.
        Figs.~\ref{Fig:TestCase3D}, \ref{Fig:TestCase3D_fill} and \ref{Fig:TestCase3D_adaptive} provide visualisations
        of the sets of triplets.
    \end{Example}
    \begin{table}
        \begin{center}
            \begin{tabular}{c|c|c|c|c|c}
            number of triplets & up to $\mathcal{M}^2$  & \texttt{iniapprox}               & \texttt{fill}                               & \texttt{expand}                           & \texttt{adapt}                              \\\hline
            approx.~of $\Gamma_{1,2}$    & &102 & 1459 & 1704 & 2460\\
            approx.~of $\Gamma_{1,3}$    & &83 & 1685 & 2058 & 2774\\
            approx.~of $\Gamma_{2,3}$   & &36 & 745 & 981 & 1510\\\hline
            function evaluations &616 & 2157 & 25510& 14232 & 5218
            \end{tabular}
            \caption{Number of triplets per $\Gamma_{i,j}$ and number of function evaluations for Example~\ref{ex:TestCase3D}.}\label{tab:TestCase3D}
        \end{center}
    \end{table}

    \section{Direct and inverse acoustic scattering problem}\label{sec:direct_acoustic_scattering}
    In this section, we briefly explain the direct and inverse scattering for the acoustic transmission problem.
    Most of the description can be found in~\cite{anachakle} and we refer the reader to~\cite{anachakle} for more details.

    Let the scatterer $D\subset \mathbb{R}^3$ be a given bounded domain with boundary $\partial D$ assumed to be of class $C^{2,\alpha}$.
    The normal unit vector $\nu$ points into the exterior $E=\mathbb{R}^3\backslash \overline{D}$ of the scatterer.
    The exterior $E$ assumed to be simply-connected is an infinite homogeneous isotropic non-absorbing acoustic medium
    which is characterised by the mass density $\varrho_e$, mean compressibility $\kappa_e$, and sound speed
    $c_e=1/\sqrt{\kappa_e\varrho_e}$.
    Likewise, the interior of $D$ is characterised by $\varrho_i$, $\kappa_i$, and $c_i=1/\sqrt{\kappa_i\varrho_i}$.
    The given scatterer is excited by a time-harmonic acoustic incident plane wave of the form
    \begin{align}
        u^\mathrm{inc}(x;\widehat{d})=\mathrm{e}^{\mathrm{i}k_e x\cdotp \widehat{d}}\,,\qquad x\in \mathbb{R}^3\,,
        \label{incidentwave}
    \end{align}
    where $k_e=\omega/c_e$ is the wave number of the acoustic wave in the host medium, $\omega>0$ the angular frequency,
    and $\widehat{d}\in \mathbb{S}^2$ the direction of incidence with $\mathbb{S}^2=\{x\in \mathbb{R}^3: \|x\|=1\}$ the
    unit sphere, where $\|\cdotp \|$ denotes the standard Euclidean norm in $\mathbb{R}^3$.
    Note that the incident field also depends on $k_e$, but this dependence is suppressed.

    The incident wave~\eqref{incidentwave} interferes with the penetrable scatterer and creates two waves.
    The first wave is the scattered field $u^{\mathrm{sca}}(x;\widehat{d})$ defined for $x\in E$ propagating outward and the second wave is the transmitted field $u^{\mathrm{int}}(x;\widehat{d})$ defined for $x\in D$.
    The total field in $E$ denoted by $u^\mathrm{ext}(x;\widehat{d})$ is the superposition of
    $u^{\mathrm{int}}(x;\widehat{d})$ and $u^{\mathrm{sca}}(x;\widehat{d})$ each of which satisfies the Helmholtz equation (the reduced wave equation) in $E$ with wave number $k_e$.
    Likewise, the transmitted field satisfies the Helmholtz equation in $D$ with wave number $k_i$.
    Precisely, we have
    \begin{align*}
        \Delta u^{\mathrm{int}}(x;\widehat{d})+k_i^2 u^{\mathrm{int}}(x;\widehat{d})=0\,,&& x\in D\,,\\
        \Delta u^{\mathrm{ext}}(x;\widehat{d})+k_e^2 u^{\mathrm{ext}}(x;\widehat{d})=0\,,&& x\in E\,.
    \end{align*}
    Due to the continuity of the acoustic pressure and the normal component of the particle velocity across $\partial D$ yields the transmission boundary conditions
    \[u^{\mathrm{int}}(x;\widehat{d})=u^{\mathrm{ext}}(x;\widehat{d})\quad\text{ and }\quad \partial_\nu u^{\mathrm{int}}(x;\widehat{d})=\tau \partial_\nu u^{\mathrm{ext}}(x;\widehat{d})\,,\qquad x\in\partial D\,,\]
    where $\tau=\varrho_i/\varrho_e >0$ is the mass density ratio of the two media.
    To ensure a well-posed boundary value problem, the scattered field
    $u^{\mathrm{sca}}(x;\widehat{d})$ needs to satisfy the Sommerfeld radiation condition
    \[\lim_{r\rightarrow \infty} r\left(\partial_r u^{\mathrm{sca}}(x;\widehat{d})-\mathrm{i}k_e u^{\mathrm{sca}}(x;\widehat{d})\right)=0\]
    with $r=\|x\|$.
    The classical acoustic transmission problem reads: find the functions $u^{\mathrm{int}}(x;\widehat{d})\in C^2(D)\cap C^1(\overline{D})$ and $u^{\mathrm{sca}}(x;\widehat{d}))\in C^2(E)\cap C^1(E)$ satisfying
    \begin{align}
        \Delta u^{\mathrm{int}}(x;\widehat{d})+k_i^2 u^{\mathrm{int}}(x;\widehat{d})=0\,,&& x\in D\,, \label{system_first_eq}\\
        \Delta u^{\mathrm{ext}}(x;\widehat{d})+k_e^2 u^{\mathrm{ext}}(x;\widehat{d})=0\,,&& x\in E\,,\\
        u^{\mathrm{int}}(x;\widehat{d})-u^{\mathrm{sca}}(x;\widehat{d})=u^{\mathrm{inc}}(x;\widehat{d})\,,&& x\in \partial D\,,\\
        \frac{1}{\tau}\partial_\nu u^{\mathrm{int}}(x;\widehat{d})- \partial_\nu u^{\mathrm{sca}}(x;\widehat{d})=\partial_\nu u^{\mathrm{inc}}(x;\widehat{d})\,,&& x\in \partial D\,,\\
        \lim_{r\rightarrow \infty} r\left(\partial_r u^{\mathrm{sca}}(x;\widehat{d})-\mathrm{i}k_e u^{\mathrm{sca}}(x;\widehat{d})\right)=0\,, && r=\|x\|
        \label{system_last_eq}
    \end{align}

    \subsection{The direct acoustic transmission problem}\label{subsec:direct_acoustic_transmission}
    Given the incident field (i.e.\ the direction of incidence $\widehat{d}$ and the wave number $k_e$), the scatterer
    $D$ (hence also its boundary $\partial D$), the wave number $k_i$, and the parameter $\tau$, one has to
    solve~\eqref{system_first_eq} --~\eqref{system_last_eq} for $u^{\mathrm{int}}(x;\widehat{d})$ and
    $u^{\mathrm{sca}}(x;\widehat{d})$.
    In the direct problem, one is only interested in the far-field $u^\infty(\widehat{x};\widehat{d})$ of
    $u^{\mathrm{sca}}(x;\widehat{d})$ which is given by
    \[u^{\mathrm{sca}}(x;\widehat{d})=\frac{\mathrm{e}^{\mathrm{i}k_e\|x\|}}{\|x\|}u^{\infty}(\widehat{x};\widehat{d})+\mathcal{O}\left(\|x\|^{-2}\right)\,,\qquad \|x\|\rightarrow \infty\]
    uniformly in all directions $\widehat{x}\in \mathbb{S}^2$.
    The far-field can be found by evaluating an integral equation over $\partial D$ given two density functions
    determined by first solving a $2\times 2$ system of integral equation of the second kind over $\partial D$
    (see~\cite[Section 4.1]{anachakle}).
    Of course, the integral equations at hand cannot be solved analytically and have to be solved numerically for
    example by the boundary element collocation method (see~\cite[Chapter 5]{kleefeld3} for more details).

    To sum up, in the direct acoustic transmission problem one is interested in $u^{\infty}(\widehat{x};\widehat{d})$
    for $\widehat{x}\in \mathbb{S}^2$ given the scatterer's boundary $\partial D$ and the direction of incidence
    $\widehat{d}\in \mathbb{S}^2$.
    The parameters $k_e$, $k_i$, and $\tau$ are given.

    \subsection{The inverse acoustic transmission problem}\label{subsec:inverse_acoustic_transmission}
    The parameters $k_e$, $k_i$, and $\tau$ are given.
    In the inverse acoustic transmission problem one tries to find/reconstruct the domain's boundary from the knowledge
    of the far-field patterns $u^{\infty}(\widehat{x};\widehat{d})$ for all $\widehat{x},\widehat{d}\in\mathbb{S}^2$.
    This can be achieved with the factorization method originally invented by Kirsch (see~\cite{Kirsch2008}).
    The theoretical justification of the factorization method for the acoustic transmission problem is given
    in~\cite[Chapter 3]{anachakle} and shown to work practically in~\cite[Chapter 4]{anachakle}.
    We briefly outline the algorithm: Assume that the far-field data are given for $\hat{x}_i$ and $\hat{d}_j$ with
    $i,j\in\{1,\ldots,m\}$ stored in the matrix $A\in \mathbb{C}^{m\times m}$.
    First, compute a singular decomposition of $A=U\Lambda V^\ast$ with $\Lambda=\mathrm{diag}(\lambda_1,\ldots,\lambda_m)$.
    For a given point $z\in \mathbb{R}^3$ compute the expansion coefficient of
    \[r_z=\left(\exp\left(-\mathrm{i}k_e z\cdotp \hat{d}_j\right)\right)_{j=1,\ldots,m}\in \mathbb{C}^m\]
    with respect to $V$ by
    \begin{equation*}
        \varrho_\ell^{(z)}=\sum_{j=1}^M V_{j,\ell}\mathrm{e}^{-\mathrm{i}k_e z\cdotp d_j}\,,\quad \ell=1,\ldots,M\,,
    \end{equation*}
    which is a matrix-vector multiplication $\varrho^{(z)}=V^\top r_z$.
    Finally, we compute
    \begin{equation*}
        W(z)=\left[\sum_{\ell=1}^M \frac{|\varrho_\ell^{(z)}|^2}{|\lambda_\ell|}\right]^{-1}
    \end{equation*}
    and plot the isosurfaces of $z\mapsto W(z)$.
    The values of $W(z)$ should be much smaller for $z\notin D$ than those lying within $D$.
    The threshold value can be approximately determined by a scatterer such as the unit sphere and then reused for other scatterers as well.
    Note that until now, an equidistant set of points $N$ within a predefined box have been used to find the values of
    $W(z)$ leading to an amount of $N\times N\times N$ function evaluations for such a ``sampling'' method
    (see also~\cite{bazan} for other sampling methods).
    This can be considerably reduced as shown in the next section.

    \section{Applications}\label{sec:applications}

    \subsection{Test cases related to the detection of faults}\label{subsec:2dtests}
\begin{figure}
        \begin{center}
            \includegraphics[height=3cm]{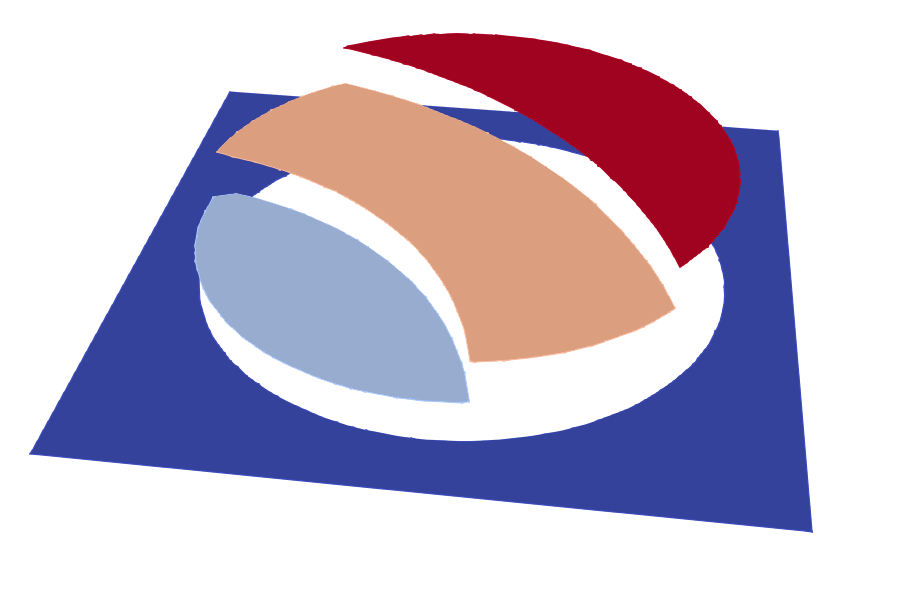}
            \includegraphics[height=3cm]{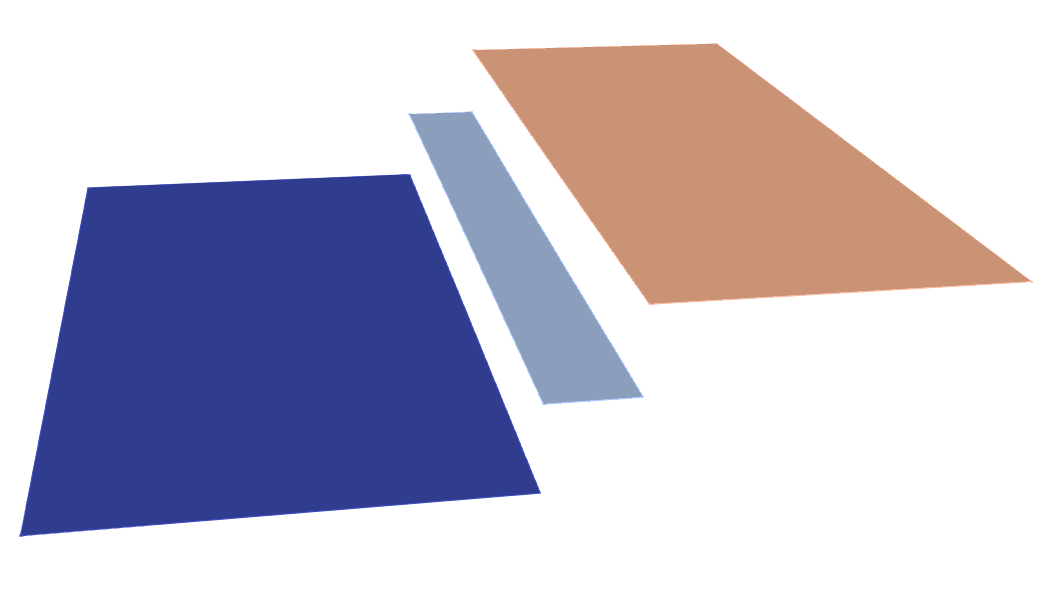}
            \includegraphics[height=3cm]{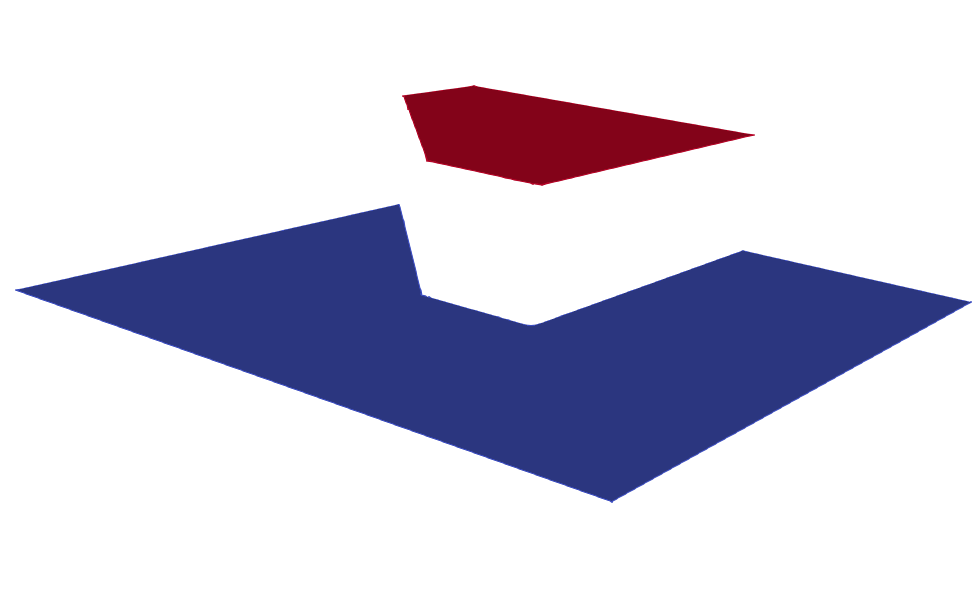}
            \caption{Reconstructed subdivisions for Test problems~\ref{testprob_Allasia01},~\ref{testprob_Allasia04},
            and~\ref{testprob_Allasia05}.}\label{fig_Allasia_final}
        \end{center}
    \end{figure}

    In this section, we consider test cases which have been defined in the context of the detection of fault lines or
    fault surfaces.
    \begin{Test problem}
        \label{testprob_Allasia01}
        Following Allasia et al.\ \cite{Allasia2010EfficientAA} and Gutzmer et at.\ \cite{Gutzmer.1997}, we set
        $x_M = (0.5, 0.5)^{\top}$, $x = (x_1, x_2)^{\top}$, $\Omega = [0,1]^2$ and consider the function
        \begin{equation}
            \label{eq:testfunc_allasia01}
            f(x_1, x_2) =
            \left\{
            \begin{array}{ccl}
                1 + 2 \lfloor 3.5 \|(x_1, x_2)\|_2 \rfloor & , & \|(x - x_M) \|_2 < 0.4 \\
                0                                          & , & \text{otherwise}
            \end{array}
            \right.
        \end{equation}
        This function is piecewise constant with smooth fault lines, it holds
        \begin{align*}
            \Gamma = \bigcup \Gamma_{i,j} = & \hspace{0.35cm}\left\{ x \in \mathbb{R}^2 \: : \: \| x - x_m\|_2 < 0.4\right\} \\
            &\cup \left\{ \| x\|^2_2 = 4/7 \: : \: \| x - x_m \|_2 < 0.4\: \right\} \\
            &\cup \left\{ \| x\|^2_2 = 6/7 \: : \| x - x_m \|_2 < 0.4\: \right\} \: .
        \end{align*}
        The
        set $\Gamma_{1,3}$ consists of two separate components (Fig.\ \ref{fig_Allasia_final}, left).
    \end{Test problem}
    \begin{Test problem}
        \label{testprob_Allasia04}
        We set
        \begin{equation}
            \label{eq:testfunc_allasia04}
            f(x_1, x_2) =
            \left\{
            \begin{array}{rcl}
                3 & , & x_1 > 0.6 \\
                1 & , & x_1 < 0.5 \\
                2 & , & \text{otherwise}
            \end{array}
            \right. \: .
        \end{equation}
        The set $\Gamma$ consists of two straight lines and coincides with the two fault lines
        from~\cite{Allasia2010EfficientAA}, Example 4 (Fig.\ \ref{fig_Allasia_final}, middle).
    \end{Test problem}
    \begin{Test problem}
        \label{testprob_Allasia05}
        We set
        \begin{equation}
            \label{eq:testfunc_allasia05}
            f(x_1, x_2) =
            \left\{
            \begin{array}{rcl}
                2 & , & x_1 > 0.4 \wedge x_2 > 0.4 \wedge x_2 < 0.2 + x_1 \\
                1 & , & \text{otherwise}
            \end{array}
            \right.
        \end{equation}
        The set $\Gamma_{1,2}$ coincides with the fault line from~\cite{Allasia2010EfficientAA}, Example 5
        (Fig.\ \ref{fig_Allasia_final}, right).
    \end{Test problem}
    \begin{figure}
        \begin{center}
            \includegraphics[height=5cm]{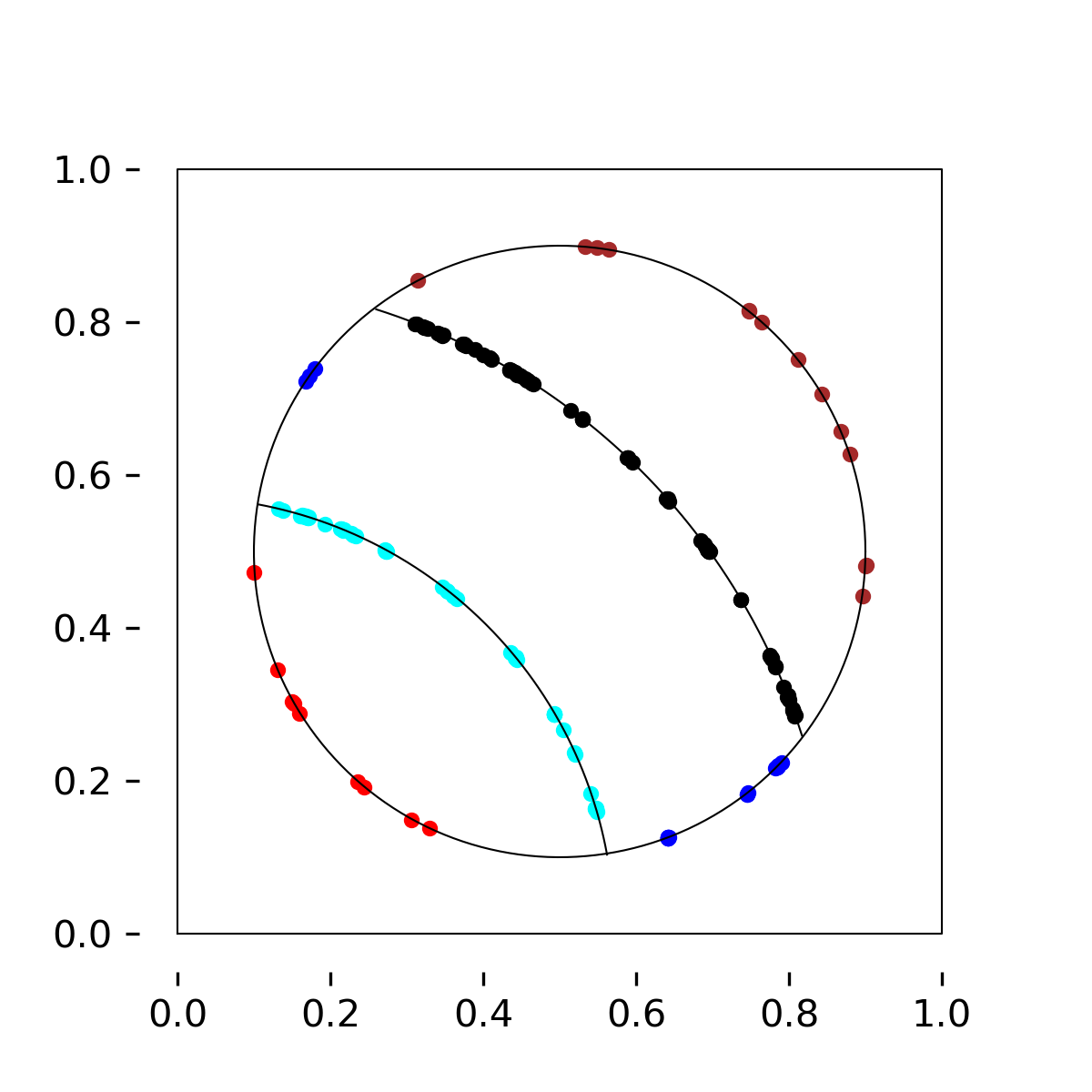}
            \includegraphics[height=5cm]{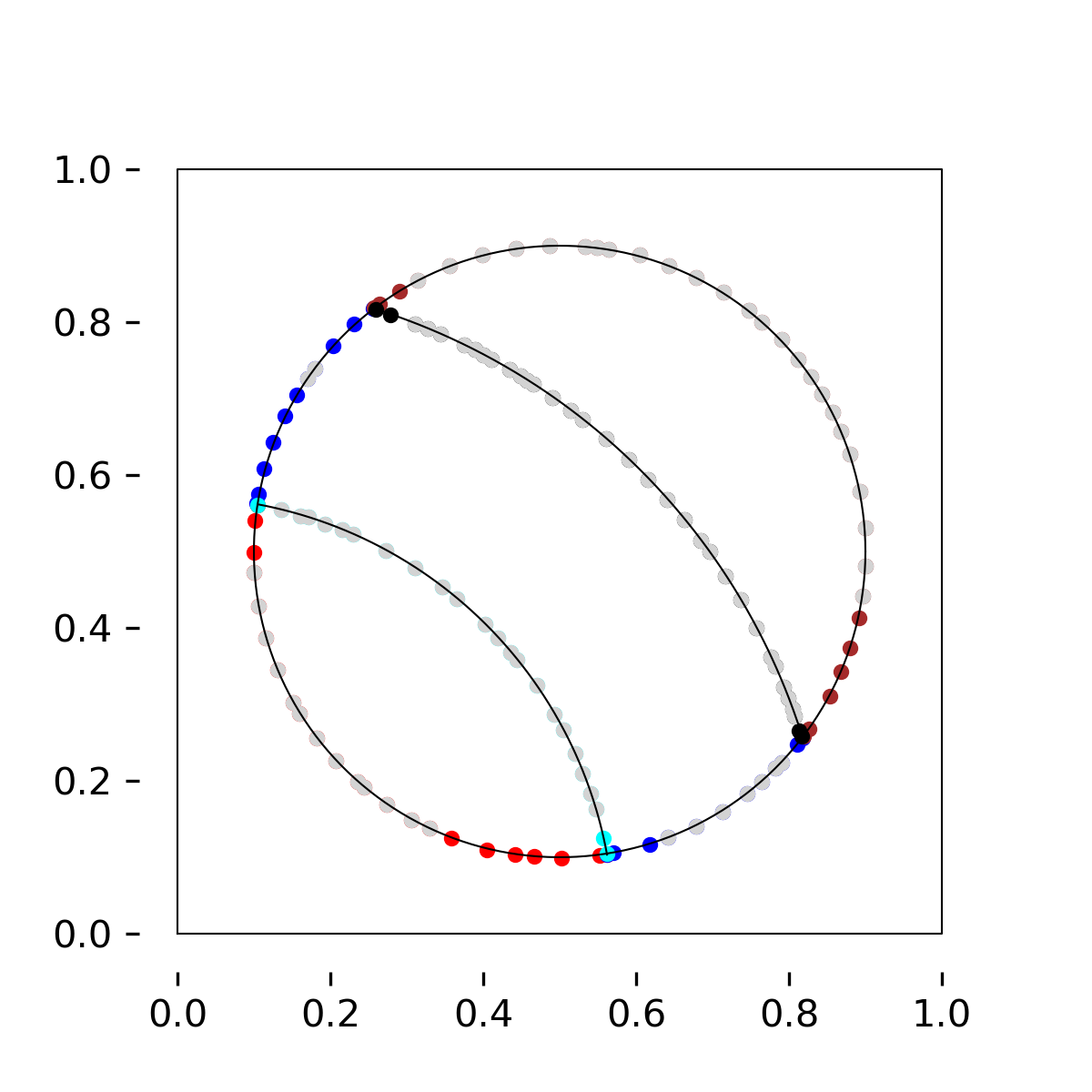}
            \includegraphics[height=5cm]{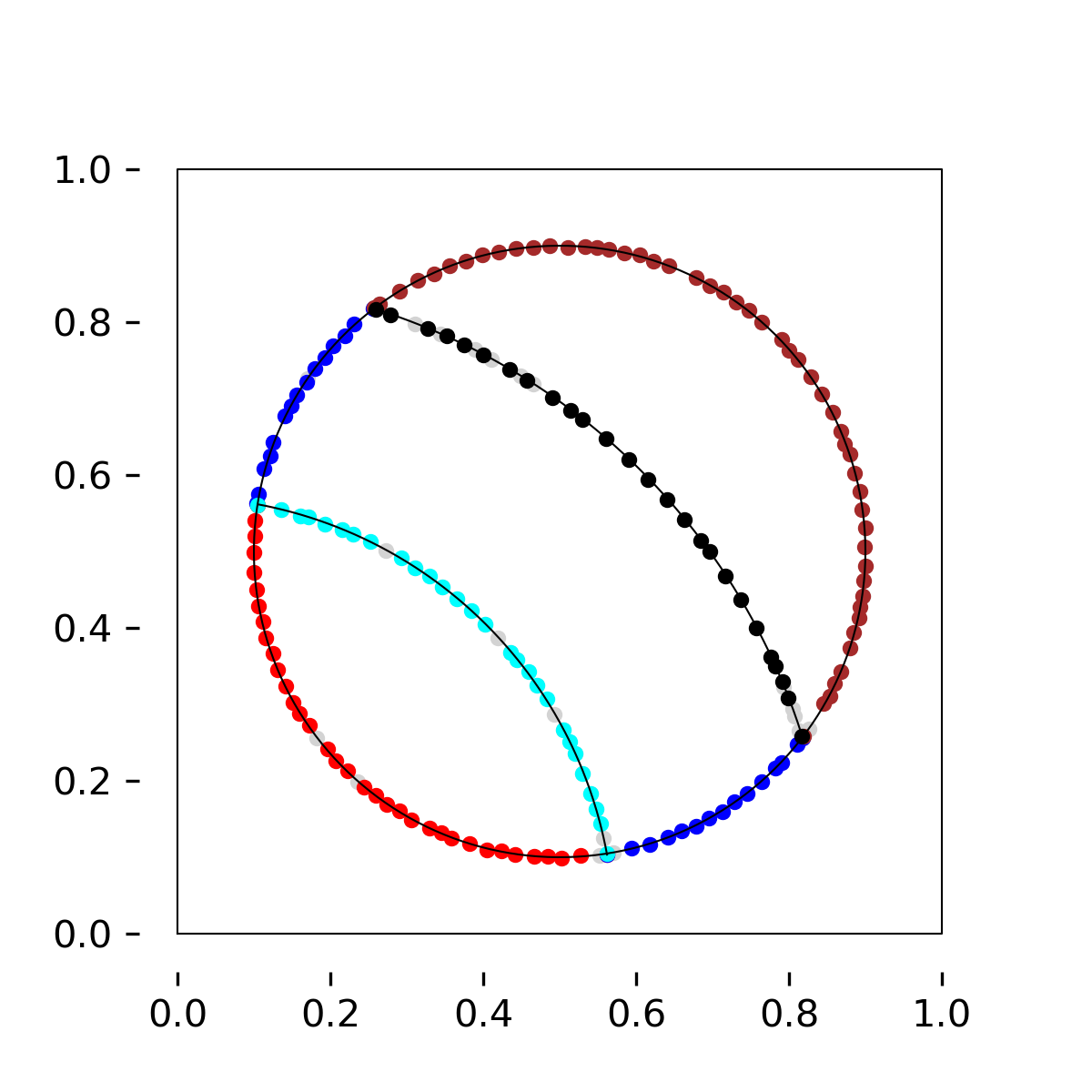}
            \caption{$S^{(i)}_{i,j}$ (left), $\breve{S}_{i,j}^{(i)}$ (middle), and
                $\hat{S}_{i,j}^{(i)}$ (right) for Test problem~\ref{testprob_Allasia01}.
                As for the computations of Example~\ref{example1} in Section~\ref{subsec:2d_description}, previously
                existing points are greyed out.}\label{fig_allasia1}
        \end{center}
    \end{figure}

    \begin{figure}
        \begin{center}
            \includegraphics[height=5cm]{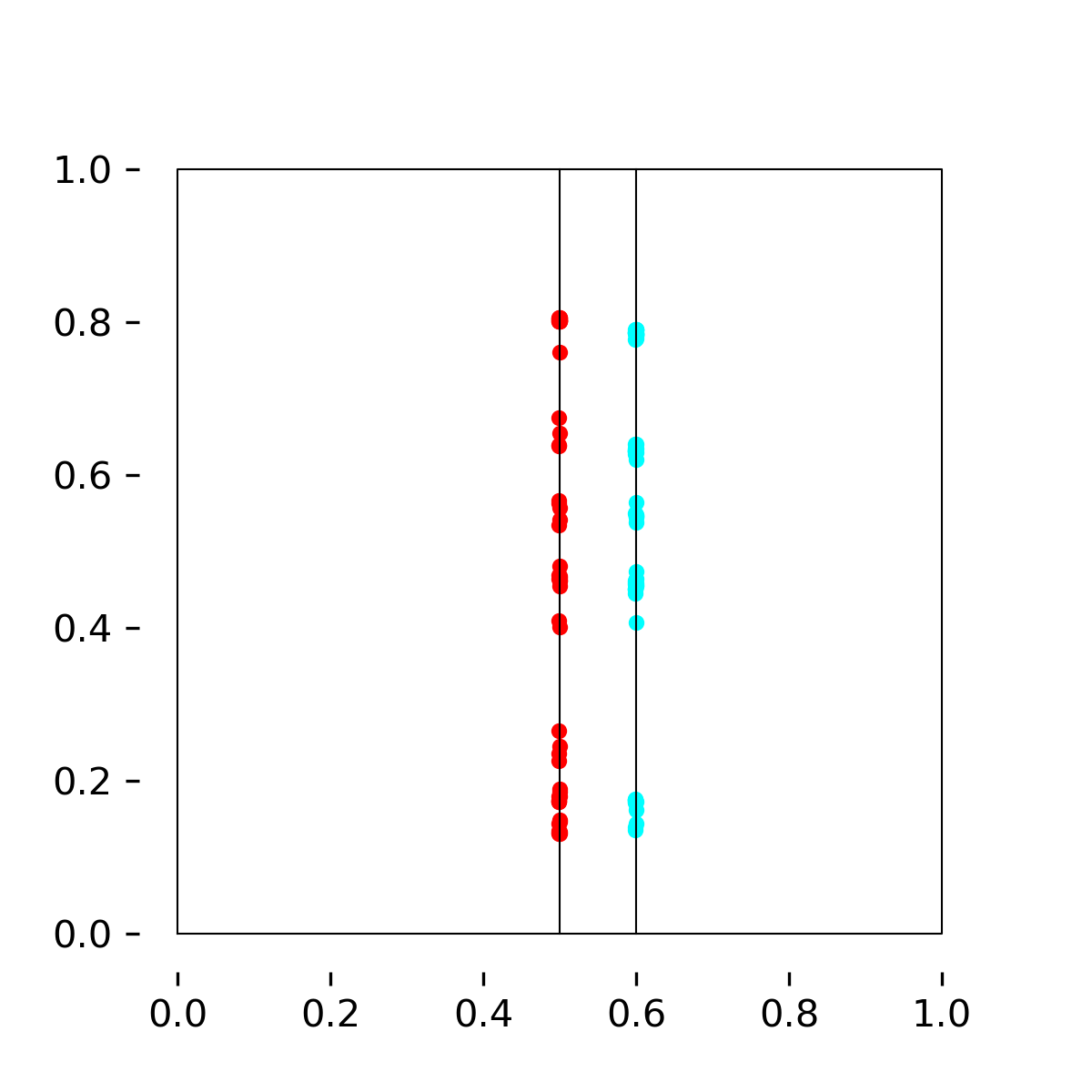}
            \includegraphics[height=5cm]{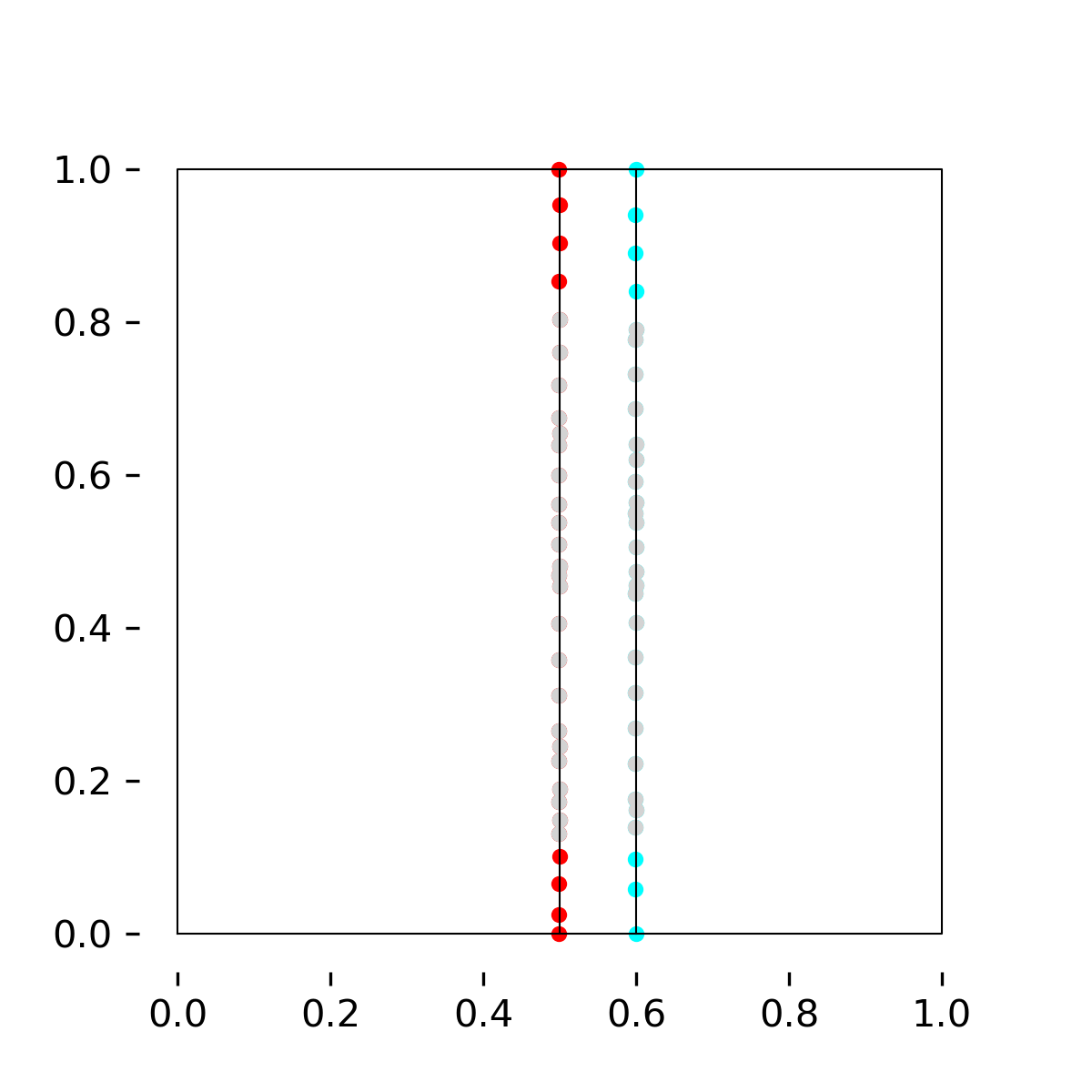}
            \includegraphics[height=5cm]{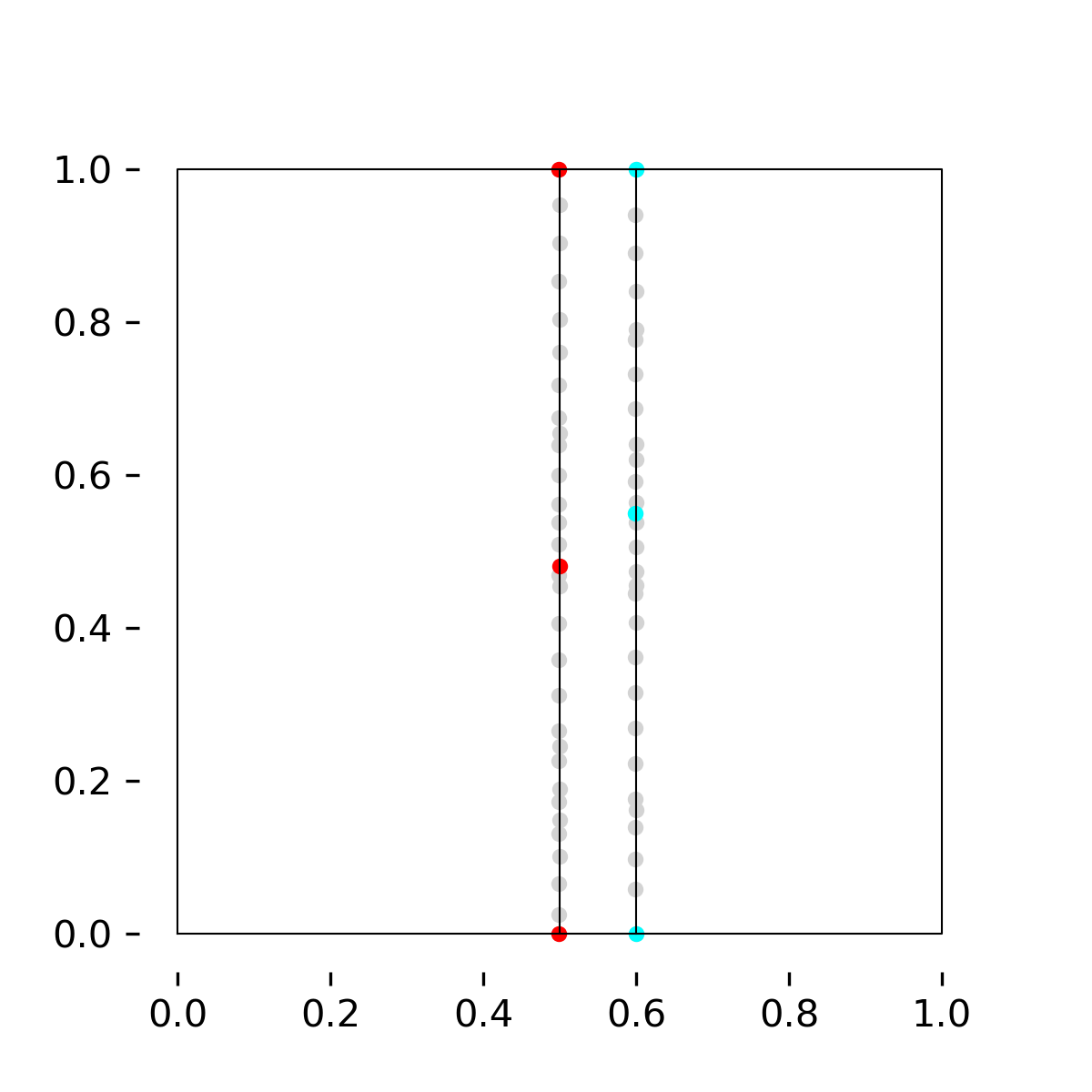}
            \caption{$S^{(i)}_{i,j}$ (left), $\breve{S}_{i,j}^{(i)}$ (middle), and
                $\hat{S}_{i,j}^{(i)}$ (right) for Test problem~\ref{testprob_Allasia04}.}\label{fig_allasia4}
        \end{center}
    \end{figure}

    \begin{figure}
        \begin{center}
            \includegraphics[height=5cm]{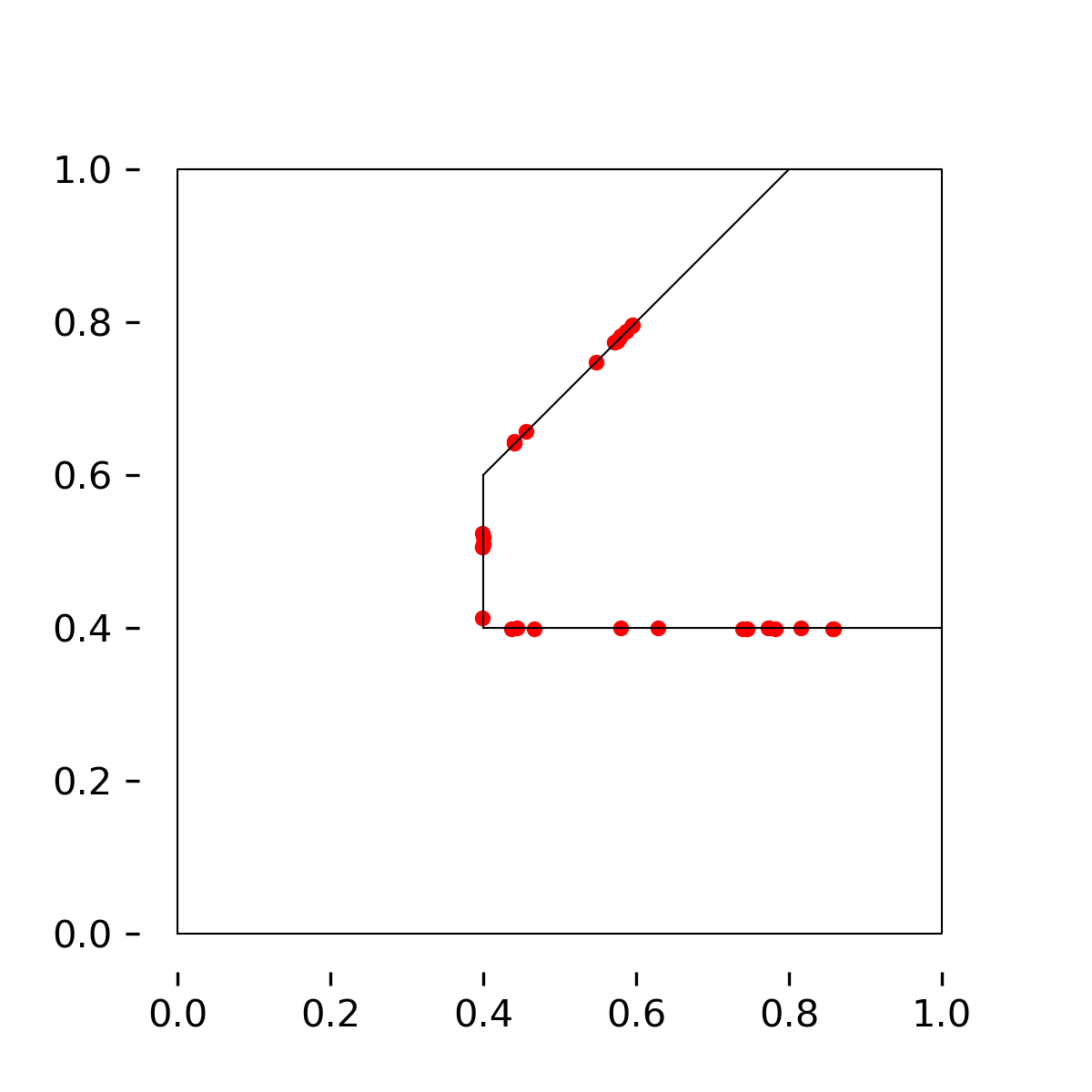}
            \includegraphics[height=5cm]{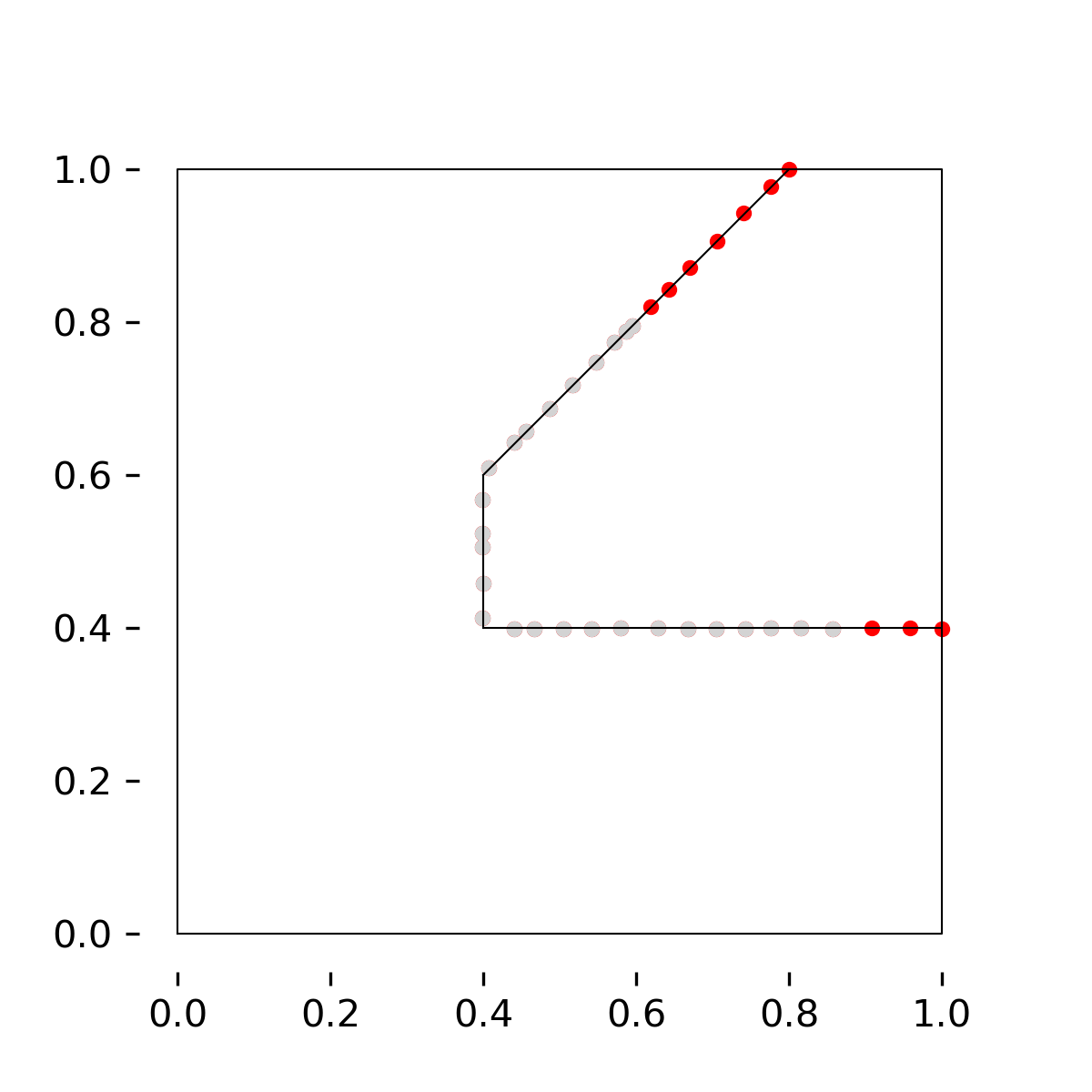}
            \includegraphics[height=5cm]{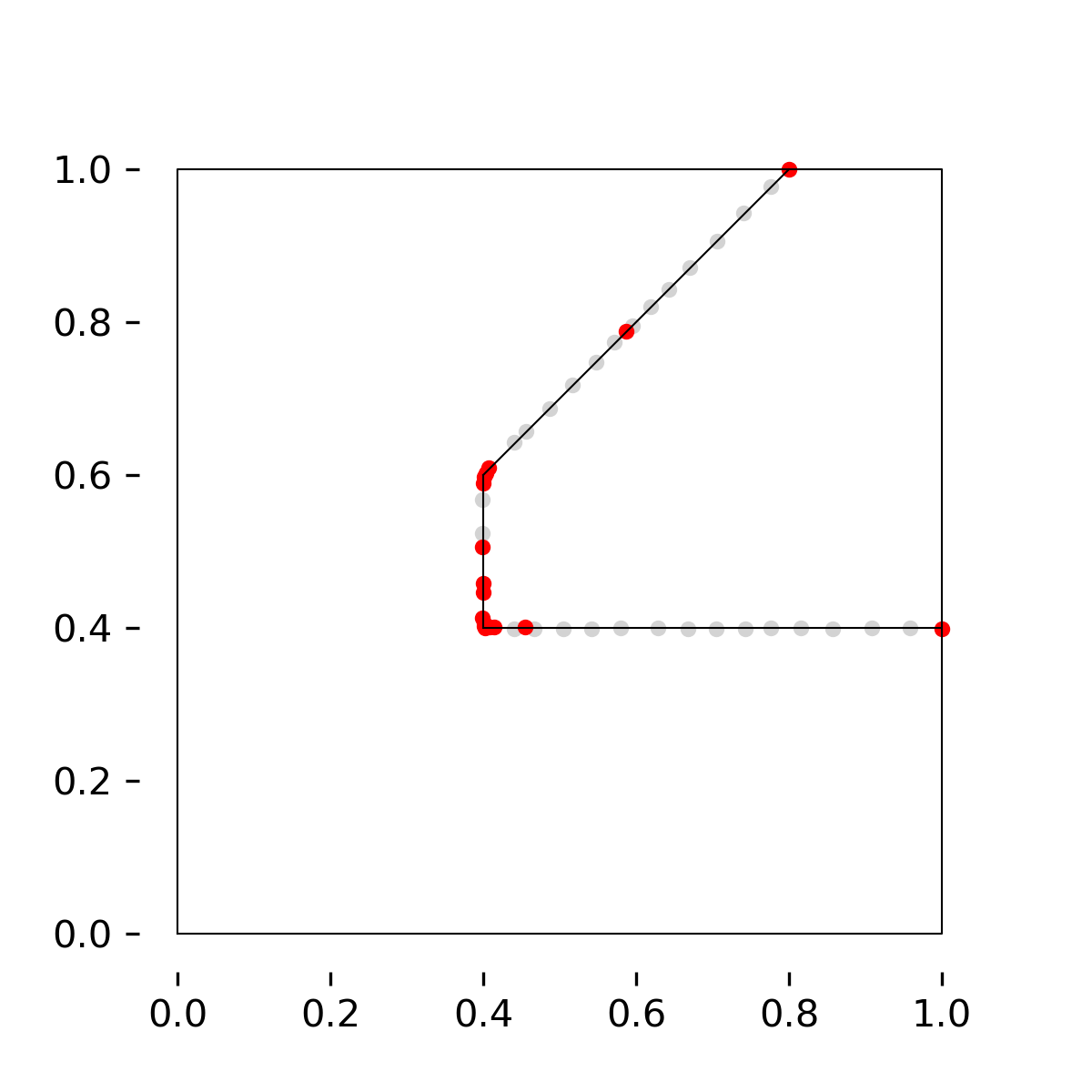}
            \caption{$S^{(1)}_{1,2}$ (left), $\breve{S}_{1,2}^{(1)}$ (middle), and
                $\hat{S}_{1,2}^{(1)}$ (right) for Test problem~\ref{testprob_Allasia05}.}\label{fig_allasia5}
        \end{center}
    \end{figure}

    We compute all three test problems starting with $X$ from Example~\ref{example1} and the
    same parameters used there and in the subsequent examples in Section~\ref{subsec:2d_description}.
    Our results (Figs.\ \ref{fig_allasia1},~\ref{fig_allasia4}, and~\ref{fig_allasia5}) demonstrate a successful approximation
    of all $\Gamma_{i,j}$.
    We provide details in Tables~\ref{tab_results} and~\ref{tab_results_npoints} combined with the results
    for Test problem~\ref{TestProb1} from the preceding section.
    Algorithm \texttt{iniapprox} is the most demanding Building block in terms of function evaluations in all examples considered.
    Averaged over Test problems~\ref{TestProb1},~\ref{testprob_Allasia04}, and~\ref{testprob_Allasia05}, it takes on
    average 4.0 classifications for adding one triplet in \texttt{fill}.
    However, it takes 12.2 classifications per triplet added for Test problem~\ref{testprob_Allasia01}.
    This is due to the fact that $\Gamma_{1,3}$ consists of two components, which is detected by a failed attempt of
    filling the large gap between these two components (Remark~\ref{rem_several_components}).
    This process adds classifications, but no triplets.
    In average over the four examples, \texttt{expand} requires $5.5$ classifications per triplet added.
    A significant part of the classifications is required for the very first and very last step, as finding them
    requires bisection on the extrapolation curve $\gamma$, which involves at least one classification per iteration
    step.
    Considering the same ratio for \texttt{adapt} would be misleading, as triplets are added and removed.

    Allasia et al.~present numbers of function evaluations for their method of surface reconstruction applied to these test
    problems in \cite{Allasia2010EfficientAA}.
    They report $5185$ evaluations of $f$ for Test problem~\ref{testprob_Allasia01} (our method: $2213$),
    $3479$ for Test Problem~\ref{testprob_Allasia04} (our method: $871$), and $2099$ for Test
    problem~\ref{testprob_Allasia05} (our method: $450$).
    It turns out that the number of function evaluations is significantly higher than for our
    algorithm.
    However, this comparison is limited by several factors.
    In contrast to us, Allasia et al.~need to employ an explicit classification algorithm, as they consider more general
    functions $f$ than we do, which may infer additional function evaluations.
    The number of function evaluations required heavily depends on the accuracy of the approximation and resolution of
    the fault line, such that comparing the number of function evaluations given a prescribed
    tolerance for both algorithms would provide a more sound basis for comparing the efficiency of the two algorithms.
    Detailed information on accuracy is however lacking in Test problem~\ref{testprob_Allasia01}.
    \begin{table}
        \begin{center}
            \begin{tabular}{c|c|c|c|c|c}
                Test P.                  & up to $\mathcal{M}^2$                       & \texttt{iniapprox}               & \texttt{fill}                               & \texttt{expand}                           & \texttt{adapt}                              \\\hline
                \ref{TestProb1}          & 146     & 458     & 94 & 185  & 205  \\
                \ref{testprob_Allasia01} & 254      & 1017      & 465 & 350  & 127\\
                \ref{testprob_Allasia04} & 176  & 595  & 26 & 74  & 0\\
                \ref{testprob_Allasia05} & 87 & 202 & 34 & 44  & 83  \\
            \end{tabular}
            \caption{Number of function evaluations per Building block for reconstructing the subdomains.}\label{tab_results}
        \end{center}
    \end{table}
    \begin{table}
        \begin{center}
            \begin{tabular}{c|c|c|c|c}
                Test P.                  & \texttt{iniapprox}         & \texttt{fill}        & \texttt{expand}      & \texttt{adapt}       \\\hline
                \ref{TestProb1}          & 40     & 55 & 89  & 125  \\
                \ref{testprob_Allasia01} & 66    & 98 & 138  & 171\\
                \ref{testprob_Allasia04} & 32  & 45 & 60  & 6\\
                \ref{testprob_Allasia05} & 17 & 26 & 36  & 18  \\
            \end{tabular}
            \caption{Total number of triplets after the respective Building block.}\label{tab_results_npoints}
        \end{center}
    \end{table}

    \subsection{Decision analysis and modeling}\label{subsec:DecisionAnalysis}
    MCDA methods depend on several parameters called ``input factors''.
    Almost all common MCDA models require to set up a \emph{performance matrix} $P = (p_{i,j})$, where $p_{i,j}$ encodes
    the objective benefit or cost of the $i$-th alternative with respect to the $j$-th criterion.
    Both finding suitable criteria for modeling the decision process and setting up $P$ requires expertise
    (e.g.~\cite{Voegele2020} among many others) and is beyond the scope of this work.
    Instead, we assume that $P$ is given and exactly known, although this assumption may be questioned in many
    practical applications.
    In addition to $P$, the decision maker needs to provide non-negative \emph{weights} $w = (w_1, \hdots w_c)^{\top}$ which
    reflect the importance of the criteria from his point of view.
    Although originally intended as a decision support tool, MCDA methods have been used for some time for decision
    analysis, e.g.~to predict decisions of actors under changed framework conditions, modelled by a changed performance matrix.
    In the context of decision analysis, $w$ is in many practical applications not exactly known and hard to obtain,
    as e.g.~surveys are time-consuming and prone to bias due to socially accepted answers and other effects.
    Moreover, the restriction to fixed weights ignores possible diversity within actors.
    All of the above motivates a robustness analysis of the decision predicted, which in many cases focuses on
    examining the robustness of the decision with respect to perturbations in $w$.

    In what follows, we apply our algorithm for computing $\Gamma_{i,j}$ to such a robustness analysis and consider
    as example an application from the scientific monitoring of the mobility turnaround in Germany.
    Ball et al.~\cite{BALL2021120925} investigate car users' attitudes towards the purchase of hybrid (HEV) and
    electric vehicles (BEV) versus conventional cars with internal combustion engine (ICE).
    For this purpose, they identify 13 criteria and divide them into 5 categories (``Ecological'',``Economic'',
    ``Social'', ``Comfort'', ``Other'').
    The authors weight both the criteria within the categories and the categories themselves.
    Taking the former weights as given, we obtain the performance matrix in Table~\ref{tab:Perf_cars} from the data
    in~\cite{BALL2021120925}.
    In contrast to SAW (Simple Additive Weighting)~\cite{SAW} as in~\cite{BALL2021120925}, we use the more complex MCDA
    method SIR-TOPSIS~\cite{Xu.2001} here, which includes the very widespread MCDA methods Promethee
    II~\cite{J.P.Brans.1986} and SAW as special cases.
    We omit all details on how TOPSIS works for the sake of brevity and refer instead to~\cite{Xu.2001}.

    \begin{table}
        \begin{center}
      \begin{tabular}{l|l|l|l|l|l}
         & Ecological & Economic & Social & Comfort & Other\\ \hline
         BEV & 0.5025 & 0.2792 & 0.6250 & 0.1497 & 0.1342 \\
         ICE & 0.1256 & 0.4167 & 0.1250 & 0.4300 & 0.6710\\
         HEV & 0.3719 & 0.3042 & 0.2500 & 0.4202 & 0.1948\\\hline\hline
          weightings & 2 & 7 & 0.1 & 8 & 0
      \end{tabular}
      \end{center}
      \caption{Performance matrix $\tilde{P}$ with corresponding criteria generated from~\cite{BALL2021120925}, and
      weightings from the car users' point of view prior to normalisation (last row).}\label{tab:Perf_cars}
    \end{table}

\begin{figure}
        \begin{center}
            \includegraphics[height=5cm]{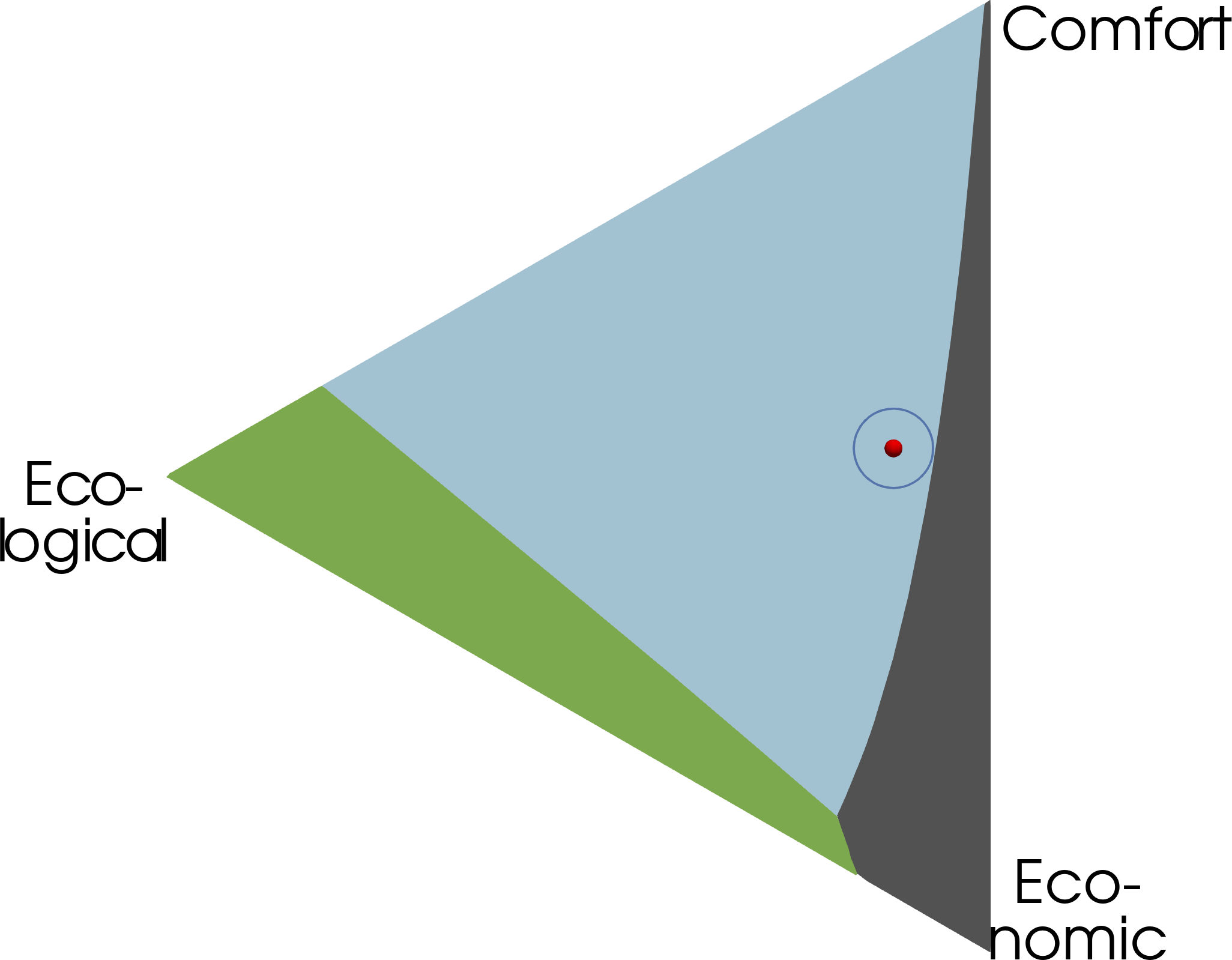}\hspace{1cm}
            \includegraphics[height=5.5cm]{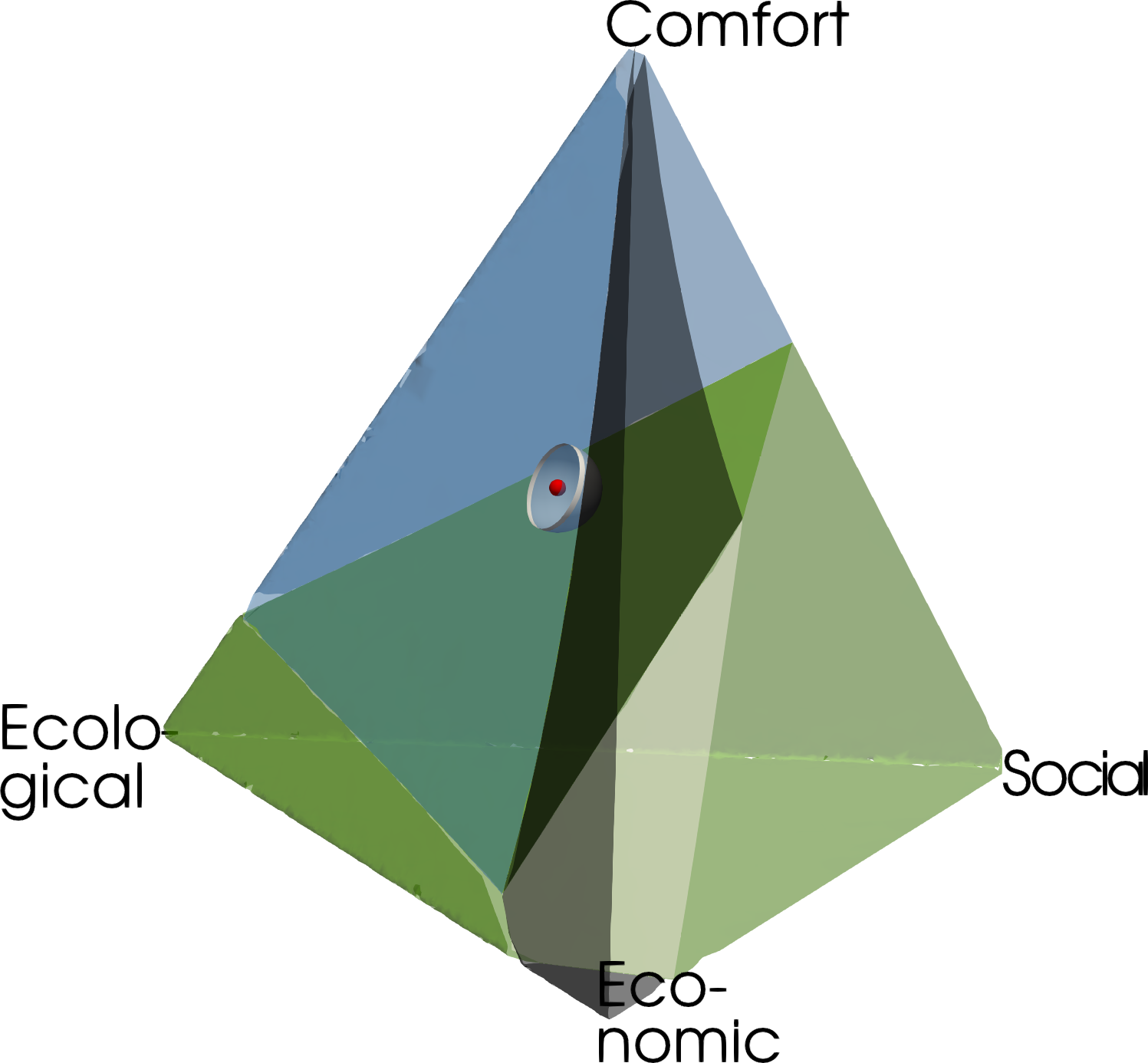}
        \end{center}
        \caption{Visualisation of the decision space for $3$ (left) and $4$ (right) variable criteria. The
        car users' weightings according to \cite{BALL2021120925} are displayed as red dot along with the circumsphere
        with radius $\rho = 0.06$.}\label{fig:MCDA}
        \end{figure}

    SIR-TOPSIS requires as many other MCDA methods that the non-negative weights are normalized, i.e. $\sum_{i=1}^n w_i = 1$.
    Therefore, the set of normalized admissible weights is the standard simplex in $\mathbb{R}^n$.
    For visualisation, we consider $m=3$ or $m=4$ of the weights to be variable, and the rest to be fixed.
    Let be $w = w_v + w_f$, where $w_v$ consists of the variable weights and zeroes elsewhere, and $w_f$ of the fixed
    ones, correspondingly.
    Therefore, normalisation implies $\sum w_{i,v} = 1 - \sum w_{i,f} := c_f$ such that the set of variable weights
    corresponds to a downscaled standard simplex in $\mathbb{R}^m$, which we embed in $\mathbb{R}^{m-1}$ by appropriate
    translation and rotation.
    This yields the equilateral triangle (for $m=3$) and the regular tetrahedron (for $m=4$) shown in Fig.~\ref{fig:MCDA}.

    For a 2D-visualisation of the decision space, we consider the weights corresponding to ``Ecological'', ``Economic'' and
    ``Comfort'' variable (Fig.~\ref{fig:MCDA}) and colour all weights leading to a decision in favour of ICEs black, for
    HEVs blue and for BEVs green.
    We divide the weights $\tilde{w}$ proposed in~\cite{BALL2021120925} as a representation of car user's
    mindset in 2020 in a fixed and variable part, setting $\tilde{w} = \tilde{w}_v + \tilde{w}_f$ and display
    $\tilde{w}_v$ as bright red dot in Fig.~\ref{fig:MCDA}.
    In contrast to~\cite{BALL2021120925}, SIR-TOPSIS seems to predict a shift to HEVs under today's conditions.
    For measuring robustness, we consider the largest sphere around $\tilde{w}_v$ still fully contained in the set of
    weightings leading to HEVs $W_{\text{HEV}}$ and propose its radius $\rho$ as simple measure of robustness.
    As we have a polygonal approximation of $W_{\mathrm{HEV}}$ at hand, iteratively approximating $\rho$ boils down to an
    intersection test of polygons, if we approximate the circle by a sufficiently fine polygon.
    We end up with $\rho \approx 0.06$, which reconciles the findings of Ball et al.~and ours: As $\rho$ is rather
    small, the decision in favour of HEVs is not very robust, as even small perturbations of the weightings may lead to
    a decision in favour of ICEs. In~\cite{BALL2021120925}, the decision towards ICEs was found not be very robust.
    However, both ours and Ball et al.'s results indicate that a broad shift towards BEVs is unlikely to happen under current
    circumstances.

\begin{figure}
        \begin{center}
        \includegraphics[height=5cm]{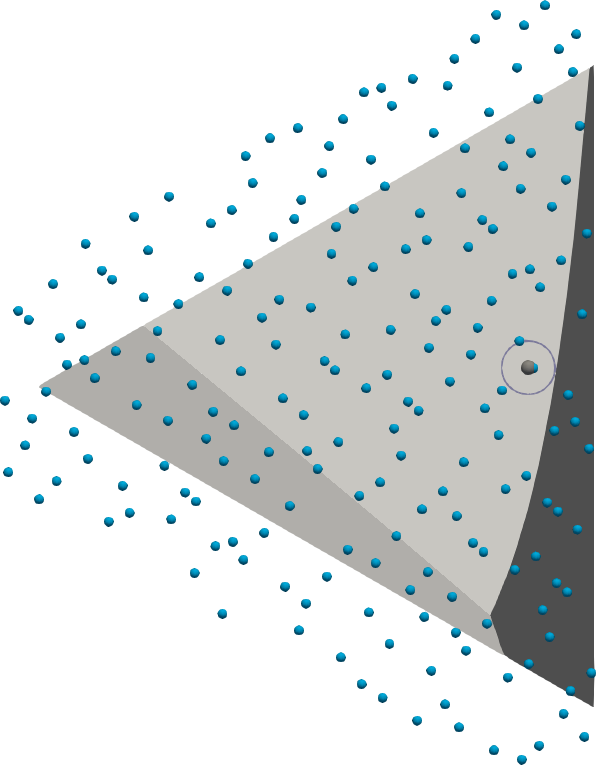}\hspace{0.5cm}
        \includegraphics[height=4.5cm]{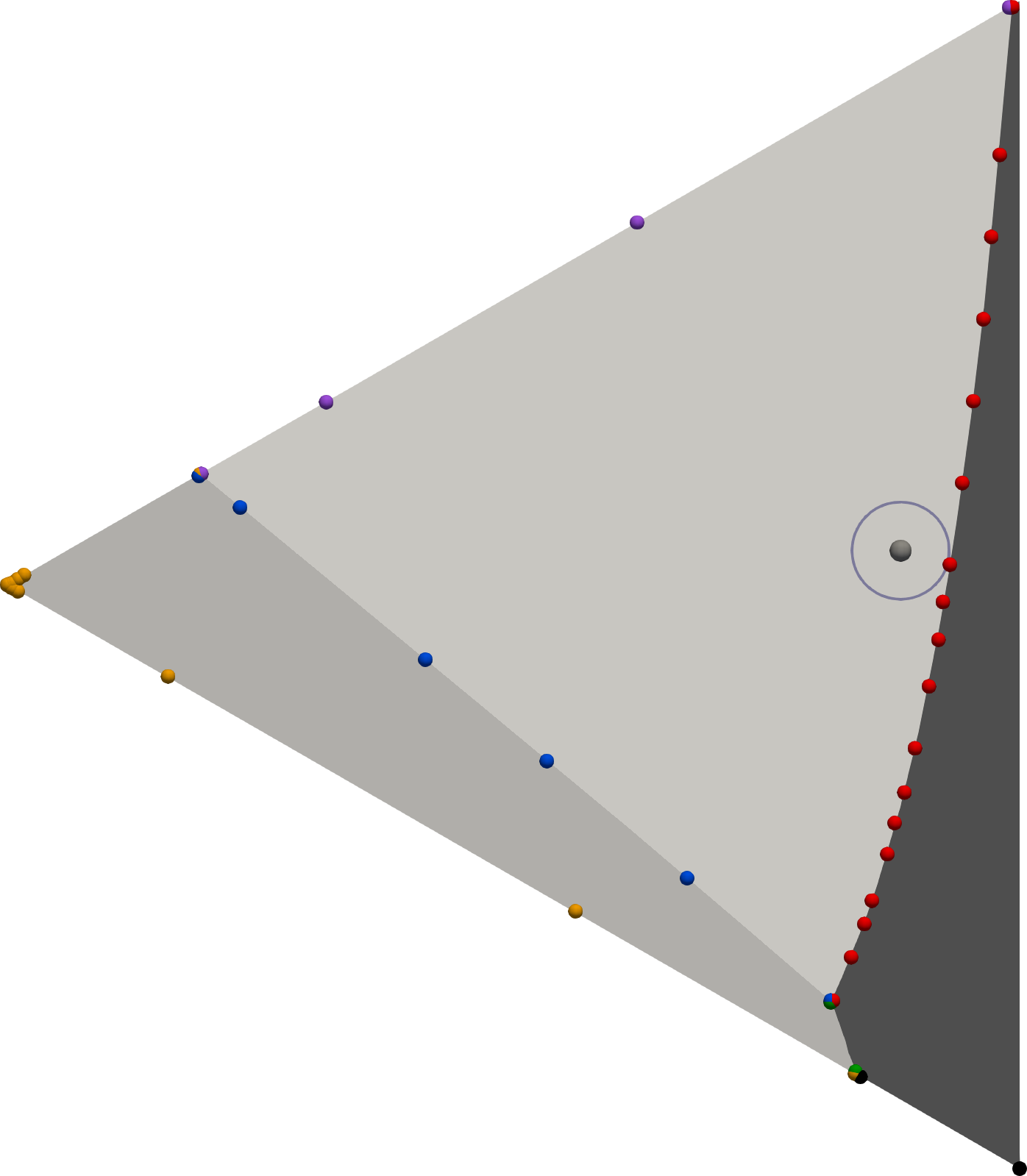}\hspace{1cm}
        \includegraphics[height=5.5cm]{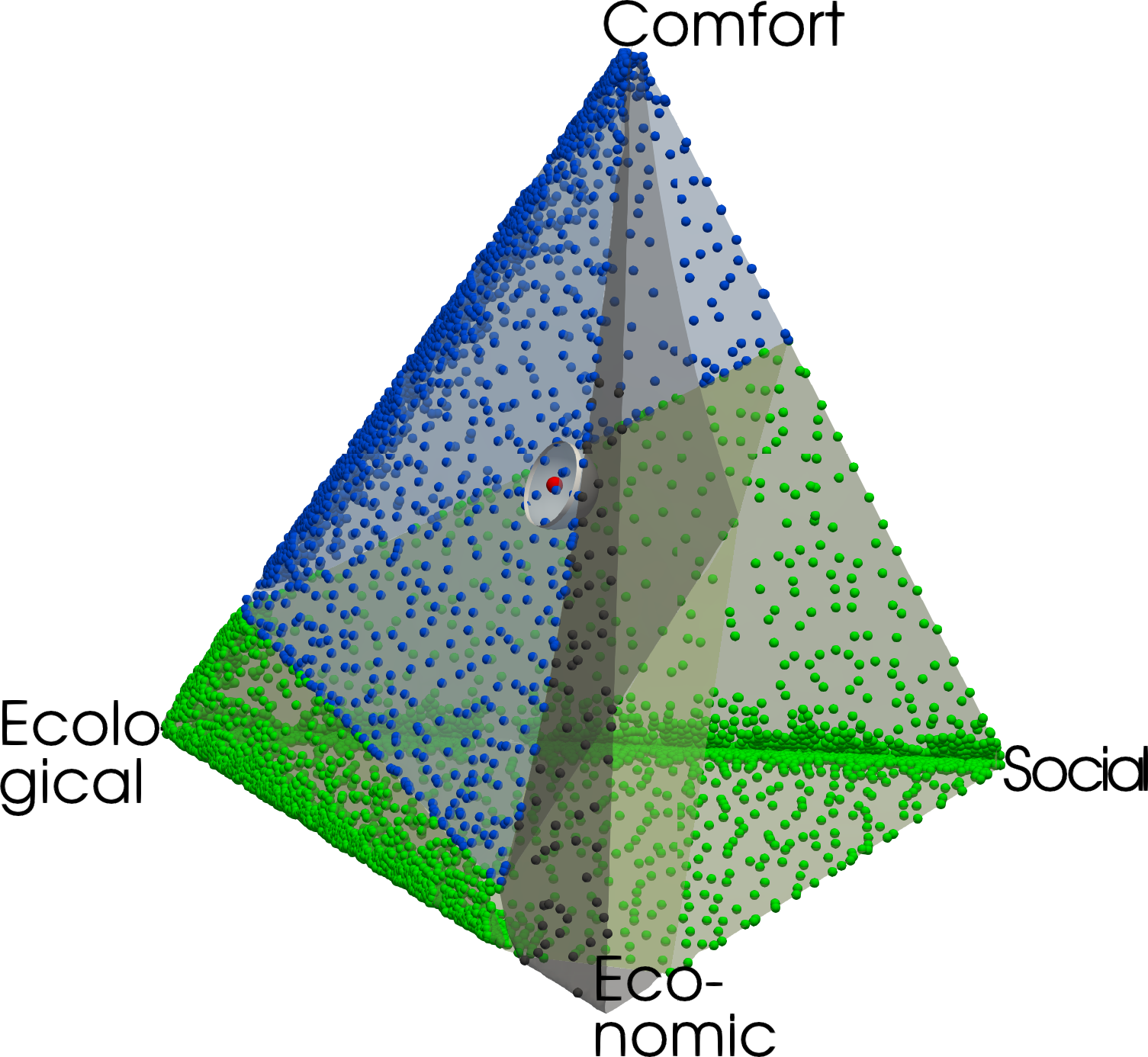}
            \end{center}
        \caption{Initial point set $X$ (left) and final sets $\hat{S}_{i,j}$ (middle) for $m=3$. Right: Selected
        final sets $\hat{S}_{i,j}$ for $m=4$.}\label{fig:MCDA_2D_points}
        \end{figure}

    For $m=3$ and $m=4$, we have $c_f = 0.9942$.
    The downscaled standard simplex is rotated and translated to the
    equilateral triangle with vertices $c_f(0.4082,-0.7071)^{\top}$, $c_f(0.4082,0.7071)$, and $c_f(-0.8165, 0)$.
    Therefore, we set $\Omega = c_f[-0.9, 0.4082] \times c_f[-0.9,0.9]$ and employ an initial point set $X$ consisting of
    $100$ Halton-distributed points, where points far away from the triangle have been discarded, as they cannot aid
    approximating $\Gamma_{i,j}$ (Fig.~\ref{fig:MCDA_2D_points}, left).
    We choose the values of all algorithm-related parameters as for the examples in Sections~\ref{subsec:2d_description}
    and~\ref{subsec:2dtests}, except of $\varepsilon_{\mathrm{gap}} = c_f\cdot 0.05$.
    For $m=4$, we start with $X$ consisting of $500$ Halton-distributed points in the vicinity of the tetrahedron and
    take all values for the parameters of the algorithm from Section~\ref{subsec:approximation_in_3D}, except of
    $k_{\text{adap}} = 3$.
    For the final sets $\hat{S}_{i,j}$, we refer to Fig.~\ref{fig:MCDA_2D_points} (right and middle) and display the
    number of function evaluations in Table~\ref{tab_resultsTOPSIS2D}.
    These sets have been used for approximating $\rho$; Fig.~\ref{fig:MCDA} was generated with the proposed
    algorithm, albeit using a finer initial point set and modified parameter settings, as analytical descriptions of
    the decision curves and surfaces, resp., are unknown.
    \begin{table}
        \begin{center}
            \begin{tabular}{c|c|c|c|c|c}
                                  & up to $\mathcal{M}^2$                       & \texttt{iniapprox}               & \texttt{fill}                               & \texttt{expand}                           & \texttt{adapt}                              \\\hline
                evaluations, $m=3$& 397     & 1325     & 198 & 283  & 44  \\
                evaluations, $m=4$& 1518     & 4553     & 13352 & 14752  & 24344  \\
                no.\ of triplets, $m=3$ & &100      & 145 & 174  & 45\\
                no.\ of triplets, $m=4$ & &335      & 1575 & 2401  & 6999\\
            \end{tabular}
            \caption{Number of function evaluations per Building block for reconstructing the decision boundaries
                (see Section~\ref{subsec:DecisionAnalysis}).}\label{tab_resultsTOPSIS2D}
        \end{center}
    \end{table}

    Due to recent geopolitical events, car users today may attach a different importance to issues of security of supply,
    e.g.~with fuel, than in 2021 when~\cite{BALL2021120925} was written.
    That kind of considerations are subsumed in the category ``Social'' which motivates to additionally consider the
    weights of that category to be variable.
    It turns out that with stronger emphasis on the category ``Social'' car users tend to prefer a BEV
    (Fig.~\ref{fig:MCDA}, right).
    One reason for this may be that the dependence on oil imports makes the purchase of an ICE or even an HEV seem
    less attractive.

    \subsection{Surface Reconstruction in 3D from scattering}\label{subsec:surface_reconstruction_3D}
    In our tests, we will use different scatterers to apply our algorithm to such as the unit sphere, the ellipsoid,
    the peanut, the acorn, the cushion, the round short cylinder, and the round long cylinder.
    Their surfaces are given parametrically in spherical coordinates
    $x=r_1 \sin(\phi)\cos(\theta)$, $y=r_2 \sin(\phi)\cos(\theta)$, and $z=r_3\cos(\phi)$ with $\theta\in [0,2\pi)$,
    $\theta\in [0,\pi]$ as $r_1=r_2=r_3=1$ for the unit sphere, $r_1=r_2=1$ and $r_3=6/5$ for the ellipsoid,
    $r_1^2=r_2^2=r_3^2=9(\cos^2(\phi)+\sin^2(\phi)/4)/4$ for the peanut, $r_1^2=r_2^2=r_3^2=9(17/4+2\cos(3\phi))/25$
    for the acorn, $r_1=r_2=r_3=1-\cos(2\phi)/2$ for the cushion,
    $r_1^{10}=r_2^{10}=r_3^{10}=1/((2\sin(\phi)/3)^{10}+\cos^{10}(\phi))$ for the round short cylinder, and
    $r_1^{10}=r_2^{10}=r_3^{10}=1/((2\cos(\phi)/3)^{10}+\sin^{10}(\phi))$ for the round long cylinder, respectively.
    We will use $m=1026$ number of incident and observation directions for the construction of the
    far-field data with the parameters $k_e=2$, $k_i=1$, and $\tau=1/2$ as also used in~\cite[p. 18]{anachakle}.
    Therefore, the factorization algorithm appears as classification function $f$.

\begin{figure}
        \begin{center}
            \includegraphics[height=3.5cm]{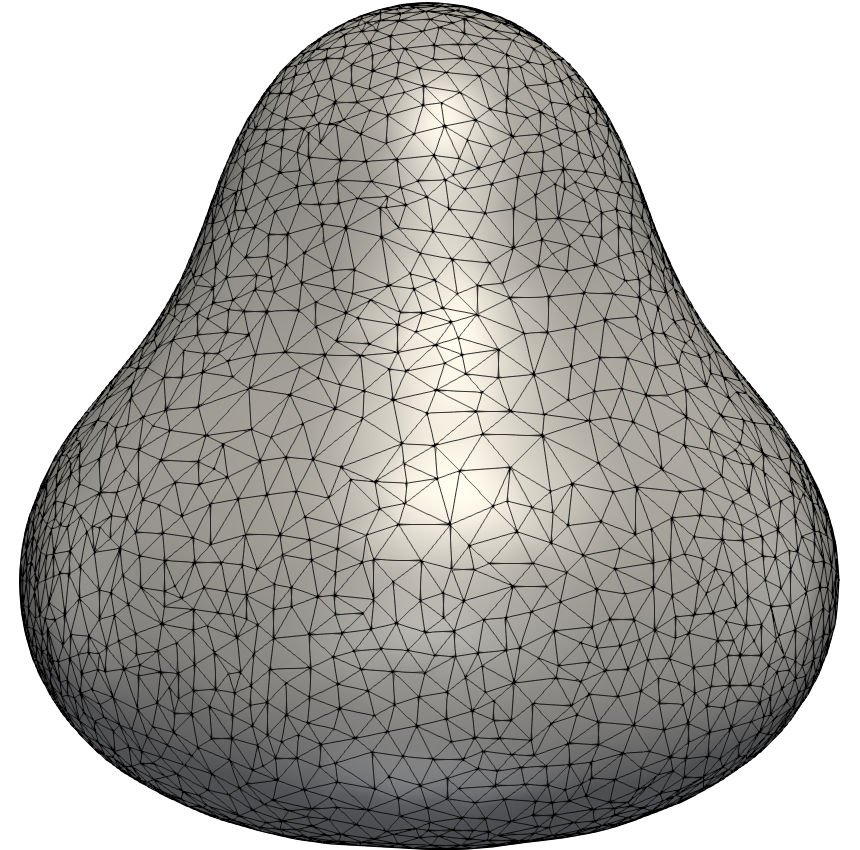}
            \includegraphics[height=3.5cm]{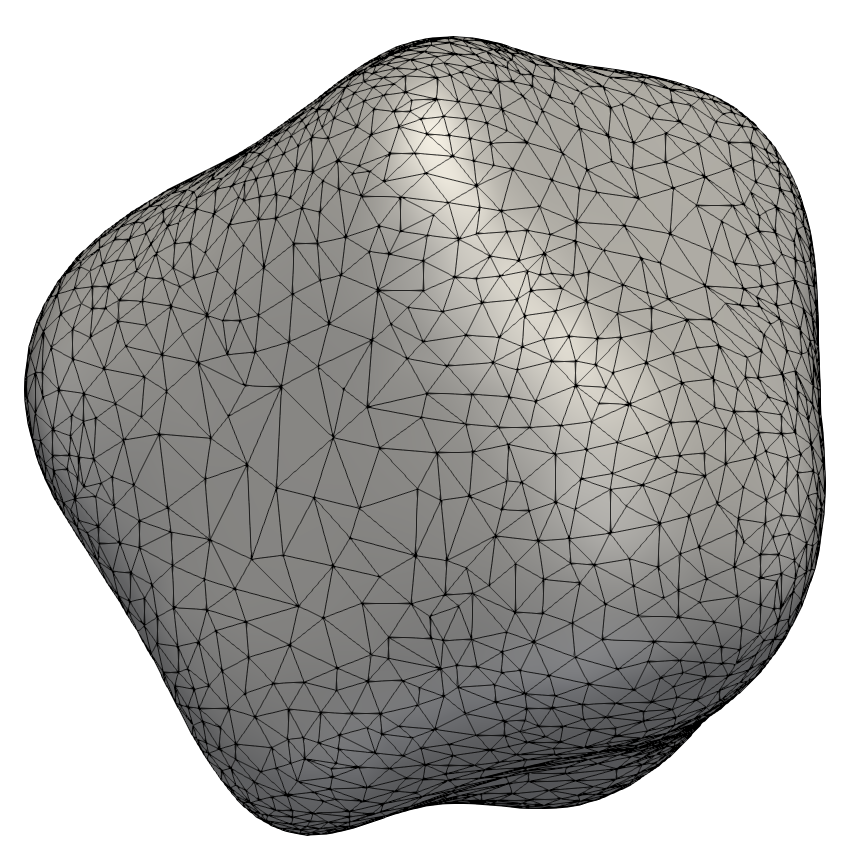}
            \includegraphics[height=3.0cm]{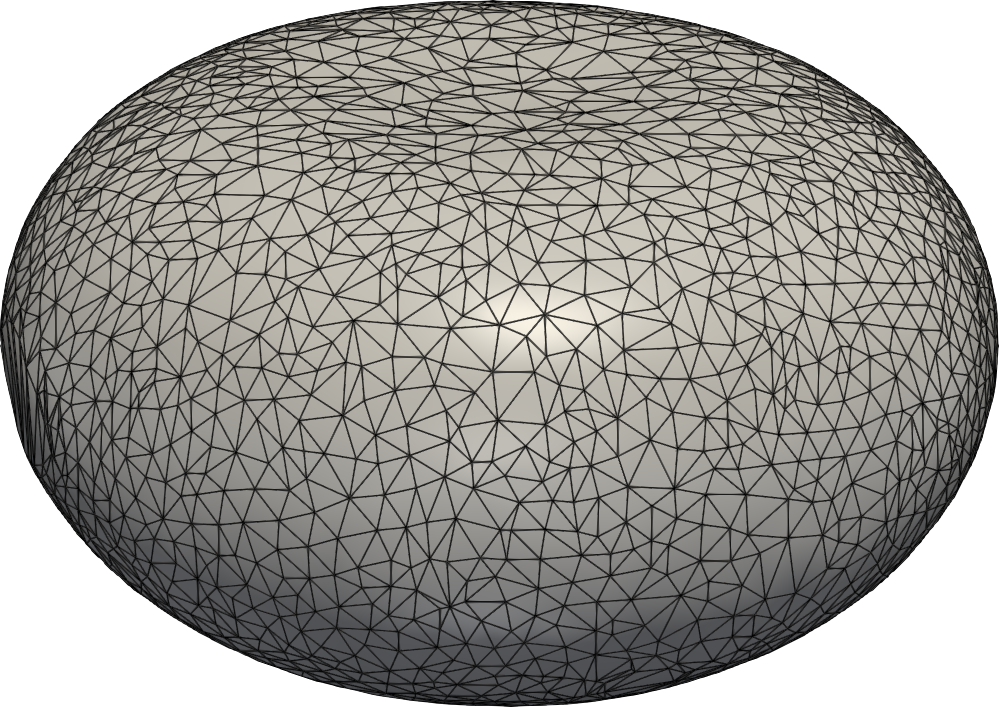}
            \includegraphics[height=3.5cm]{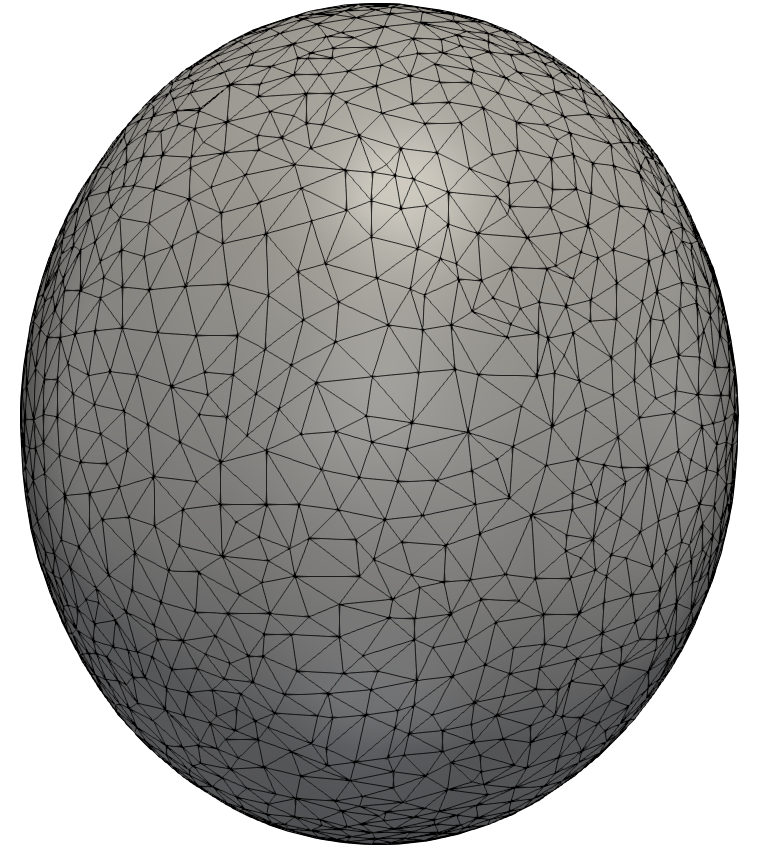}
            \caption{Reconstructed scatterers acorn, bumpy sphere, cushion and ellipsoid.}
            \label{3D-scatterers_1}
        \end{center}
    \end{figure}

\begin{figure}
        \begin{center}
            \includegraphics[height=3.5cm]{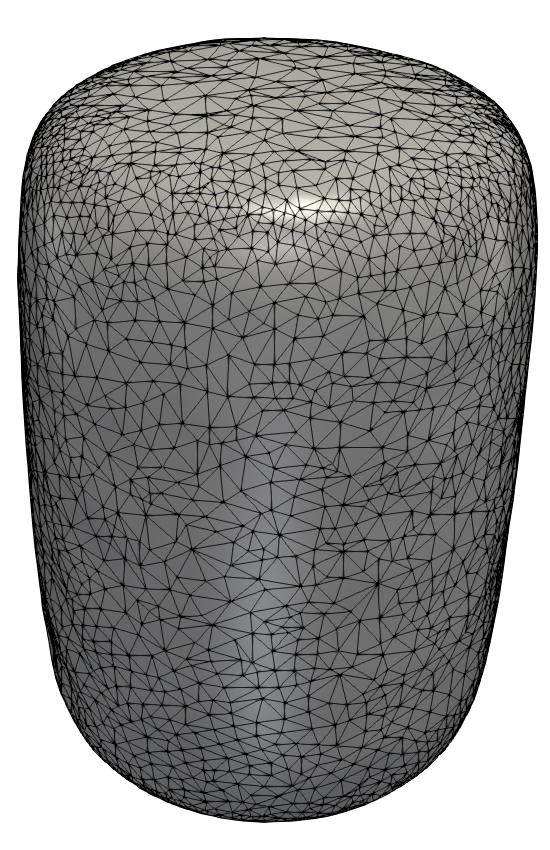}\hspace{3mm}
            \includegraphics[height=3.5cm]{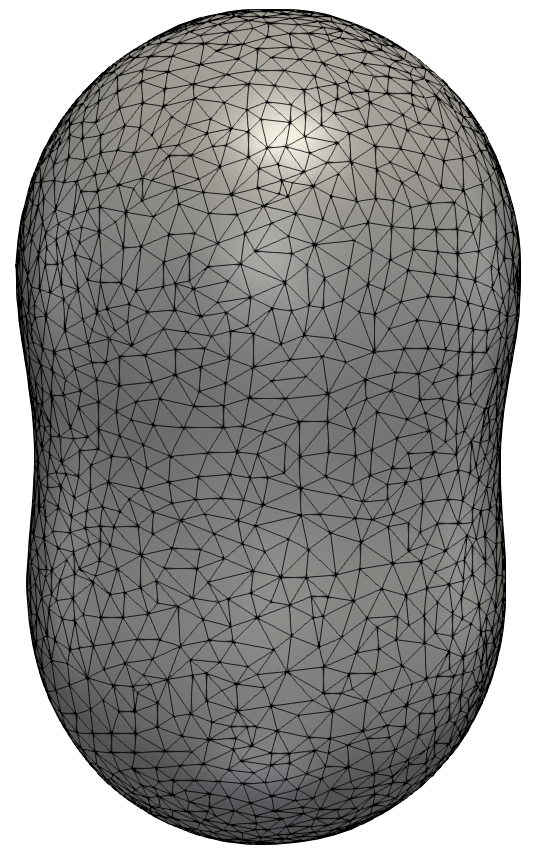}\hspace{3mm}
            \includegraphics[height=3.5cm]{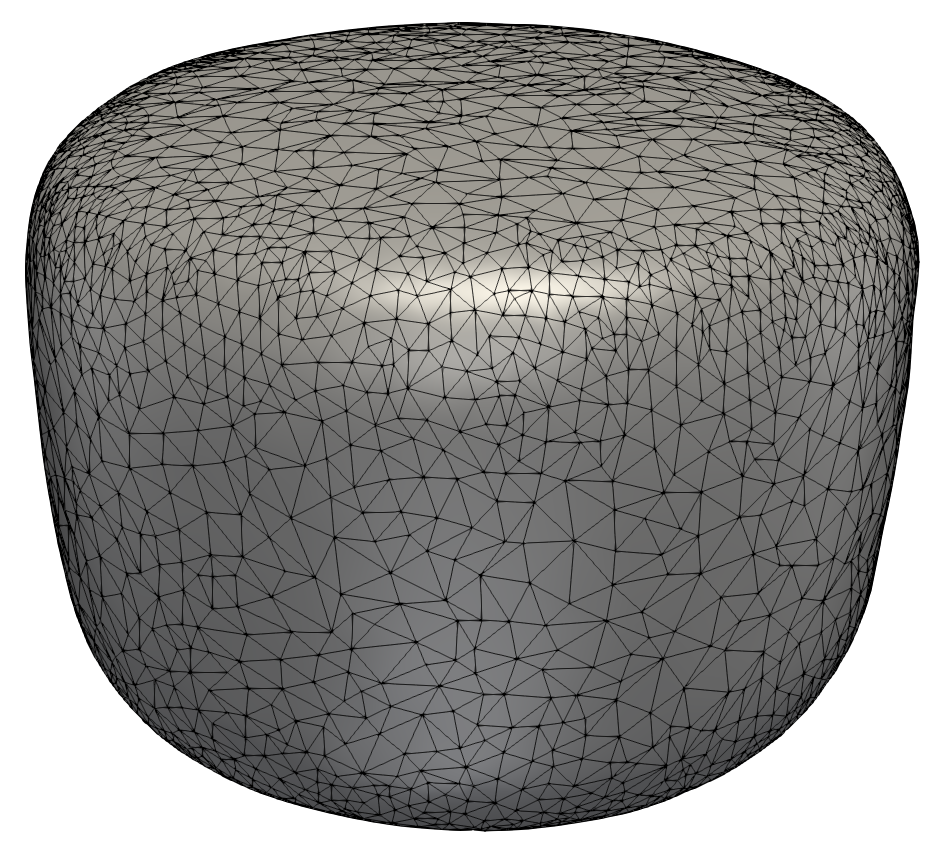}\hspace{3mm}
            \includegraphics[height=3.5cm]{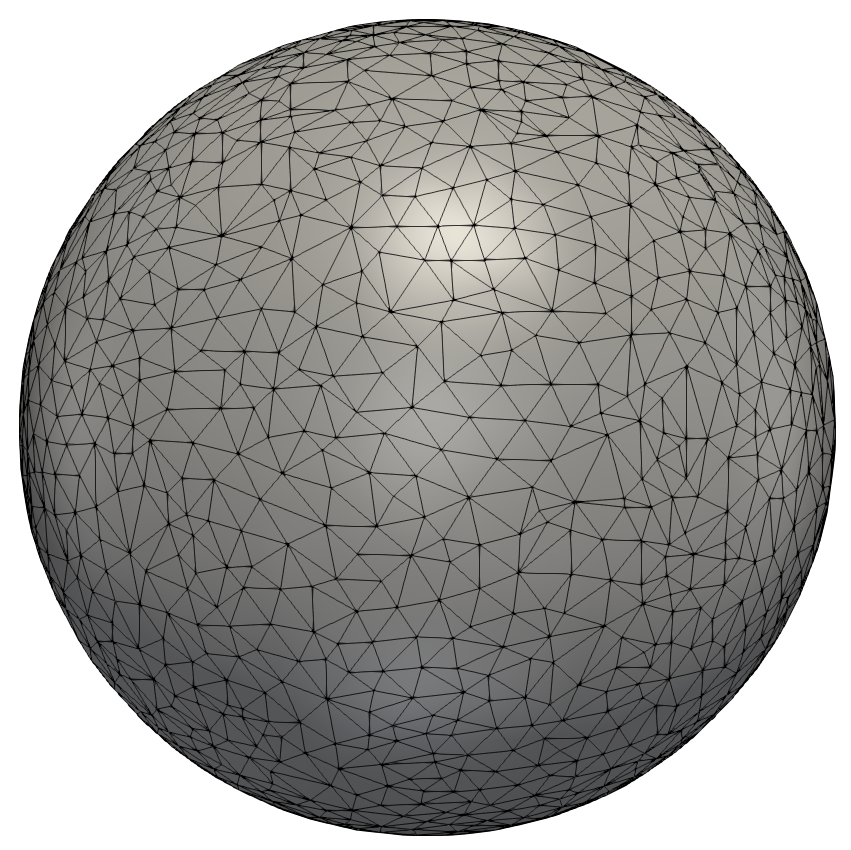}
            \caption{Reconstructed Scatterers long cylinder, peanut, short cylinder, and sphere.}
            \label{3D-scatterers_2}
        \end{center}
    \end{figure}

    For our experiments, we set all algorithm-related parameters as in Section~\ref{subsec:approximation_in_3D}
    except of $\varepsilon_{\mathrm{err}} = 0.01$ and $\varepsilon_{\text{gap}} = 0.25$.
    The initial set $X$ consists of $200$ Halton-distributed points in $\Omega = [-1.5, 1.5]^3$ apart of the long and the
    short cylinder, where we set $\Omega = [-2, 2]^3$.
    The results~(Figs.\ \ref{3D-scatterers_1} and~\ref{3D-scatterers_2}) indicate successful reconstructions.
    For the number of triplets after each part of the algorithm and the number of function evaluations, we refer to
    Tables~\ref{tab:nevals_scatterers} and~\ref{tab:ntriplets_scatterers}.
    As no triplet in these tests fulfilled condition~\eqref{eq:angle_cond}, no expansion took place, such that we
    omit the corresponding column in our tables.
    As all scatters feature one closed surface, this behaviour of our algorithm is as desired.
    In order to obtain a visually appealing reconstruction of the scatterers, Anagnostopoulos et al.~\cite{anachakle}
    create a tensor product set of $55^3 = 166,375$ points and use the corresponding classifications to compute isosurfaces
    based upon these data for visualisation purposes.
    Our approach, on the other hand, requires only a fraction of these function evaluations.
    Since about half of the total computation time was used by Kirsch's factorization method in the surface reconstruction
    with a Matlab implementation of our method, our algorithm enables a significant speedup compared to~\cite{anachakle}.

    Moreover, in contrast to the level set approach, our method provides a set of points near the scatterer, which can
    be used for computing a further refined surface representation like triangulation or higher order interpolation.
    There are two major sources of inaccuracy in surface reconstruction: the inaccuracy of the factorization method
    itself and the error induced by the representation of the reconstructed surface, be it an isosurface or a set of
    points.
    As our algorithm controls the latter one, comparing with the ground truth aka the analytical description of the true
    surface allows for analysing the former error source far more easily than having only implicit surfaces available.
    \begin{remark}
        In our tests, we do not commit to inverse crime.
        The far field data used for reconstruction later one have been produced using boundary element collocation with
        high accuracy.
        The inverse problem is solved using Kirsch's factorization method combined with the algorithm proposed in this
        work.
    \end{remark}

    \begin{table}
        \begin{center}
            \begin{tabular}{c|c|c|c|c}
                scatterer    & up to $\mathcal{M}^2$         & \texttt{iniapprox}         & \texttt{fill}     & \texttt{adapt}   \\\hline
                acorn        & 463   & 1891   & 7242   & 8272 \\
                bumpy sphere & 384   & 1326   & 4844   & 8911 \\
                cushion      & 473 & 1956 & 6237 & 7881\\
                ellipsoid    & 383     & 1329     & 4694     & 5458 \\
                long cyl.    & 378      & 1377      & 7119      & 9816 \\
                peanut       & 394       & 1381     & 5149       & 6922  \\
                short cyl.   & 465      & 2033      & 9021      & 9667 \\
                sphere       & 384       & 1314     & 3954       & 5102  \\
            \end{tabular}
            \caption{Number of function evaluations per Building block for reconstructing the scatterers.}\label{tab:nevals_scatterers}
        \end{center}
    \end{table}

    \begin{table}
        \begin{center}
            \begin{tabular}{c|c|c|c}
                scatterer    & \texttt{iniapprox}      & \texttt{fill}                       & \texttt{adapt}                      \\\hline
                acorn        & 222   & 1291   & 3839   \\
                bumpy sphere & 158   & 827   & 3103   \\
                cushion      & 232 & 1125 & 3613 \\
                ellipsoid    & 157     & 843     & 2305     \\
                long cyl.    & 149      & 1188      & 4040      \\
                peanut       & 164       & 878       & 2864       \\
                short cyl.   & 217      & 1594      & 4543      \\
                sphere       & 159     & 715       & 2120       \\
            \end{tabular}
            \caption{Total number of triplets after the respective Building block for reconstructing the scatterers.}\label{tab:ntriplets_scatterers}
        \end{center}
    \end{table}

    \section{Conclusions and Outlook}\label{sec:conclusions}
    In this article, we presented a method for approximating manifolds of discontinuity of a function $f$ in 2D and 3D
    and demonstrated successful applications of our method to Multi-criteria Decision Aid, to generic test cases
    connected with the detection of faults and to an inverse acoustic scattering problem.
    In all cases, our method requires significantly fewer evaluations of $f$ than previously existing algorithms we
    compared our method to.
    At least for the inverse acoustic scattering problem, this leads to a significant acceleration of the overall
    reconstruction, as computing the factorization method dominates the computational time even in our computations,
    where the number of such computations could be reduced by a factor of approx.~7 to 12.

    Our algorithm could be easily enhanced with classification algorithms as in~\cite{Allasia2010EfficientAA} in order
    to consider more general functions than we did, provided the classification is certain.
    This enables for tackling fault detection problems and for applying our method for interpolation of piecewise
    smooth functions assuming that $f$ is smooth on $\Omega_1, \hdots, \Omega_n$ with
    $\Omega = \overline{\Omega_1} \cup \hdots \cup \overline{\Omega_n}$, but globally discontinuous.
    In~\cite{LENARDUZZI2017113}, the authors propose interpolation with radial basis functions on each $\Omega_i$ in
    this setting.
    This however requires knowledge about the boundaries of each $\Omega_i$ which could be obtained using our method.
    For that purpose,~\cite{LENARDUZZI2017113} provides a fault detection algorithm based on local approximation
    properties.
    Combining these two approaches is subject of our current research.\\

    \noindent
        {\bf Acknowledgements:} The authors want to thank Rodin Eybesh and Luis Hasenauer for porting a significant part
    of our Matlab implementation to python.

    \section*{Declarations}
    \begin{itemize}
    \item {\bf Conflict of Interest}: The authors declare no competing interests.
    \item {\bf Code Availability:} The codes used for producing the results presented are available at
    \url{https://github.com/mgrajewski/faultapprox-matlab} and
    \url{https://github.com/mgrajewski/faultapprox-python}
    \end{itemize}
    \bibliographystyle{acm}
    \bibliography{literature}
\end{document}